\documentclass[11pt]{amsart}
\usepackage{graphics, curves, epsfig, verbatim, graphicx}
\usepackage[pass]{geometry}
\usepackage{enumerate} 
\usepackage{pst-plot, pst-grad}
\usepackage{latexsym,eucal,amsfonts,amssymb,amsmath,amscd}
\usepackage{mathrsfs}
\usepackage{tikz}
\usetikzlibrary{matrix}
\numberwithin{equation}{section}
\newtheorem{theorem}{Theorem}[section]

\newtheorem{lemma}[theorem]{Lemma}
\newtheorem{prop}[theorem]{Proposition}

\newtheorem{definition}[theorem]{Definition}

\newtheorem{remark}[theorem]{Remark}
\def \bpf {\begin{proof}}
\def \epf {\end{proof}}
\def \beq {\begin{equation*}}
\def \eeq {\end{equation*}}
\def \bsp{\begin{split}}
\def \esp{\end{split}}
\def \beqq {\begin{equation}}
\def \eeqq {\end{equation}}
\def \mca {{\mathscr A}}
\def \mcb {{\mathscr B}}

\def \mcd {{\mathscr D}}
\def \mce {{\mathscr E}}
\def \mcf {{\mathscr F}}
\def \mcg {{\mathscr G}}
\def \mch {{\mathscr H}}
\def \mci {{\mathscr I}}
\def \mcj {{\mathscr J}}
\def \mck {{\mathscr K}}
\def \mcl {{\mathscr L}}
\def \mcm {{\mathscr M}}

\def \mco {{\mathscr O}}
\def \mcp {{\mathscr P}}

\def \mcr {{\mathscr R}}
\def \mcs {{\mathscr S}}

\def \mcu {{\mathscr U}}
\def \mcv {{\mathscr V}}

\def \mcx {{\mathscr X}}
\def \mcy {{\mathscr Y}}

\def \mbr {{\mathbb R}}

\def \comp {\operatorname{comp}}

\def \div {\operatorname{div}}
\def \det {\operatorname{det}}

\def \tr {\operatorname{Tr}}
\def \ric {\textrm{Ric}}
\def \supp {\text{supp }}

\def\Id {\operatorname{Id}}
\def \id {\operatorname{Id}}
\def \WF {\operatorname{WF}}
\def \loc {\operatorname{loc}}
\def \ein{\operatorname{Ein}}

\def \diag {\operatorname{Diag}}
\def \eps {\epsilon}   
\def \la {\lambda}   
\def \La {\Lambda}

\def \lap {\Delta}
\def \p {\partial}

\def \det {\text{det}}

\def \ha {\frac{1}{2}}

\def \bfq {\textbf{Q}}

\def \ba {\begin{eqnarray*}}
\def \ea {\end{eqnarray*}}
\def \sym {\text{Sym}}

\def \hat {\widehat}

\def \bt {\textbf{t}}
\def \fnf {\frac{n}{4}}

\begin{document}
\title[Determination of space-time structures]{Determination of space-time structures from gravitational perturbations}
\author{Gunther Uhlmann} 
\address{Gunther Uhlmann
\newline
\indent Department of Mathematics, University of Washington 
\newline
\indent and Institute for Advanced Study, the Hong Kong University of Science and Technology} 
\email{gunther@math.washington.edu}

\author{Yiran Wang}
\address{Yiran Wang
\newline
\indent Department of Mathematics, University of Washington 
}
\email{wangy257@math.washington.edu}
\begin{abstract} 
We study inverse problems for the Einstein equations with source fields in a general form. Under a microlocal linearization stability condition, we show that by generating small gravitational perturbations and measuring the responses near a freely falling observer, one can uniquely determine the background Lorentzian metric up to isometries in a region where the gravitational perturbations can travel to and return. We apply the result to two concrete examples when the source fields are scalar fields (i.e.\ Einstein-scalar field equations) and electromagnetic fields (i.e.\ Einstein-Maxwell equations).
\end{abstract}

\maketitle
\tableofcontents

\section{Introduction}
Consider Einstein equations on a smooth four dimensional manifold $M$:
\beq
\ein(g) = T_{sour},
\eeq
where $\ein(g)$ denotes the Einstein tensor  and $T_{sour}$ is the source stress-energy tensor describing the distribution of source fields, e.g.\ scalar fields, electromagnetic fields, perfect fluids etc. The solution $g$ describes the space-time structure, or equivalently the  gravitational field. Let $\hat g$ be a background field with background stress-energy tensor $\hat T$ such that $\ein(\hat g) = \hat T$. The problem we study in this work is to determine the background metric $\hat g$ by making small perturbations of $\hat T$ and  measuring the perturbations of $\hat g$. In particular, we consider the Einstein equations with stress-energy tensors of an abstract form. We formulate a microlocal linearization stability condition, under which we prove that the topological, differentiable structure and the metric (up to an isometry) of space time can be determined by making active measurement near the world line of an observer. See Theorem 2.1 for the precise statements. Then we apply the method to two concrete examples of stress-energy tensors: the scalar fields (i.e.\ Einstein-scalar field equations) and the electromagnetic fields (i.e.\ Einstein-Maxwell equations).  Comparing with previous work, we obtain stronger results (i.e.\ determination of the metric up to isometries  rather than conformal classes) in a  general setup. Also, we take a thorough microlocal approach, which makes the arguments more transparent.

The inverse problem for Einstein equations was proposed and studied in Kurylev-Lassas-Uhlmann \cite{KLU1} for the Einstein-scalar field equations. The observation of gravitational waves done first by the LIGO project \cite{Ab} has opened up the possibility of many more observations of space-times.  For Einstein-Maxwell equations, the problem is studied in Lassas-Uhlmann-Wang \cite{LUW1}.  The authors considered the Einstein-Maxwell system in vacuum and used the measurements of the response to electromagnetic perturbations. In this article we consider a general stress-energy tensor  but we consider the response to both gravitational and electromagnetic perturbations. Similar inverse problems have been studied for semilinear wave equations with quadratic nonlinearity in Kurylev-Lassas-Uhlmann \cite{KLU}, with general nonlinear term in Lassas-Uhlmann-Wang \cite{LUW} and with quadratic derivative terms in Wang-Zhou \cite{WZ}. In particular in \cite{LUW}, the authors improved the previous results in \cite{KLU} so that the isometry class of the metric and some information about the nonlinear terms can be determined in many cases. In de Hoop-Uhlmann-Wang \cite{DUW}, the nonlinear responses of two scalar waves at an interface of different media was considered and the related inverse problem was addressed. Also, in de Hoop-Uhlmann-Wang \cite{DUW1}, the problem for elastic systems is considered. In particular, the nonlinear interaction of two elastic waves is carefully analyzed and used to determine elastic parameters from boundary measurements.

\section{The main results}\label{sec-main}
We introduce notations and assumptions to be used throughout the paper unless otherwise specified. Let $M$ be a $4$-dimensional smooth manifold and $g$ a time oriented Lorentzian metric on $M$. We take the signature of  $g$ as $(-, +, +, +)$. For $p\in M$, the set of light-like vectors at $p$ is denoted by $L_{g, p}M = \{\theta \in T_pM\backslash \{0\}:  g(\theta, \theta) = 0\}$ and the corresponding bundle by $L_gM = \bigcup_{p\in M}L_{g, p}M$. The set of  future (past) light-like vectors are denoted by $L^+_{g, p}M$ ($L^-_{g, p}M$), and the corresponding bundle $L_g^\pm M = \bigcup_{p\in M}L_{g, p}^\pm M$. Similarly, the set of light-like co-vectors at $p\in M$ is denoted by $L^*_{g, p}M$ and the bundles $L_g^{*}M, L_g^{*,\pm}M$. Since the metric $g$ is non-degenerate, it induces a natural isomorphism $i_g: T_pM\rightarrow T_p^*M$ for any $p \in M$.  {\em With this isomorphism, we use vectors and co-vectors interchangeably.} When it becomes necessary, we use the standard raising and lowering index operation. More explicitly,
\beq
\xi^\sharp = i_g(\xi)\in T_p^*M, \ \ \xi \in T_p M; \ \ \eta^\flat = i_g^{-1}(\eta)\in T_p M,\ \ \eta \in T_p^*M.
\eeq
We emphasize that these operators depend on the metric $g$ as well, though the dependency is not showing up in the notations $\sharp, \flat$. 

Consider causal relations on $(M, g)$. For $p, q\in M$, we denote by $p\ll q$ ($p< q$) if $p \neq q$ and there is a future pointing time-like (causal) curve from $p$ to $q$. We denote by $p\leq q$ if $p = q$ or $p<q$. The chronological future of $p\in M$ is  $I_g^+(p) = \{q\in M: p\ll q\}$ and the causal future of $p\in M$ is $J_g^+(p) = \{q\in M:  p \leq q\}$. The chronological past and causal past are defined by  $I_g^-(p) = \{q\in M: q\ll p\}, J_g^+(p) = \{q\in M:  q \leq p\}$ respectively. For any set $A\subset M$, we denote the causal future, past by $J_g^\bullet(A) = \bigcup_{p\in A}J_g^\bullet(p)$ with $\bullet = +, -$ respectively. Finally, we denote 
\beq
J_g(p, q) = J_g^+(p)\cap J_g^-(q) \text{ and } I_g(p, q) = I_g^+(p)\cap I_g^-(q).
\eeq
We remark that all of these sets depend on   the metric $g$. Let $\exp_{g, p}: T_pM  \rightarrow M$ be the exponential map on $(M, g)$. The geodesic from $p$ with  direction $\theta$ is denoted by $\gamma_{g, \theta}: \mbr \rightarrow M, \gamma_{g, \theta}(t) = \exp_{g, p}(t\theta)$. Then we denote the forward light-cone at $p\in M$ by
\beq
\mcl^+_{g, p} = \{\gamma_{g, \theta}(t): \theta \in L^+_{g, p}M, \ \ t> 0\}, 
\eeq 
which is a submanifold of $M$.

Throughout the paper, we use the standard Einstein summation convention unless otherwise specified.  For  $(M, g)$, we denote by $\ric(g)$ the Ricci curvature tensor. Let $\sym^2$ be the vector bundle of symmetric covariant $2$ tensors on $M$. For any $T\in \sym^2$, the trace of $T$ with respect to $g$ is   $\tr_g(T) = g^{\alpha\beta}T_{\alpha\beta}$. The scalar curvature $S(g)=\tr_g(\ric(g))$. The Einstein tensor is 
\beq
\ein(g) = \ric(g) - \ha S(g) g. 
\eeq
For $T\in \sym^2$, the Einstein field equations are $\ein(g) = T.$ It follows from the Bianchi identity that $T$ must satisfy the conservation law $\div_g T = 0$, where $\div_g$ denotes the divergence operator.  \\

We start with the assumptions on the background manifold $(M, \hat g)$. We assume that $\hat g$ is {\em globally hyperbolic}. According to \cite{BS}, this means that there is no closed causal paths in $(M, \hat g)$ and for any $p, q\in M$ and $p<q$, the set $J_{\hat g}(p, q)$ is compact.  In \cite{BS0}, it is proved that such $(M, \hat g)$ is isometric to the product manifold
\beq
\mbr\times \mcm \text{ with metric } \tilde g =-\beta(t, y)dt^2 + \kappa(t, y),
\eeq
where $\mcm$ is a $3$-dimensional smooth manifold, $\beta: \mbr\times \mcm\rightarrow \mbr_+$ is smooth and $\kappa$ is a smooth family of Riemannian metric on $\mcm$. From now on, we shall identify $(M, \hat g)$ with this isometric image. We use $x = (t, y) = (x^0, x^1, x^2, x^3)$ as the local coordinates on $M$ where $y = (x^1, x^2, x^3)$ is the local coordinates on $\mcm$. For $\bt\in \mbr$, we denote 
\beq
M(\bt) = (-\infty, \bt)\times \mcm.
\eeq
Also, we assume that $\hat g$ is time oriented and 
\begin{gather}\label{as1}
 \text{$(M, \hat g)$ has no conjugate point.} \tag{A1}
\end{gather}
This assumption is made to avoid some technicalities and to make it easier to follow the proofs in this article. This assumption can be removed by following the arguments in \cite[Section 5]{KLU1}.

Next, we want to assume that $\hat g$ satisfies the Einstein equation 
\beq
\ein(\hat g) = \hat T, \ \ \hat T\in C^\infty(M; \sym^2)
\eeq
for some background field $\hat T$.  In general, a stress-energy tensor $T$ is a function of the physical fields $\psi= (\psi_l)_{l = 1}^L, L > 0$ in presence and the metric $g$ i.e.\ $T = T(g, \psi)$. 
 Let $\bold{B}^L, L\geq 0$ be a vector bundle on $M$ with fiber dimension $L \geq 1.$ We assume that $\psi$ are governed by the wave equation 
\beq
\square_g \psi_l  + V_l(\psi) = \mcf^\Psi_l, \quad l = 1, 2, \cdots, L, 
\eeq
where $V_l$ are smooth potential functions and $\mcf^\Psi$ is the source term. We shall assume that 
\begin{gather}
 \hat T =  \hat T(g, \psi) \text{ is a smooth function of $g \in \sym^2$ and $\psi \in \bold{B}^L$ } \tag{A2}\ 
\end{gather}
so the stress energy tensor is essentially described by the field $\psi$ or equivalently the source $\mcf^\Psi$ through the wave equation. We remark that $\hat T$ only depends on $g$ and not on the derivatives $\p g$.  Let $\hat \psi$ be a background field  and we assume that 
\begin{gather}
\text{$\hat g$ is a solution to } \ein(g) = \hat T(g, \hat \psi) \text{ for some } \hat \psi \in C^\infty(M; \bold{B}^L)\\
  \text{satisfying} \quad \square_{\hat g} \hat \psi_l  + V_l(\hat \psi) = 0 \text{ for } t > 0 . \tag{A3}
\end{gather}
Here, we regard $\hat \psi$ as the background field that produces $\hat g$, but sometimes we also call $\hat T = \hat T(\hat g, \hat \psi)$ the background stress-energy tensor. Assumptions (A1)-(A3) are the assumptions for the unperturbed background fields.\\

Now we consider the Einstein equations near a freely falling observer, which is represented by a time-like geodesic on a Lorentzian manifold. Let $\mu_{\hat g}(t), t\in [-1, 1]$ be a time-like geodesic (segment) on $(M, \hat g)$. Without loss of generality, we assume $\mu_{\hat g}(-1) \in \{0\}\times \mcm$. Let $V$ be an open relatively compact neighborhood of $\mu_{\hat g}([-1, 1])$ and  $V\subset M(\bt_0)$ for some $\bt_0>0$. Let $\mcf \in C^4(M; \sym^2), \mcf^\Psi \in C^4(M; \bold{B}^L)$ be compactly supported in $V$ and sufficiently small. We consider the solution $g$ to 
\beqq\label{einsour00}
\begin{gathered}
\ein(g) =  \hat T(g, \psi) + \mcf, \text{ in } M(\bt_0),\\
\square_g \psi_l + V_l(\psi) = \mcf^\Psi_l, \quad l = 1, 2, \cdots, L, \\
g = \hat g, \quad \psi = \hat \psi \quad \text{ in } M(\bt_0)\backslash J_g^+(\supp (\overline \mcf)),
\end{gathered}
\eeqq
where $\overline \mcf = (\mcf, \mcf^\Psi)$. We shall denote by $T_{sour} = \hat T(g, \psi) + \mcf$  the source stress energy tensor. \eqref{einsour00} is the full version of the forward problem or coupled Einstein-field equations and we establish the local well-posedness in Section \ref{sec-redein}. There is a complexity coming from the fact that the source $T_{sour}$ is subject to the conservation law $\div_g T_{sour} = 0$ as a consequence of Bianchi's identity.  Actually, the background system   is \eqref{einsour00} coupled with the conservation law. We address the difficulty later using the concept of {\em microlocal linearization stability} in Section \ref{sec-muls}.  In our formulation, we regard $\mcf$ as the active source perturbation that one can control and that generates the gravitational waves, and $\psi$ is some adaptive field so that the conservation law holds. 

For our treatment of the inverse problem, sometimes it is convenient to regard $\psi, \mcf$ as the source fields. Once we establish the well-posedness of \eqref{einsour0}, we can ignore the field equations for $\psi_l$ in \eqref{einsour0} and just consider the Einstein equations 
\beqq\label{einsour0}
\begin{gathered}
\ein(g) =  \hat T(g, \psi) + \mcf, \text{ in } M(\bt_0),\\
g = \hat g,  \text{ in } M(\bt_0)\backslash J_g^+(\supp (\mcf)).
\end{gathered}
\eeqq
We define the source-to-solution set for the Einstein equation as 
\beqq\label{eqdata}
\begin{split}
\mcd(\delta; V) = \{(g, \psi, \mcf):   &\text{ $\mcf \in C^4(V; \sym^2)$ is compactly supported, } \|\mcf\|_{C^4(M; \sym^2)} < \delta; \\
&\psi \in C^4(M; \bold{B}^L), \ \ \|\psi - \hat \psi\|_{C^4(M(\bt_0))} < \delta;\\
 & \text{$g \in C^4(M(\bt_0); \sym^2)$ is a solution of \eqref{einsour0} with $\psi, \mcf$ } \}.
\end{split}
\eeqq
Of course, one can equivalently describe the set using the source $\mcf^\Psi \in C^4(M; \bold{B}^L)$. Here, to define semi-norms for $C^m(M)$,  we shall take a complete Riemannian metric $\widehat g^+$ on $M$, whose existence is guaranteed by \cite{NO}. Then we can introduce distances on $M$ and $TM$, and Sobolev spaces on $M$ using $\widehat g^+$. The  regularity requirements will be explained in Section \ref{sec-redein}. We remark that the data set also depends on $\hat g, \hat \psi, \bt_0$ which does not show up in the notation. Later, we use the abbreviated notation $\mcd(\delta)$ when there is no confusion. 

\begin{figure}[htbp]
\centering
\includegraphics[scale=0.75]{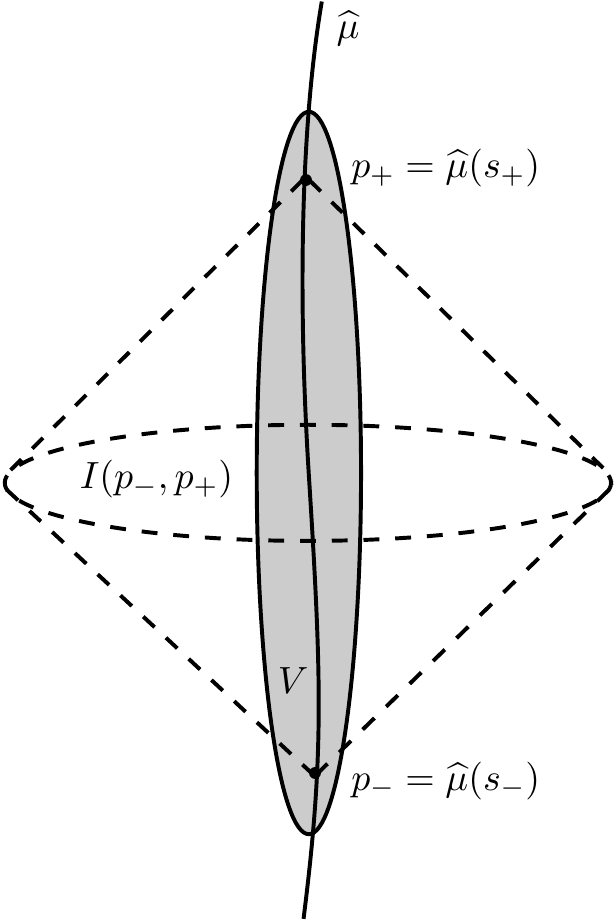}
\caption{Illustration of the inverse problem}
\label{figinv}
\end{figure}

Let's recall that solutions to the Einstein equations are defined uniquely only as equivalent classes under diffeomorphisms (of certain regularity). Hence in \eqref{eqdata}, the data $(g, \psi, \mcf)$ should be regarded as the equivalent class under such diffeomorphisms. In any fixed coordinate system or gauge, we obtain a unique representation of the data, so the formulation of the data set should depend on the gauge where the measurements are taken. In this work, we consider Fermi coordinates associated with freely falling observers. The exact formulation and physical explanations of this coordinates is postponed to Section \ref{sec-fermi}. Roughly speaking, this means that we take the data set $\mcd(\delta; V)$ in Fermi coordinates and we denote it by $\mcd_{F}(\delta; V)$, see \eqref{dataF}.  

The inverse problem we study is follows. Let $p_\pm =  \mu_{\hat g}(s_\pm),  -1<s_-<s_+<1$ be two points on the time-like geodesic. Suppose that we can control all possible perturbations $\mcf$ in $V$ and measure the metric perturbations $g$ in $V$, can we determine $\hat g$ in $I(p_-, p_+)$ from these information? See Figure \ref{figinv}.  Our main result is that these data sets determine the metric $g$ uniquely up to isometries in $I(p_-, p_+)$. More precisely, 

\begin{theorem}\label{main}
Let $(M, \hat g^{(i)}), i = 1, 2$ be two four-dimensional time oriented, globally hyperbolic smooth Lorentzian manifolds and $\hat T^{(i)} = \hat T^{(i)}(g^{(i)}, \psi^{(i)}), \psi^{(i)}\in \bold{B}^L$ such that (A1)-(A3) hold. 
Let $\mu_{\hat g^{(i)}}(t), t\in [-1, 1]$ be time-like geodesics on $(M, \hat g^{(i)})$. Let $V$ be open relatively compact neighborhoods of both $\mu_{\hat g^{(i)}}([-1, 1])$ and  $V\subset M^{(i)}(\bt_0)$ for some $\bt_0>0$.  Let $-1<s_-<s_+ < 1$ and $p_\pm^{(i)} = \mu_{\hat g^{(i)}}(s_\pm)$ such that $p^{(1)}_\pm = p^{(2)}_\pm$. 
For  $\delta > 0$, we consider solutions $g^{(i)}$ to the Einstein equations
\beqq\label{einm}
\begin{gathered}
\ein(g^{(i)}) = \hat T^{(i)}(g^{(i)}, \psi^{(i)}) + \mcf^{(i)}, \text{ in } M^{(i)}(\bt_0),\\
g^{(i)} = \hat g^{(i)}, \text{ in } M^{(i)}(\bt_0)\backslash J^+_{g^{(i)}}(\supp (\mcf^{(i)})),
\end{gathered}
\eeqq
where $(\psi^{(i)}, \mcf^{(i)}) \in \mcd^{(i)}(\delta; V)$ defined as \eqref{eqdata}. Suppose that 
 \begin{enumerate}
  \item  the microlocal linearization stability condition (Definition \ref{mlstab}) holds for $\mcd^{(i)}(\delta; V)$.
 \item the source-to-solution maps satisfy 
\[
 \mcd^{(1)}_F(\delta; V) = \mcd^{(2)}_F(\delta; V).
\]
\end{enumerate}
Then there is a diffeomorphism $\Psi: I_{\hat g^{(1)}}(p^{(1)}_-, p^{(1)}_+)\rightarrow I_{\hat g^{(2)}}(p^{(2)}_-, p^{(2)}_+)$ such that  $\Psi^*\hat g^{(2)} = \hat g^{(1)}$ in $I_{\hat g^{(1)}}(p^{(1)}_-, p^{(1)}_+)$. 
\end{theorem}

In Section \ref{sec-wgauge}, we prove the theorem for the data sets formulated in wave gauge. In Section \ref{sec-einscal}, we prove the theorem (see Theorem \ref{mainsca}) when $\hat T$ is given by scalar fields,  improving the result of \cite{KLU1}. In Section \ref{sec-einmax}, we prove the theorem (see Theorem \ref{mainem}) when $\hat T$ is given by electromagnetic fields. It is instructive to compare this result with the one in \cite{LUW1} where only electromagnetic source perturbations are allowed. \\

The way we formulate and solve the problem consists of three parts, according to which we organized the contents of the paper. 

{\em(1) Microlocal linearization stability: Section \ref{sec-redein} and \ref{sec-muls}}. This part essentially concerns the forward problem, especially the microlocal behavior of the solutions to the Einstein equations. The way we solve the inverse problem  is to make use of the nonlinear interaction of progressing gravitational  waves. In particular, we consider progressing waves with conormal singularities (called distorted plane waves in Section \ref{sec-muls}) and study the nonlinear interaction of such waves. So it is desirable that the data set $\mcd(\delta)$ contains sources which can generate such waves. Because the source must satisfy the conservation law $\div_g T_{sour} = 0$, it is appropriate to consider the coupled system  
\beqq\label{einsour1}
\begin{gathered}
\ein(g) = T_{sour}, \\
\div_g T_{sour} = 0,
\end{gathered}
\eeqq
instead of \eqref{einsour0}. In this work, we do not pursue this point but formulate a {\em microlocal linearization stability condition}  in Section \ref{sec-muls} after we introduce the linearized Einstein equations and distorted plane waves. The condition holds for the Einstein-scalar fields under certain conditions as shown in \cite{KLU1}. One possible approach to verify the condition is to study the coupled system \eqref{einsour1}, especially the propagation of singularities. 

 {\em (2) Interaction of singularities:  Section \ref{sec-wgauge}}. We study the interaction of four progressive gravitational waves and we show that their nonlinear interaction generates a point source of gravitational waves.   All the analysis are carried out in wave gauge, because thanks to the work of Choquet-Bruhat, we know that the reduced Einstein equation is a second order quasilinear hyperbolic system which is convenient for analyzing the propagation of singularities. In Section \ref{sec-redein}, we review the reduced Einstein equations in wave gauge, especially the solvability and stability. In Section \ref{sec-wgauge}, we analyze the interactions of four distorted plane waves using methods developed in \cite{LUW} and \cite{LUW1} to show that a new point source is generated. We determine the leading singularities and find the  principal symbols of the new wave.  We make specific choices of the sources so that the new singularity is not vanishing and can be detected. We remark that the subject of propagation and nonlinear interaction of gravitational waves is interesting in its own right in General Relativity theory, see \cite{Gri} and more recently \cite{LuR, LuR1}. 

 {\em (3) The observation sets:  Section \ref{sec-wave} and \ref{sec-fermi}}. As we already mentioned, there is a problem of choosing observation gauges for the gravitational waves produced in part (2). From a mathematical point of view, the most convenient gauge is the wave gauge. In Section \ref{sec-wave}, we follow the methods in \cite{KLU1} to prove the determination of the conformal class of the metric.  To determine the conformal factor, we analyze the conformal transformation of the principal symbols of the newly generated waves as in \cite{LUW} and prove Theorem \ref{main} when the measurements are in wave gauge. For other gauges e.g.\ the Fermi coordinates associated with freely falling observers, we need an additional step to show that such singularities are observable after the gauge transformation.  This is done in Section \ref{sec-fermi}, where we complete the proof of Theorem \ref{main}.

\section{Local well-posedness of the Einstein equations}\label{sec-redein}
We discuss the well-posedness result needed for our analysis. As this is a well-studied subject, see for example \cite{Cb}, we shall be very brief. It is convenient to use an equivalent form of the Einstein equations. By taking the trace of $\ein(g) = T_{sour}$, we have $S(g) = -\tr_g(T_{sour})$. The Einstein equations are equivalent to 
\beqq\label{ein1}
\ric(g) = \rho(g, T_{sour}), \ \ \rho(g, T) \doteq T - \ha (\tr_g(T)) g.
\eeqq
The Laplace-Beltrami operator is given in local coordinates by
\beqq\label{eqbeltra}
\square_g = -|\det g|^{-\ha}  \frac{\p} {\p x^\beta} (|\det g|^\ha g^{\beta\alpha}  \frac{\p} {\p x^\alpha}) = -g^{\alpha\beta}  \frac{\p^2} {\p x^\alpha\p x^\beta} + \Gamma^\alpha  \frac{\p} {\p x^\alpha},
\eeqq
where $\Gamma^\alpha$ is the contracted Christoffel symbol given by
\beq
\Gamma^\mu = g^{\alpha\beta}\Gamma_{\alpha\beta}^\mu = -|\det g|^{-\ha} \frac{\p} {\p x^\la}(|\det g|^\ha g^{\la\mu}),
\eeq
see e.g. Section 3 of \cite{FM} or \cite{KLU1}. Here, the Christoffel symbols are given by
\beq
\Gamma_{\alpha\beta}^\mu = \ha g^{\mu\la}(\frac{\p g_{\la\alpha}}{\p x^\beta} + \frac{\p g_{\la\beta}}{\p x^\alpha}- \frac{\p g_{\alpha\beta}}{\p x^\la}).
\eeq

We start with the well-posedness of the Einstein equations in wave gauge where they are reduced to a quasilinear hyperbolic system.  Let $g'$ be a solution to the Einstein equations \eqref{ein1} and $g'$ be close to $\hat g$ in $C^m(M)$ with $m$ to be specified later. Consider the wave map $f: (M, g')\rightarrow (M, \hat g)$. We refer the reader to e.g.\ \cite[Appendix III]{Cb} for a general discussion of wave maps and \cite[Appendix A.3]{KLU1} for the results we need below. For $m\geq 5$ an integer, let $g'\in C^m(M(\bt_0')), \bt_0'>0$ be sufficiently close to $\hat g$ in $C^m(M(\bt_0'))$. Then there is a unique wave map 
\beq
f\in C^0([0, \bt_0']; H^{m-1}(\mcm))\cap C^1([0, \bt_0']; H^{m-2}(\mcm)). 
\eeq
Moreover, introduce  
\beq
E^m(\bt; \bullet) = \bigcap_{p = 0}^{m-1} C^p([0, \bt]; H^{m-1-p}(\mcm; \bullet)),
\eeq
where $\bullet$ denotes a vector bundle on $M$. When there is no confusion, we abbreviate the notation as $E^m(\bt), E^m(\bullet)$ or even $E^m.$ For $m\geq 4$ and $m$ an even integer, the wave map
\beq
f\in E^m(\bt_0') \subset C^{m-3}(M(\bt_0'))
\eeq
and $f$ depends on $g' \in C^m(M(\bt_0'))$ continuously. Here, we will choose $\bt_0'<\bt_0$ so that $f(M(\bt_0'))\subset M(\bt_0)$. 

If $f$ is a wave map with respect to $(g', \hat g)$ and set $g=f_*g'$, then the identity map $\id$ is a wave map with respect to $(g, \hat g)$ and the wave map equation for $\id$ is equivalent to the harmonicity condition $\Gamma^n = \hat \Gamma^n,$ where  $\hat \Gamma^n = g^{\alpha\beta}\hat \Gamma^n_{\alpha\beta}$. Let $F^n = \Gamma^n - \hat \Gamma^n$. The Ricci curvature can be written as (see \cite[equation (23)]{KLU1})
\beqq\label{eqredric}
\ric_{\mu\nu}(g) = (\ric_{\hat g}(g))_{\mu\nu} + \ha (g_{\mu n} \hat \nabla_\nu F^n + g_{\nu n}\hat \nabla_\mu F^n),
\eeqq
where $\ric_{\hat g}(g)$ denotes the reduced Ricci tensor given by (see \cite[equation (23) and (22)]{KLU1})
\beqq\label{ricre}
\begin{gathered}
(\ric_{\hat g}(g))_{\mu\nu}  = -\ha g^{pq} \frac{\p^2 g_{\mu\nu}}{\p x^p \p x^q} + Q_{\mu\nu}, \text{ where }\\
Q_{\mu\nu } = g^{ab}g_{ps}\Gamma^p_{\mu b}\Gamma^s_{\nu a} + \ha (\frac{\p g_{\mu\nu}}{\p x^a}\hat \Gamma^a + g_{\nu l} \Gamma^l_{ab} g^{aq}g^{bd} \frac{\p g_{qd}}{\p x^\mu} + g_{\mu l} \Gamma^l_{ab} g^{aq} g^{bd} \frac{\p g_{qd}}{\p x^\nu}) \\
+ \ha (g_{\mu q}\frac{\p \hat\Gamma^q}{\p x^\nu}  + g_{\nu q}\frac{\p \hat \Gamma^q}{\p x^\mu}).
 \end{gathered}
\eeqq 
So in wave gauge where $F^n=0$, the Ricci tensor becomes the reduced Ricci tensor. 
Finally, the Einstein equation \eqref{einsour0}, which is $\ric(g') = \rho(g', T'_{sour})$ in $M(\bt_0')$ using the current notations,  becomes the reduced Einstein equation in wave gauge. Thus the coupled system \eqref{einsour00} in wave gauge becomes 
\beqq\label{ein2}
\begin{gathered}
\ric_{\hat g}(g) = \rho(g, T_{sour}) \text{ in } M(\bt_0),\\
\square_g \psi_l + V_l(\psi) = \mcf^\Psi, \quad l = 1, 2, \cdots, L, \\
g = \hat g, \quad \psi = \hat \psi \text{ in } M(\bt_0)\backslash J_g^+(\supp(\overline \mcf)),
\end{gathered}
\eeqq
where $g = f_*g', T_{sour} = f_*T'_{sour}$, see \cite[Appendix A.4]{KLU1} and $\psi, \mcf^\Psi$ are also regarded as the push forward in wave gauge. Notice that \eqref{ein2} is a second order quasilinear hyperbolic system for $g$. 
Let $u = g - \hat g, \phi = \psi - \hat \psi$ be the perturbed fields. The equations for $u, \phi$ are
\beqq\label{pernon}
\begin{gathered}
-g^{pq}\frac{\p^2 u }{\p x^p x^q} + \mca(x, u, \p u) =   2\rho(g, \mcf) + 2\rho(g, \hat T(g, \phi + \hat \psi)) - 2\rho(\hat g, \hat T(\hat g, \hat \psi)), \text{ in } M(\bt_0), \\
\square_g \phi_l + \mca_l^\Psi(x, g, \phi) = \mcf^\Psi, \quad l = 1, 2, \cdots, L, \\
u =\phi = 0, \ \ \text{ in } M(\bt_0)\backslash J_g^+(\supp(\overline \mcf)),
\end{gathered}
\eeqq
where $\mca, \mca^\Psi$ are smooth functions and with values in $\sym^2, \bold{B}^L$. One can see that the term $\mca^\Psi$ does not depend on  $\p g$ by using the harmonicity conditions in wave gauge. This is a second order quasilinear hyperbolic system. Applying the result in \cite[Appendix B]{KLU1}, we obtain
\begin{prop}\label{stabest}
Let $m_0\geq 4$ be an even integer, $\bt_0>0$ and $\hat \psi\in C^\infty(M(\bt_0); \bold{B}^L)$. Assume that $\mcf \in E^{m_0}(\bt_0; \sym^2)$ and $\mcf^\Psi \in E^{m_0}(\bt_0; \bold{B}^L)$ are compactly supported and $\|\mcf\|_{E^{m_0}(\bt_0; \sym^2)},$ $\|\mcf^\Psi\|_{E^{m_0}(\bt_0; \bold{B}^L)}< c_0$. 
For $c_0$ sufficiently small, there exists a unique solution $(u, \phi)$ satisfying the equation \eqref{pernon} and 
\beqq\label{eqstab}
\|u\|_{E^{m_0}(\bt_0; \sym^2)} \leq Cc_0, \ \ \|\phi\|_{E^{m_0}(\bt_0; \bold{B}^L)} \leq Cc_0 
\eeqq
where $C$ denotes a generic constant.
\end{prop} 

Finally, we consider the well-posedness of the full Einstein equations \eqref{einsour0} on the data set $\mcd(\delta)$. It is known (see e.g.\ \cite{KLU1}) that  solutions of the reduced Einstein equations give solutions to the Einstein equation only if the conservation law holds. Assume this is the case. Our purpose is to clarify the regularity requirements of the data set. For our later analysis, we would like to take $\mcf \in E^{m_0}(M; \sym^2)$ with $m_0 \geq 7$ in wave gauge so that $\mcf \in C^4(M; \sym^2)$. From Prop.\ \ref{stabest}, we need $m_0$ to be an even integer ($m_0\geq 8$) and get $g\in E^{m_0}(M; \sym^2)$. Now consider the wave map $f$. If $g' \in C^{m}$ then $f\in E^{m}$. Since we have $g = f_* g'$ and $\mcf  = f_* \mcf'$, we need $m \geq m_0 + 1 \geq 9$ and $\mcf' \in E^7(M; \sym^2) \subset C^4(M; \sym^2)$.

\section{Microlocal linearization stability}\label{sec-muls}
Let us consider the linearization of the reduced Einstein equations \eqref{pernon} in wave gauge. Suppose that the source perturbations $\mcf, \mcf^\Psi$  in \eqref{pernon} are $\mcf_\eps, \mcf^\Psi_{\eps}$ depending on a small parameter $\eps >0$ and 
 \begin{enumerate} 
 \item[(L1)]  $\mcf_\eps|_{\eps = 0} = 0, \ \  \p_{\eps} \mcf_\eps |_{\eps = 0} = f \in E^{m_0}(\bt_0; \sym^2), \quad m_0\geq 8$. 
 \item[(L2)] $\mcf^\Psi_\eps|_{\eps = 0} = 0, \ \  \p_{\eps} \mcf^\Psi_\eps |_{\eps = 0} = f^\Psi \in E^{m_0}(\bt_0; \bold{B}^L), \quad m_0\geq 8$. 
 \end{enumerate}
 The solution $u_\eps, \phi_\eps$ of \eqref{pernon} with such $\mcf_\eps, \mcf^\Psi_\eps$ depends on $\eps$ by Prop.\ \ref{stabest} and it is clear that $u_\eps|_{\eps=0} = \phi_\eps|_{\eps = 0} = 0$. Let $v = (\dot g, \dot \phi) =  (\p_\eps u_\eps |_{\eps = 0}, \p_\eps \phi_\eps|_{\eps=0})$ be the linearized field. Then $v$ satisfies the linearized reduced Einstein equations of \eqref{pernon}, which are
\beqq\label{perlin}
\begin{gathered}
-\hat g^{pq} \frac{\p^2 }{\p x^p x^q} \dot g + \mcb(x, v , \p \dot g ) =  2 \rho(\hat g, f) + R,\text{ in } M(\bt_0),\\ 
\square_{\hat g} \dot \phi +  \mcb^\Psi(x, v) = f^\Psi, \\
v  = 0 \text{ in } M(\bt_0)\backslash J^+_{\hat g}(\supp(\bar f)),
\end{gathered}
\eeqq
where  $\bar f = (f, f^\Psi)$;  $\mcb $ is a first order linear differential operator  with smooth coefficients and section valued in $\sym^2$; $\mcb^\Psi$ is smooth in its arguments valued in $\bold{B}^L$; and $R$ is  smooth  valued in $\sym^2$. It is convenient to write the system \eqref{perlin} in matrix form. One can identify $\sym^2\oplus\bold{B}^L$ as a vector bundle $\mcs$ with fiber isometric to  a $10+ L$-dimensional vector space by renumbering the components. From \eqref{eqbeltra}, we see that the term $-\hat g^{pq} \p_{pq}^2$ is comparable to the wave operator $\square_{\hat g}$. Then the system \eqref{perlin} can be written as
\beqq\label{lineq}
\bold{P}_{\hat g}v \doteq 
\begin{bmatrix}
\square_{\hat g} \bold{Id}_{10} + \mcb_1 & \mcb_2 \\[5pt]
\mcb^\Psi_2 & \square_{\hat g} \bold{Id}_L + \mcb^\Psi_1 
\end{bmatrix} v = 
\begin{bmatrix} 
2\rho(\hat g, f) + R \\[5pt]
f^\Psi
\end{bmatrix} = F,
\eeqq
where $\bold{Id}_\bullet$ denotes the $\bullet\times \bullet$ identity matrix, $\mcb_1$ is a first order differential operator and $\mcb_2, \mcb^\Psi_1,  \mcb^\Psi_2$ are smooth function of $v$. It is known that if $(M, \hat g)$ is globally hyperbolic, there exists a causal inverse for $\square_{\hat g}$, see e.g. \cite{Bar}. We next describe the microlocal structure of the causal inverse.

Let $X$ be an $n$-dimensional smooth manifold. For two Lagrangians $\La_0, \La_1 \subset T^*X$ intersecting cleanly at a co-dimension $k$ submanifold i.e. $T_q\La_0\cap T_q\La_1  = T_q(\La_0\cap \La_1),\ \ \forall q\in \La_0\cap \La_1,$ 
the paired Lagrangian distribution associated with $(\La_0, \La_1)$ is denoted by $I^{p, l}(\La_0, \La_1)$. We refer the reader to \cite{DUV, MU, GrU93, GrU90} for the precise definition and details. Let $\mcp(x, \xi) = |\xi|^2_{\widehat g^*}$ be the symbol of $-\widehat g^{pq}  \p_p \p_q$ with $\widehat g^* = \widehat g^{-1}$ the dual metric. Let $\Sigma_{\widehat g} = \{(x, \xi)\in T^*M: \mcp(x, \xi) = 0\}$ be the characteristic set. Notice that $\Sigma_{\widehat g}$ consists of light-like co-vectors. The Hamilton vector field of $\mcp$ is denoted by $H_\mcp$ and in local coordinates
\beq
H_\mcp = \sum_{i = 1}^4( \frac{\p \mcp}{\p \xi_i}\frac{\p }{\p x_i} - \frac{\p \mcp}{\p x_i}\frac{\p }{\p \xi_i}).
\eeq
The integral curves of $H_\mcp$ in $\Sigma_{\widehat g}$ are called null bi-characteristics and their projections to $M$ are geodesics. Let $\diag = \{(z, z')\in M\times M: z = z'\}$ be the diagonal and 
\beq
N^*\diag = \{(z, \zeta, z', \zeta')\in T^*(M\times M)\backslash 0: z = z', \zeta' = -\zeta\}
\eeq 
be the conormal bundle of $\diag$ minus the zero section. We let $\La_{\widehat g}$ be the Lagrangian obtained by flowing out $N^*\diag\cap \Sigma_{\widehat g}$ under $H_\mcp$. Here,  we regard $\Sigma_{\widehat g}, H_\mcp$ as objects on  the product manifold $T^*M\times T^*M$ by lifting from the left factor. From \cite[Lemma 3.1]{KLU1} we conclude that there exists a causal inverse $\bfq_{\hat g} \in I^{-\frac{3}{2}, -\ha}(N^*\diag, \La_{\hat g})$ of $\bold{P}_{\hat g}$ such that $\bold{P}_{\hat g}\bfq_{\hat g} = \bold{Id}$ on $\mce'(M; \mcs)$. Hereafter, whenever the dependence on $\hat g$ is not important, we shall abbreviate the notations as $\bold{P}, \bold{Q}.$ Recall that we define Sobolev spaces on $M$  using the Riemannian metric $\hat g^+$. From Prop.\ 5.6 of \cite{DUV} or Theorem 3.3 of \cite{GrU90}, we know that $\bfq: H_{\comp}^{m}(M; \mcs)\rightarrow H^{m+1}_{\loc}(M; \mcs)$ is continuous for $m\in \mbr$.\\

Next we construct distorted plane wave solutions as in \cite{KLU}. These are conormal distributions, see H\"ormander \cite{Ho3, Ho4}. 
Suppose that $\La$ is a smooth conic Lagrangian submanifold of $T^*X\backslash 0$. Following standard notation, we denote by $I^\mu(\La)$ the Lagrangian distribution of order $\mu$ associated with $\La$. 
For $u\in I^\mu(\La)$, we know that the wave front set $\WF(u)\subset \La$ and $u\in H^s(X)$ for any $s< -\mu-\fnf$. The principal symbol of $u$ is well-defined as a half-density bundle tensored with the Maslov bundle on $\La$, see \cite[Section 25.1]{Ho4}. For a submanifold $Y\subset M$, the conormal distributions to $Y$ are denoted by $I^\mu(N^*Y)$. 

Let $(x_0, \theta_0)\in L^+M(\bt_0)$.  For $s_0>0$ a small parameter, we let 
\beq
K(x_0, \theta_0; s_0) = \{\gamma_{x_0, \theta'}(t)\in M(T_0); \theta' \in \mco(s_0), t\in (0, \infty)\},
\eeq
where $\mco(s_0) \subset L^+_{x_0}M$ is a open neighborhood of $\theta_0$ consisting of $\zeta\in L_{x_0}^+M$ such that $\|\zeta - \theta_0\|_{\hat g^+}< s_0$. Notice that as $s_0\rightarrow 0$, $K(x_0, \theta_0;  s_0)$ tends to the geodesic $\gamma_{x_0, \theta_0}$.  Next, let 
\beqq\label{ymani}
Y(x_0, \theta_0; s_0) = \{\gamma_{x_0, \theta'}(t)\in M(T_0); \theta' \in \mco(s_0), t=\|\theta_0\|_{\hat g^+}\} 
\eeqq
be a $2$-dimensional surface. We let $\La(x_0, \theta_0; s_0)$ be the Lagrangian submanifold obtained by flowing out $N^*K(x_0, \theta_0;  s_0)\cap N^*Y(x_0, \theta_0; s_0)$ under the Hamilton vector field of $\sigma(\square_{\hat g})$ in $\Sigma_{\hat g}$. More concisely, we have
\beqq\label{lamani}
\La(x_0, \theta_0;  s_0) = \La_{\hat g}\circ (N^*K(x_0, \theta_0;  s_0)\cap N^*Y(x_0, \theta_0;  s_0) \cap \Sigma_{\hat g}),
\eeqq
where we used $\circ$ for composition of sets as relations. By our assumption that there are no conjugation points on $(M, \hat g)$, $K(x_0, \theta_0;  s_0)$ is a co-dimension one  submanifold near $\gamma_{x_0, \theta_0}(t)$ and 
\beq
\La = N^*K(x_0, \theta_0; s_0) \text{ near } \gamma_{x_0, \theta_0}(t).
\eeq
Applying \cite[Prop.\ 2.1]{GrU93}, we obtain the following result, see also \cite[Lemma 3.1]{KLU}.

\begin{prop}\label{distor}
Suppose $Y(x_0, \theta_0; s_0)$ and $\La(x_0, \theta_0;  s_0)$ are defined as in \eqref{ymani} and \eqref{lamani} respectively. For $\bar f = (f, f^\Psi)\in I^{\mu+1}(N^*Y; \mcs)$ with $\mu$ an integer, the solution to the linearized Einstein equations \eqref{lineq} is 
\beq
v = \bfq_{\hat g}({F}) \in I^{\mu - \ha, -\ha}(N^*Y, \La; \mcs).
\eeq
In particular, away from $Y$, we have $v\in I^{\mu-\ha}(\La; \mcs)$ and microlocally away from $\La$, we have $v\in I^{\mu - 1}(Y; \mcs)$. 
Moreover, for $(p, \xi)\in N^*Y\cap \Sigma_{\hat g}$ and $(q, \eta)\in \La\backslash N^*Y$ which lie on the same bi-charactersitics, the symbol of $v$ is given by
\beq
\sigma(v)(q, \eta) = \sigma(\bfq_{\hat g})(q, \eta; p, \xi) \sigma({F})(p, \xi),
\eeq
where $\sigma(\bfq_{\hat g})$ is an invertible linear map.
\end{prop}

Now we are ready to define the microlocal linearization stability condition.  This concept for the Einstein equation was introduced in Choquet-Bruhat and Deser \cite{CbD} and developed in Fischer-Marsden \cite{FM1} and Moncrief \cite{Mo}, see also \cite{Cb, GB} for additional material background and references. In the context of Einstein equations with matter, the condition is formulated by Girbau and Bruna \cite{GB}. 

Consider the conservation law  
\beq
\div_g (\hat T(g, \psi) + \mcf) = 0
\eeq
for the Einstein equation \eqref{einsour0}. Now we ignore the field equations for $\psi$ and think of $\psi, \mcf$ are the sources for the Einstein equations \eqref{einsour0}. We assume that 
\begin{gather*}
 \psi_\eps|_{\eps = 0} = \hat \psi , \ \ \p_{\eps} \psi_{\eps}|_{\eps = 0} = \dot \psi \in C^4(M(\bt_0); \bold{B}^L),
\end{gather*}
which is a derived version of (L2). Let $\mcf$ be as in (L1). Then from Section \ref{sec-redein} we know that the solution $g$ of the reduced Einstein equations in wave gauge can be written as $g_\eps = \hat g + \eps \dot g + O(\eps^2)$ in $H^4(M(\bt_0))$. We obtain the linearized conservation law
\beqq\label{mcons} 
 \hat g^{pk} \hat \nabla_p f_{kj} + P_{j}(\dot g, \p\dot g, \dot \psi, \p \dot \psi) = 0, \ \ j = 0, 1, 2, 3,
\eeqq
where $P_j$ is a smooth function and linear in $\dot g, \p\dot g, \dot \psi, \p \dot \psi$.  In our context, the Einstein equation \eqref{einsour0} is called {\em linearization stable} at $(\widehat g, \widehat \psi)$ if for any $f, \dot \psi$ and $\dot g$ satisfying \eqref{mcons}, one can find a family $\mcf_\eps, \psi_\eps$ depending on $\eps$  such that 
\begin{enumerate}
\item  $\mcf_\eps|_{\eps = 0} = 0, \p_\eps \mcf_\eps|_{\eps = 0} = f $ and $\psi_\eps|_{\eps = 0} = \widehat \psi,  \p_\eps \psi_{\eps}|_{\eps = 0} = \dot \psi $;
\item  there is a solution  $g_\eps$ to \eqref{einsour0} so that $(g_\eps, \psi_\eps, \mcf_\eps)\in \mcd(\delta)$ and $g_\eps|_{\eps = 0} = \widehat g,  \p_\eps g_\eps|_{\eps = 0} = \dot g$.
\end{enumerate}

In the case that $f$ has conormal singularities, the condition can be microlocalized as follows. Let $f\in I^{\mu+ 1}(N^*Y; \sym^2), f^\Psi \in  I^{\mu+ 1}(N^*Y; \bold{B}^L)$. By Prop.\ \ref{distor}, we know that $\dot g\in I^{\mu - 1}(N^*Y; \sym^2), \dot \psi \in I^{\mu - 1}(N^*Y; \bold{B}^L)$ microlocally away from $\La$. By Lemma 5.13 of \cite{DUV}, we know that $\p \dot g \in I^{\mu}(N^*Y, \sym^2)$ away from $\La$. Since $P_j(\bullet)$ is linear in $\dot g, \p \dot g$ and $\dot \psi, \p \dot \psi$, we see that $P_{j} \in I^{\mu}(N^*Y; \sym^2)$ away from $\La$. However, $\hat \nabla_p f_{kj} \in I^{\mu+ 2}(N^*Y; \sym^2)$ and $\hat \nabla_p f_{kj} - \p_p f_{kj} \in  I^{\mu+ 1}(N^*Y; \sym^2)$. Therefore, the leading singularity of the linearized conservation law is in $\hat g^{pk}\p_p f_{kj}$. Hence we obtain the microlocal linearized conservation law, that is, the principal symbols of \eqref{mcons} should satisfy 
\beqq\label{mcons1}
\hat g^{pk} \eta_p \sigma(f_{kj})(y, \eta) = 0, \ \ j = 0, 1, 2, 3,
\eeqq
for any $(y, \eta)\in N^*Y\backslash \La$. 
Now we define the linearization stability condition. 
\begin{definition}[Microlocal linearization stability]\label{mlstab}
We say that the microlocal linearization stability condition holds for $\mcd(\delta; V)$ if the followings are true: 
 
Let $Y$ be a two-dimensional submanifold in a space-like hypersurface of $M$. Given any $(y, \eta)\in  N^*Y \cap \Sigma_{\hat g}$, compact neighborhood $W$ of $y$ and $4\times 4$ symmetric matrix $A$  satisfying
\beqq\label{mlseq}
\hat g^{pk}(y)  \eta_p A_{kj} = 0, \ \ j = 1, 2, 3, 4,
\eeqq
there is a family $\mcf_\eps, \mcf^\Psi_\eps$ for $\eps$ small and solutions $g_\eps, \psi_\eps$ to the Einstein equations \eqref{einsour00}  such that 
$(g_\eps, \psi_\eps, \mcf_\eps) \in \mcd(\delta; V)$. In addition, for some $\mu < -17$, 
\beq
\begin{gathered}
\mcf_\eps|_{\eps = 0} = 0, \ \ f \doteq \p_\eps \mcf_\eps|_{\eps = 0} \in I^{\mu+1}(N^*Y; \sym^2), \\
\mcf^\Psi_\eps|_{\eps = 0} = 0, \ \  \p_\eps \mcf^\Psi_\eps|_{\eps = 0} \in I^{\mu+1}(N^*Y; \bold{B}^L),
\end{gathered}
\eeq
and  the principal symbol $\sigma(f)(y, \eta) = A$.  
In addition, we assume that 
\begin{enumerate}
\item  The family $\mcf_\eps, \mcf^\Psi_\eps$ are supported in $W$.
\item The conservation law $\div_g T_{sour} = 0$ holds where $\mcf, \mcf^\Psi$ vanish. 
\end{enumerate}
\end{definition}  

\begin{remark}
The condition is non-trivial. The additional assumptions (1) and (2) are made for convenience in the consideration of wave interactions. In fact, we only need them in Section \ref{subsec-sour} to construct sources of distorted plane waves associated with four small parameters. For the examples studied in Sections \ref{sec-einscal} and \ref{sec-einmax}, (2) is automatically satisfied.  Examples for which (1) holds are constructed in \cite{KLU1} for the Einstein-scalar field equations. One needs sufficiently many scalar fields $\psi$ to construct the family of sources. These additional assumptions can be dropped if the source construction in Sections \ref{sec-einscal} can be analyzed using other methods. 
\end{remark}

Finally, we discuss the symbol spaces of the sources and distorted plane waves. 
We know that the symbols of the sources should satisfy the microlocal linearized conservation law \eqref{mcons1}. 
In the following, we  use the notation $I^\mu_C(N^*Y; \sym^2)$ to denote distributions in $I^\mu(N^*Y; \sym^2)$ whose principal symbols satisfy \eqref{mcons1} for $(y, \eta)\in N^*Y$. Let $\sym^2_y$ be the fiber of $\sym^2$ at $y\in M$ and $\mcy_{y, \eta}$ be the space of elements in $\sym^2_y$ satisfying \eqref{mcons1}. We observe that $\mcy_{y, \eta}$ is a vector space of dimension $6$. 

Next, since we work in wave gauge, the metric $g$ should satisfy the harmonic gauge condition $\hat \Gamma^k = \Gamma^k$. So $\dot g$ should satisfy the linearized harmonicity condition, which is 
\beqq\label{mgauge}
  -\hat g^{\alpha\beta}\p_{\alpha} \dot g_{\beta \mu} + \ha \hat g^{\alpha\beta} \p_{\mu} \dot g_{\alpha\beta} = h^{\alpha\beta}_\mu\dot g_{\alpha\beta}, \ \ \mu = 0, 1, 2, 3,
\eeqq
where $h$ depends on $\hat g$ and its derivatives. This is found in \cite[Section 3.2.3]{KLU1}. Now consider $f\in I^{\mu+1}(N^*Y; \sym^2)$ whose symbols satisfy the microlocal linearized conservation law, then the symbol of  $\dot g$ should satisfy the microlocal linearized gauge condition, which is 
\beqq\label{mgauge1}
  -\hat g^{\alpha\beta} \xi_\alpha \sigma(\dot g_{\beta \mu})(x, \xi) + \ha \hat g^{\alpha\beta}  \xi_\mu \sigma(\dot g_{\alpha\beta})(x, \xi) = 0, \ \ \mu = 0, 1, 2, 3
\eeqq
for $(x, \xi)\in \La$. Here, we used that the leading singularities are on the left hand side of \eqref{mgauge}. We let $\mcx_{x, \xi}$ be the space of elements in $\sym^2_x$ which satisfy \eqref{mgauge1}. We notice that $\mcx_{x, \xi}$ is of dimension $6$. Now in view of Prop.\ \ref{distor} and \eqref{lineq}, we conclude that the principal symbol 
\beq
\sigma(\bfq): \mcy_{y, \eta} \rightarrow \mcx_{x, \xi} \text{ is a linear bijection.}
\eeq 
Here the symbol map is regarded as the restriction from $\sym^2_y$ to $\sym^2_x$. 
We take this conclusion as an important consequence of the microlocal linearization stability condition. It indicates that for the inverse problem we formulated in this work, we need all the possible source perturbations for the Einstein equations in a six dimensional space and measure all possible gravitational perturbations. An interesting question is whether one can solve the problem with fewer sources e.g.\ in a subspace of $\mcy_{y, \eta}$ and/or with measurements of the field $\psi$. This is studied in \cite{LUW1} for the Einstein equations with electromagnetic sources in vacuum, in which case the sources are essentially taken only in a three dimensional space.  This is the key difference between the formulations in \cite{LUW1} and the current work, see Section \ref{sec-einmax}.

\section{Nonlinear interactions of singularities}\label{sec-wgauge} 

We construct sources $\mcf$ in the reduced Einstein equations \eqref{pernon} which generates four distorted plane waves for the linearized Einstein equation \eqref{perlin}. Due to the nonlinearity of the Einstein equations, these waves will interact and we analyze such interactions in this section. 

\subsubsection{Construction of sources and distorted plane waves.}\label{subsec-sour}
Take four points $x^{(i)} \in V, i = 1, 2, 3, 4$ and light-like vectors $\theta^{(i)} \in L^+_{\hat g, x^{(i)}} M$. We denote $\vec x = (x^{(i)})_{i = 1}^4$ and $\vec \theta = (\theta^{(i)})_{i = 1}^4$.  We assume that 
\beqq\label{fourpts}
 x^{(j)} \notin J_g^+(x^{(k)}), \ \ j, k = 1, 2, 3, 4, \ \  j\neq k, 
\eeqq
which means that the points $x^{(i)} $ are  causally unrelated. For $ i= 1, 2, 3, 4$ and $s_0>0$ small, we define
\beqq\label{eqKi}
\begin{gathered}
K_i = K(x^{(i)}, \theta^{(i)};  s_0),\ \ Y_i = Y(x^{(i)}, \theta^{(i)};  s_0) \text{ and } \La_i = \La(x^{(i)}, \theta^{(i)}; s_0),
\end{gathered}
\eeqq
similar to the construction of distorted plane waves in Section \ref{sec-muls}. Then we have $\La_i = N^*K_i$. For $ i = 1, 2, 3, 4$, let 
\begin{enumerate}
\item $f^{(i)} \in I_C^{\mu+1}(N^*Y_i; \sym^2), f^{\Psi, (i)} \in I^{\mu+1}(N^*Y_i; \bold{B}^L), \mu < -17 $ be the linearized sources,
 \item $v^{(i)} = \bfq (F^{(i)}) \in I^{\mu-\ha}(N^*K_i)$  be the distorted plane waves.
 \end{enumerate}
Finally, we take four small parameters  $\eps_i, i = 1, 2, 3, 4$ and set  $\eps = (\eps_1, \eps_2, \eps_3, \eps_4)$. Let $F^{(i)}$ be the $F$ term in \eqref{lineq} with $f$ replaced by $f^{(i)}$. Let $F_{\eps} = \sum_{i = 1}^4 \eps_i F^{(i)}$ and $v_{\eps} = \sum_{i = 1}^4 \eps_i v^{(i)}$ then we have (see \eqref{lineq})
\beq
\bold{P} v_{\eps} = F_{\eps}.
\eeq

\begin{figure}[htbp]
\centering
 
\includegraphics[scale=0.75]{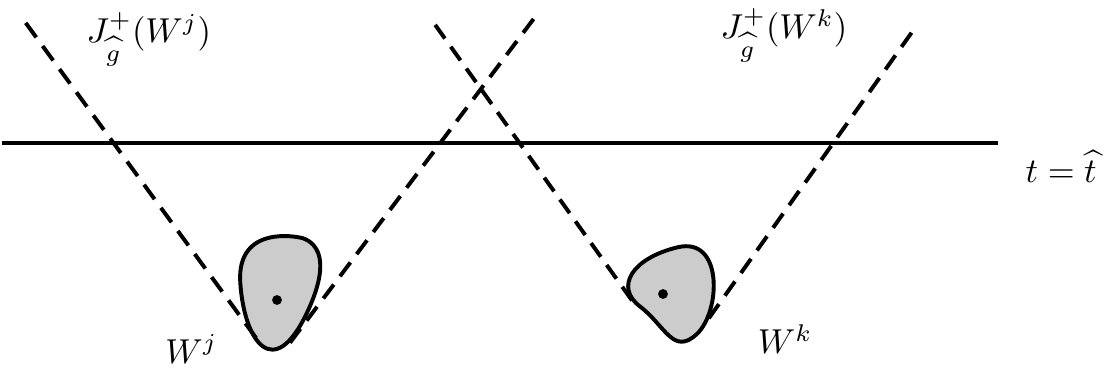}
\caption{The two shaded regions are $W^j, W^k$, which are neighborhoods of the two points inside $x^{(j)}, x^{(k)}$ respectively. }
\label{figcausal}
\end{figure}

Now consider the nonlinear equation \eqref{pernon}. By the microlocal linearization stability condition, there exists $\psi_{\eps_i}^{(i)}, \overline\mcf_{\eps_i}^{(i)}, i = 1, 2, 3, 4$ such that 
\beq
\p_{\eps_i} \overline\mcf^{(i)}_{\eps_i}|_{\eps_i = 0} = 
\begin{bmatrix}
f^{(i)} \\[5pt]
f^{\Psi, (i)} 
\end{bmatrix}.
\eeq
Actually, since the points $x^{(i)}$ are causally unrelated by \eqref{fourpts}, we can choose compact neighborhoods $W^i$ of $x^{(i)}$ which are causally unrelated, namely
\beqq\label{eqWj}
\supp(W^{j})\cap J^+_{\widehat g} (\supp(W^{k})) = \emptyset, \ \ j, k = 1, 2, 3, 4, \ \  j\neq k. 
\eeqq
Let $\overline\mcf_{\eps_i}^{(i)}$ be supported in $W^i$. We take $\overline\mcf_\eps = \sum_{i = 1}^4 \overline\mcf_{\eps_i}^{(i)}$ as the source in \eqref{pernon}. To see that the data $ \overline\mcf_\eps$ admits a solution to the Einstein equation \eqref{einsour00}, 
we proceed as follows. 
We use \eqref{eqWj}. For some $j\neq k$, there exists some $\hat t \in \mbr$ such that $ W^j, W^k \subset M(\hat t)$ and 
\beq
J^+_{\hat g}(W^j)\cap J^+_{\hat g}(W^k)\cap M(\hat t) = \emptyset, 
\eeq
see Figure \ref{figcausal}. 

It is clear that there exists a solution $g_\eps$ to the Einstein equations \eqref{einsour0} or \eqref{ein2} for $t\leq  \hat t$, which is the sum of the two solutions corresponding to $\overline\mcf_{\eps_j}^{(j)}$ and $\overline\mcf_{\eps_k}^{(k)}$. Then one can extend $g_\eps$ to $t >\hat t$ by solving a Cauchy problem for the vacuum Einstein equation with data on the surface $t = \hat t$, because in this region the sources all vanish and by condition (2) of the microlocal linearization stability condition Def.\ \ref{mlstab}, the conservation law holds. This procedure can be continued to include all the sources in consideration. Thus, we  complete the construction of the sources. 

\begin{remark}
In general, the construction of the sources involving four parameters would lead to expansions of $\mcf_\eps$ containing $\eps_i\eps_j, \eps_i\eps_j\eps_k$ and higher order terms. This does not happen here because we used the assumptions (1) and (2) in Definition \ref{mlstab}. However, if the singularities in those expansion terms are known, the analysis below for the interactions can be carried out as well. 
\end{remark}

\subsubsection{The asymptotic expansion}
We let $u_{\eps}$ be the solution of the nonlinear equation \eqref{pernon} with $\mcf_\eps, \psi_\eps$ 
and derive the asymptotic expansion of $u_\eps$ to find the nonlinear interaction terms. We shall see that because $T(g, \psi)$ does not depend on the derivative of $\psi$ and below we only look at the leading order singularities which can only come from nonlinear terms with two derivatives, the $\psi$ components of the distorted plane waves $v$ do not play a role (the derivative of the $g$ components takes the leading order singularity). Thus, from now on, we assume that the field $\psi_\eps$ are smooth and only look at the Einstein equations in \eqref{pernon}. This helps us to focus on the gravitational wave interaction from the Einstein equations, which is the main issue here. Also, one can reuse the analysis and calculations for other models as shown in Sections \ref{sec-einscal} and \ref{sec-einmax}. The derivation bears some similarities to the Einstein-Maxwell case \cite{LUW1}, however, we repeat some arguments for the reader's convenience. 

From equations \eqref{pernon} and \eqref{perlin} which $u_{\eps}, v_{\eps}$ satisfy, we derive
\beqq\label{eqexpan}
\begin{gathered}
\bold{P}(u_\eps - v_\eps) + \sum_{i = 2}^4 P_i(x, u_\eps) + \sum_{i = 2}^4H_i(x, u_\eps) + \mcr_\eps = 0,
\end{gathered}
\eeqq
where $\mcr_\eps$ represents the collection of higher order nonlinear terms in $O(u_\eps^5)$, and the rest are lower order terms. We explain them in detail. First of all, the $P_i, i = 2, 3, 4$ terms come from the quasilinear term $-g^{pq} \p_p\p_q  u_\eps$ in \eqref{pernon}. In a given local coordinates, we can express two tensors by $4\times 4$ matrices. It is easy to see that 
 \beqq
 \begin{split}
 g^{-1} &= (\widehat g + u_\eps)^{-1} = (\Id + \widehat g^{-1} u_\eps)^{-1} \widehat g^{-1}\\
 & = \widehat g^{-1} - \widehat g^{-1} u_\eps \widehat g^{-1} + (\widehat g^{-1} u_\eps)^2 \widehat g^{-1} + (\widehat g^{-1} u_\eps)^3 \widehat g^{-1} + \cdots.
 \end{split}
 \eeqq
Then the $P_i$ terms are
 \beqq\label{eqpi}
 \begin{gathered}
 P_2(x, u_\eps) = (\widehat g^{-1} u_\eps \widehat g^{-1})^{pq} \frac{\p^2 u_\eps}{\p x^p \p x^q}, \\
  P_3(x, u_\eps) = -(\widehat g^{-1} u_\eps \widehat g^{-1} u_\eps \widehat g^{-1})^{pq} \frac{\p^2 u_\eps}{\p x^p \p x^q}, \\
 P_4(x, u_\eps) = (\widehat g^{-1} u_\eps\widehat g^{-1} u_\eps \widehat g^{-1} u_\eps \widehat g^{-1})^{pq} \frac{\p^2 u_\eps}{\p x^p \p x^q}. 
\end{gathered}
 \eeqq
More precisely, we  write down the components of these terms as 
\beqq\label{eqformp}
\begin{split}
P_{2, i}(x,  u_\eps) &= \widehat g^{pa}u_{\eps, ab}\widehat g^{bq} \frac{\p^2 u_{\eps, i}}{\p x^p \p x^q}, \\
P_{3, i}(x,  u_\eps) & = -\widehat g^{pa} u_{\eps, ab} \widehat g^{bc} u_{\eps, cd} \widehat g^{dq} \frac{\p^2  u_{\eps, i}}{\p x^p \p x^q} ,\\
P_{4, i}(x,  u_\eps) & = \widehat g^{pa} u_{\eps, ab} \widehat g^{bc} u_{\eps, cd} \widehat g^{de} u_{\eps, ef} \widehat g^{fq}\frac{\p^2 u_{\eps, i}}{\p x^p \p x^q}, \ \  i = 1, 2, \cdots, 14.
\end{split}
\eeqq
Notice that each term has two derivatives and the coefficients of the derivatives are polynomials of the metric component. For convenience, we also regard $P_{j, i}, j = 2, 3, 4, i = 1, 2, 3, 4$ as multilinear functions. For example, we can write
\beq
\begin{gathered}
P_{2, i}(x,  u^{(1)}, u^{(2)})  = \widehat g^{pa}u^{(1)}_{ab}\widehat g^{bq} \frac{\p^2 u^{(2)}_{i}}{\p x^p \p x^q}, i = 1, 2, 3, 4 \\
\text{ and  } P_{2, i}(x, u_\eps) = P_{2, i}(x,  u_\eps, u_\eps).
\end{gathered}
\eeq
This type of notations are also used below.

Next, the terms $H_i, i = 2, 3, 4$ in \eqref{eqexpan} come from the semilinear terms of \eqref{pernon}.  Each component of $H_i$ is a sum of $i$-th order monomials of $u_{\eps}, \p u_{\eps}$. Here, we slightly abused the notation so that $u_{\eps}$ stands for their components $u_{\eps, ij}, i, j = 0, 1, 2, 3$. This is also used in the following for simplicity.  We point out that $H_i$ are at most quadratic in $\p u_{\eps}$ which can be seen from the expression of the reduced Ricci tensor \eqref{ricre} and the fact that the Christoffel symbols only have first derivatives of the metric $g$, and the assumption on $\hat T(g, \phi)$ that it does not involve $\p g$. We write the component of $H_i$ as multilinear functions as
\beq
\begin{split}
H_{2, \theta}(x, w^{(1)}, w^{(2)}) &= \sum_{i, j = 1}^{14} \sum_{ \alpha, \beta = 1}^4 \sum_{a, b=0}^1 H^\theta_{2, ij\alpha\beta ab} \p_\alpha^a  w^{(1)}_i \p_\beta^b  w^{(2)}_j, \\
H_{3, \theta}(x, w^{(1)}, w^{(2)}, w^{(3)}) &= \sum_{i, j,k = 1}^{14} \sum_{ \alpha, \beta, \gamma = 1}^4 \sum_{a, b, c =0, 1; a+b+c \leq 2} H^\theta_{3, ijk\alpha\beta\gamma abc} \p_\alpha^a w^{(1)}_i \p_\beta^b w^{(2)}_j \p_\theta^c w^{(3)}_k,\\
H_{4, \theta}(x, w^{(1)}, w^{(2)}, w^{(3)}, w^{(4)}) &= \sum_{i, j,k,l  = 1}^{14}\sum_{\alpha, \beta, \gamma, \delta = 1}^4 \sum_{(a, b, c, d) \in \mca} H^\theta_{4, ijkl\alpha\beta\gamma\delta abcd} \p_\alpha^a w^{(1)}_i \p_\beta^b  w^{(2)}_j \p_\gamma^c w^{(3)}_k \p_\delta^d w^{(4)}_l, 
\end{split}
\eeq 
where the set $\mca = \{(a, b, c, d) : a, b, c, d =0, 1; a+b+c+d \leq 2\}$, $\theta = 1, 2, \cdots, 14$ and the coefficients are all smooth. We emphasis that the derivatives  in these terms appear at most twice and such terms are especially important for the analysis below. We denote such terms i.e.\ terms in $H_i$ with two derivatives, by $\widehat H_i, i = 2, 3, 4.$   

For convenience, we denote
\beq
\begin{gathered}
G_i(x, u_\eps) = P_i(x, u_\eps) + H_{i}(x, u_\eps), \ \ i = 2, 3, 4,
\end{gathered}
\eeq
then $G_i$ are polynomials in $u_\eps, \p u_\eps$ of order $i$. We also denote $\widehat G_i = P_i + \widehat H_i.$ Then the asymptotic expansion \eqref{eqexpan} can be written as
\beqq\label{eqite1}
\begin{gathered}
u_\eps =  v_\eps - \bfq( G_2 + G_3 + G_4) + \mcr_\eps.
\end{gathered}
\eeqq
By the stability estimate, we have $\|u_\eps\|_{E^{4}(\bt_0)} \leq C\sum_{i = 1}^4 \eps_i$. Note $E^{4}(\bt_0)\subset H^4(M(\bt_0))$ is an algebra. We have 
\beq
\| u_{\eps, a} u_{\eps, b}\|_{H^{4}(M(\bt_0))} \leq C(\sum_{i = 1}^4 \eps_i)^2, \ \ a, b = 0, 1, 2, 3.
\eeq
Therefore, by the continuity of $\bfq$, we obtained the first term in \eqref{eqite1} i.e.\ $v_\eps$. Notice that all the terms in $\mcr_\eps$ are in $H^{4}(M(\bt_0))$. To get other terms, we will iterate the formula \eqref{eqite1} i.e.\ plug \eqref{eqite1} to the right hand side of \eqref{eqite1} to get 
\beq
u_\eps = v_\eps + \sum_{0\leq a< b\leq 3} \eps_a\eps_b \mcu^{(2)}_{ab} + \sum_{0\leq a< b< c\leq 3} \eps_a\eps_b\eps_c \mcu^{(3)}_{abc} + \eps_1\eps_2\eps_3\eps_4 \mcu^{(4)} + \mcr_\eps,
\eeq
where $\mcr_\eps$ denotes terms in $H^4(M(\bt_0))$ and $\bigcup_{i = 1}^4 O(\eps_i^2)$. The process generates many terms. However, we only need term $\mcu^{(4)}$ of the order $\eps_1\eps_2\eps_3\eps_4$ which can only be obtained from the multiplication of four terms of $v_\eps$ and $\p v_\eps$.

\subsubsection{The nonlinear interactions} 
We recall some results from \cite{LUW, LUW1, WZ} to analyze the singularities in $\mcu^{(4)}$. We make some assumptions on  $K_i, i = 1, 2, 3, 4$ defined in \eqref{eqKi}. Recall that two submanifolds $X, Y$ of $M$ intersect transversally if $T_q X + T_q Y = T_q M, \ \ \forall q\in X\cap Y.$ We assume that
\begin{enumerate}
\item For $1\leq i< j\leq 4$, $K_i$ intersects $K_j$ transversally at co-dimension $2$ submanifolds $K_{ij}$;
\item For $1\leq i<j <k\leq 4$, $K_{ij}$ intersects $K_k$ transversally at co-dimension $3$ submanifolds $K_{ijk}$;
\item $K_{123}$ intersect $K_4$ transversally at a point $q_0$. 
\end{enumerate}
In particular, the last condition means that the four submanifolds intersect at a point $q_0$ and the normal co-vectors $\zeta_i$ to $K_i$ at $q_0$ are linearly independent. We shall denote 
\beq
\begin{gathered}
 \La_{ij} = N^*K_{ij}\backslash 0, \ \ 1\leq i< j\leq 4; \ \ \La_{ijk} = N^*K_{ijk}\backslash 0, \ \ 1\leq i<j<k\leq 4; \ \ \La_{q_0} =  T^*_{q_0}M\backslash 0,
\end{gathered}
\eeq
where $0$ stands for the zero section of $T^*M$.  All of these are conic Lagrangian submanifolds. We introduce the following notations
\beq
\begin{gathered}
\La^{(1)} = \bigcup_{1\leq i \leq 4}  \La_i; \ \ \La^{(3)} = \bigcup_{ 1\leq i<j<k\leq 4} \La_{ijk}.
\end{gathered}
\eeq
For convenience, we introduce a notation that for any $\Gamma \subset T^*M$, $\Gamma^{\hat g}$ denotes the flow out of $\Gamma$ under $\La_{\hat g}$ i.e.\ $\Gamma^{\hat g} = \La_{\hat g}\circ (\Gamma \cap \Sigma_{\hat g})$, where the composition is understood as composition of relations. Finally, 
let $\Theta = \La^{(1)}\cup \La^{(3), \hat g}$ and $\mck = \pi (\Theta)$ where $\pi : T^*(M\times M)\rightarrow M\times M$ denotes the standard projection. 

Assume that $v_i \in I^{\mu}(N^*K_i), i = 1, 2, 3, 4$ are scalar valued distorted plane waves. Let $Q_{\hat g} \in I^{-\frac{3}{2}, -\ha}(N^*\diag, \La_{\hat g})$ be the causal inverse of $\square_{\hat g}$. In the study of semilinear wave equations \cite{LUW}, the singularities of the following terms are carefully analyzed
\beqq\label{fourtho}
\begin{split}
\mcy_1 &= Q_{\hat g}(c(x)v_1v_2v_3v_4), \\  
\mcy_2 &= Q_{\hat g}(a(x)v_1Q_{\hat g}(b(x)v_2v_3v_4)), \\
\mcy_3 &= Q_{\hat g}(b(x)v_1v_2Q_{\hat g}(a(x)v_3v_4)), \\
\mcy_4 &= Q_{\hat g}(a(x)v_1Q_{\hat g}(a(x)v_2Q_{\hat g}(a(x)v_3v_4))), \\
\mcy_5 &= Q_{\hat g}(a(x)Q_{\hat g}(a(x)v_1v_2)Q_{\hat g}(a(x)v_3v_4)),
\end{split}
\eeqq
where $a(x), b(x), c(x)$ are smooth functions.  These terms involve multiplication of four conormal distributions whose singular support intersect at $q_0$. In fact, all the terms in $\mcu^{(4)}$ are combinations of such terms. It is proved in \cite[Prop.\ 3.9]{LUW} that under appropriate assumptions $\mcy_i$ have conormal singularities to $\La_{q_0}$ i.e.\ $q_0$ becomes a point source. The result is generalized in \cite{LUW1} to include the case when the order of $v_i$ are different and we recall it below.

\begin{prop}[Prop.\ 6.1 of \cite{LUW1}]\label{porder}
Let $v_i\in I^{\mu_i}(\La_i), i = 1, 2, 3, 4$ and $\tilde \mu = \sum_{i = 1}^4 \mu_i$. Let $\mcv$ be a smooth vector field. Microlocally away from $\Theta$, we have the following 
\beq\label{internew}
\begin{split}
& (1)\ \ Q_{\hat g}(c(x)v_1v_2v_3v_4) \in I^{\tilde \mu+\frac{3}{2}}(\La_{q_0}^{\hat g}\backslash \Theta)  \\
& (2)\ \ Q_{\hat g}(a(x)v_1\mcv Q_{\hat g}(b(x)v_2v_3v_4)) \in I^{\tilde\mu+ \ha}(\La_{q_0}^{\hat g}\backslash \Theta) \\
& (3)\ \ Q_{\hat g}(b(x)v_1v_2 \mcv Q_{\hat g}(a(x)v_3v_4)) \in I^{\tilde\mu+ \ha}(\La_{q_0}^{\hat g}\backslash \Theta) \\ 
& (4)\ \ Q_{\hat g}(a(x)v_1 \mcv Q_{\hat g}(a(x)v_2 \mcv Q_{\hat g}(a(x)v_3v_4))) \in I^{\tilde \mu-\frac{1}{2}}(\La_{q_0}^{\hat g}\backslash \Theta)\\
& (5)\ \ Q_{\hat g}(a(x)\mcv Q_{\hat g}(a(x)v_1v_2)\mcv Q_{\hat g}(a(x)v_3v_4)) \in I^{\tilde \mu-\frac{1}{2}}(\La_{q_0}^{\hat g}\backslash \Theta).
\end{split}
\eeq
\end{prop}
 
We point out that the set $\Theta$ contains the singularities produced from the interaction of one, two and three distorted plane waves and $\mck$ is the singular support. So in the above proposition (as well as the rest of the analysis), we basically only look at the new singularities produced by the interaction of four distorted plane waves. 
  
In Section 3.5 of \cite{LUW}, the principal symbols of the terms \eqref{fourtho} are found explicitly, and we recall them here. Consider the symbols at $(q, \eta)\in \La_{q_0}^{\hat g}\backslash \Theta$ which is joined with $(q_0, \zeta)\in \La_{q_0}$ by bi-characteristics. We can write $\zeta = \sum_{i = 1}^4\zeta_i$ where $\zeta_i\in N^*_{q_0}K_i$. Let $A_i$ be the principal symbols of $v_i$. By Prop.\ 3.12 of \cite{LUW}, we get 
\beqq\label{eqpsym}
\begin{split}
 \sigma_{\La_{q_0}}(\mcy_j)(q, \eta) &= \sigma_{\La_{\hat g}}(Q_{\hat g})(q, \eta, q_0, \zeta)\mcp_j(\zeta_1, \zeta_2, \zeta_3, \zeta_4) \prod_{i = 1}^4A_i(q_0, \zeta_i), \ \ j = 1, 2, 3, 4, 5, \\
\text{where }& \mcp_1 = (2\pi)^{-3}c(q_0),\\
 & \mcp_2 =   (2\pi)^{-3}a(q_0)b(q_0)   |\zeta_3+\zeta_4|_{{\hat g}^*(q_0)}^{-2},\\
 &  \mcp_3 = (2\pi)^{-3}a(q_0)     b(q_0) |\zeta_2+\zeta_3+\zeta_4|_{{\hat g}^*(q_0)}^{-2},\\
&  \mcp_4 =  (2\pi)^{-3}a^3(q_0)   |\zeta_2+\zeta_3+\zeta_4|_{{\hat g}^*(q_0)}^{-2}  |\zeta_3+\zeta_4|_{{\hat g}^*(q_0)}^{-2},\\
 &  \mcp_5  = (2\pi)^{-3}a^3(q_0) |\zeta_3+\zeta_4|_{{\hat g}^*(q_0)}^{-2}   |\zeta_1+\zeta_2|_{{\hat g}^*(q_0)}^{-2}.
  \end{split}
 \eeqq 
We remark that these are the principal symbols near the intersection point $q_0$, where we can trivialize half density bundles and Maslov bundles using local coordinates. The generalization of the above formulas to the case with derivatives is quite straightforward, see \cite{LUW1} and \cite{WZ}.

With these preparations, we are ready to state and prove the main result of this section.
\begin{prop}\label{intert}
Let $v^{(i)}\in I^{\mu-\ha}(\La_i; \mcs), i= 1, 2, 3, 4$ be distorted plane waves. Suppose that $\bigcap_{j = 1}^4\gamma_{x^{(j)}, \theta^{(j)}}(\mbr_+) = q_0$ and the corresponding tangent vectors at $q_0$ are linearly independent. Let $s_0> 0$ be sufficiently small  such that $K_j$ intersect at $q_0$ only. Let $\sigma(4)$ denote  the  set of permutations of $(1, 2, 3, 4)$. Then the fourth order interaction term $\mcu^{(4)} = \bfq(\mch + \hat \mch)$ such that   
 \beq
 \bfq(\mch) \in I^{4\mu + \frac{3}{2}}(\La^{\hat g}_{q_0}\backslash \Theta, \mcs) \text{ and } \bfq(\hat \mch) \in I^{4\mu+ \ha}(\La^{\hat g}_{q_0}\backslash \Theta, \mcs),
 \eeq
where $\mch = \sum_{i = 1}^5 \mch_{i}$ and $\mch_i$ are given by
\begin{equation}\label{eqmch1}
\mch_{1} = -\sum_{(i, j, k, l) \in \sigma(4)} \widehat G_4(x, v^{(i)}, v^{(j)}, v^{(k)}, v^{(l)})
\end{equation} 
\begin{equation}\label{eqmch2}
\begin{gathered}
\mch_{2} = \sum_{(i, j, k, l) \in \sigma(4)} \bigg( \widehat G_3(x, v^{(i)}, v^{(j)}, \bfq (\widehat G_2(x, v^{(k)}, v^{(l)}))) + \widehat G_3(x, v^{(i)}, \bfq (\widehat G_2(x, v^{(j)}, v^{(k)})), v^{(j)}) \\
+ \widehat G_3(x, \bfq (\widehat G_2(x, v^{(i)}, v^{(j)})), v^{(k)}, v^{(l)}) \bigg) 
\end{gathered}
\end{equation}
\begin{equation}\label{eqmch3}
  \mch_{3} =  \sum_{(i, j, k, l) \in \sigma(4)} \bigg(\widehat G_2(x, \bfq(\widehat G_3(x, v^{(i)}, v^{(j)}, v^{(k)})), v^{(l)}) + \widehat G_2(x, v^{(i)}, \bfq(\widehat G_3(x, v^{(j)},  v^{(k)},  v^{(l)}))) \bigg) 
 \end{equation}
\begin{equation}\label{eqmch4}
 \mch_{4}  = -  \sum_{(i, j, k, l) \in \sigma(4)} \widehat G_2(x, \bfq(\widehat G_2(x, v^{(i)}, v^{(j)})), \bfq(\widehat G_2(x, v^{(k)},  v^{(l)})))
 \end{equation}
 \begin{equation}\label{eqmch5}
 \begin{gathered}
  \mch_{5}  = -\sum_{(i, j, k, l) \in \sigma(4)} \bigg(\widehat G_2(x, v^{(i)}, \bfq(\widehat G_2(x,  v^{(j)}, \bfq(\widehat G_2(x, v^{(k)},  v^{(l)}))))) \\
+ \widehat G_2(x, v^{(i)}, \bfq(\widehat G_2(x, \bfq(\widehat G_2(x, v^{(j)},v^{(k)})), v^{(l)}))) \\ 
  + \widehat G_2(x, \bfq(\widehat G_2(x, v^{(i)}, \bfq(\widehat G_2(x, v^{(j)}, v^{(k)})))), v^{(l)})\\
+ \widehat G_2(x, \bfq(\widehat G_2(x, \bfq(\widehat G_2(x, v^{(i)}, v^{(j)})), v^{(k)})), v^{(l)})\bigg).
\end{gathered}
\end{equation}  
\end{prop} 

\bpf
First of all, we put \eqref{eqite1} to the right hand side of \eqref{eqite1} to get one term from $G_4(x, u_\eps)$: 
\beq
\begin{split}
&G_4(x, u_\eps) = G_4(x, v_\eps) + \mcr_\eps = \sum_{(i, j, k, l) \in \sigma(4)} G_4(x, v^{(i)}, v^{(j)}, v^{(k)}, v^{(l)}) + \mcr_\eps.
\end{split}
\eeq
The summation terms consist of two types of terms because $G_4 = P_4 + H_4$. The terms from $P_4$ can be found in \eqref{eqformp} and they all have two derivatives. Using Prop.\ \ref{porder}, we conclude that $\bfq (P_4(x, v^{(i)},  v^{(j)},  v^{(k)},  v^{(l)})) \in I^{4\mu + \frac{3}{2}}(\La^{\widehat g}_{q_0}\backslash \Theta)$. The terms from $H_4$ are of the form
\beq
A_{ijkl\alpha\beta m n}v^{(a)}_i v^{(b)}_j \p^m_\alpha   v^{(c)}_k \p_\beta^n  v^{(d)}_l, \ \ m, n\leq 1; \ \ i, j, k, l = 1, \cdots, 14; \ \ \alpha, \beta = 1, 2, 3, 4,
\eeq
and $a, b, c, d$ are permutations of $1,2,3,4$ (Note this is not in Einstein summation.) We can apply Prop.\ \ref{porder} to conclude that when $m, n = 1$, the term after applying $\bfq$ is in $I^{4\mu + \frac{3}{2}}(\La^{\widehat g}_{q_0}\backslash \Theta)$ and otherwise in $I^{4\mu + \ha}(\La^{\widehat g}_{q_0}\backslash \Theta)$. When $m = n = 1$, the terms only come from $\widehat H_4$. Thus we obtain the leading term $\mch_{1}$.

Second, we consider the other terms in the asymptotic expansion. To get order $\eps_1\eps_2\eps_3\eps_4$ terms, we need to iterate twice or three times using \eqref{eqite1}. From the term $G_3(x, u_\eps)$, we get 
\beqq\label{eqh3}
\begin{gathered}
 G_3(x, u_\eps)  = -\sum_{(i, j, k, l) \in \sigma(4)} \bigg( G_3(x, v^{(i)}, v^{(j)}, \bfq (G_2(x,  v^{(k)},  v^{(l)}))) \\
+ G_3(x,  v^{(i)}, \bfq (G_2(x,   v^{(j)},  v^{(k)})),  v^{(j)}) + G_3(x,  \bfq (G_2(x,  v^{(i)},   v^{(j)})),   v^{(k)},   v^{(l)})\bigg) + \mcr_\eps.
\end{gathered}
\eeqq
From Prop.\ \ref{porder}, we know that the term after applying $\bfq$ is in $I^{4\mu+\frac{3}{2}}(\La^{\widehat g}_{q_0}\backslash \Theta)$ if the $G_i, i = 2, 3$ terms involved have two derivatives i.e. they are $\widehat G_i, i = 2, 3.$ Otherwise, the terms are in $I^{4\mu +\ha}(\La^{\widehat g}_{q_0}\backslash \Theta)$ which are less singular. So we get one piece in $\mch_{2}$.  

Finally, from the term $G_2(x, u_\eps),$ we obtain from the iteration using \eqref{eqite1} that 
\begin{equation}\label{eqg2}
\begin{split}
 &G_2(x, u_\eps)  \\
 = &G_2(x, v_\eps - \bfq(G_2(x, u_\eps) + G_3(x, u_\eps))) + \mcr_\eps \\
= & G_2(x, \bfq(G_2(x,  v_\eps)), \bfq(G_2(x, v_\eps))) -  G_2(x, v_\eps, \bfq(G_2(x, u_\eps) + G_3(x, v_\eps))) \\
&  - G_2(x, \bfq(G_2(x, u_\eps) + G_3(x, v_\eps)), v_\eps)  + \mcr_\eps,\\
 = &  G_2(x, \bfq(G_2(x, v_\eps)), \bfq(G_2(x,  v_\eps))) -  G_2(x, v_\eps, \bfq(G_3(x, v_\eps)))   - G_2(x, \bfq(G_3(x, v_\eps)), v_\eps)  \\
 &  -  G_2(x, v_\eps, \bfq(G_2(x, v_\eps -  \bfq(G_2(x,  v_\eps)))))  - G_2(x, \bfq(G_2(x, v_\eps - \bfq(G_2(x, v_\eps))), v_\eps) + \mcr_\eps\\
  = &   G_2(x, \bfq(G_2(x,  v_\eps)), \bfq(G_2(x, v_\eps))) -  G_2(x, v_\eps, \bfq(G_3(x,  v_\eps)))   - G_2(x, \bfq(G_3(x, v_\eps)), v_\eps)  \\
 &  +  G_2(x, v_\eps, \bfq(G_2(x, v_\eps, \bfq(G_2(x, v_\eps))))) + G_2(x, \bfq(G_2(x, v_\eps, \bfq(G_2(x, v_\eps)))), v_\eps)  \\
 & +  G_2(x, \bfq(G_2(x, v_\eps, \bfq(G_2(x, v_\eps))), v_\eps) + G_2(x, \bfq(G_2(x, \bfq(G_2(x, v_\eps)), v_\eps), v_\eps) + \mcr_\eps. 
\end{split}
\end{equation}
Using $v_\eps = \sum_{i = 1}^4 \eps_i v^{(i)}$, we find that 
\beq
\begin{gathered}
 G_2(x, v_\eps, \bfq(G_3(x, v_\eps)))  + G_2(x, \bfq(G_3(x, v_\eps)), v_\eps)\\
 =  \sum_{(i ,j, k, l) \in \sigma(4)} \bigg( G_2(x,  v^{(i)}, \bfq(G_3(x, v^{(j)},  v^{(k)},  v^{(l)})))   + G_2(x, \bfq(G_3(x,  v^{(i)},  v^{(j)},  v^{(k)})),  v^{(l)}) \bigg).
\end{gathered}
\eeq
Applying Prop.\ \ref{porder}, we see that the leading order term is achieved when the $G_i, i = 2, 3$ are $\widehat G_i$ and the leading terms are in $I^{4\mu + \frac{3}{2}}(\La^{\widehat g}_{q_0}\backslash \Theta)$. So we get the other piece of $\mch_2$. Using the same argument, we can obtain the terms in $\mch_3$ from the rest of terms in \eqref{eqg2}. The details are omitted here. 
\epf
 
If $K_i, i = 1, 2, 3, 4$ do not intersect, the singularities of $\mcu^{(4)}$ are at most conic and they live on $\mck\subset M$. In particular, $\mck$ is the set in $M$ carrying the singularities produced by three wave interactions. The proof of the proposition below is the the same as that of Prop.\ 4.1 of \cite{LUW} for the scalar case, so we just state the result without proving it. 
\begin{prop}\label{inter5}
If $\bigcap_{j = 1}^4\gamma_{x^{(j)}, \theta^{(j)}}(\mbr_+) = \emptyset$, then $\mcu^{(4)}$ is smooth away from the set $\mck$. If $\bigcap_{j = 1}^4\gamma_{x^{(j)}, \theta^{(j)}}(\mbr_+)  = q_0$  but the tangent vectors at $q_0$ are linearly dependent,  then $\mcu^{(4)}$ is smooth away from $J_{g}^+(q_0)$. 
\end{prop}

 \subsubsection{Non-vanishing of the principal symbol}
Finally, we show that the newly generated singularities are non-vanishing or more precisely the principal symbol of $\mcu^{(4)}$ is not vanishing. We follow the strategy in \cite[Section 4]{KLU1} and \cite{LUW, LUW1}. 

\begin{prop}\label{symb}
Under the assumption of Prop.\ \ref{intert}, the principal symbol of the components of $\sigma(\mcu^{(4)})(\zeta)$ for  $\zeta\in \La_{q_0} \backslash \Theta$ are analytic functions defined on 
\beq
\begin{gathered}
\mcx = \{(\zeta^{(1)}, \zeta^{(2)}, \zeta^{(3)},\zeta^{(4)}, \zeta, A^{(1)}, A^{(2)},A^{(3)},A^{(4)}):
\zeta\in \La_{q_0}\backslash \Theta, \\
\zeta^{(i)} \in L_{q_0}^*M \text{ are linearly independent},  i =  1, 2, 3, 4,\\
\text{ and } A^{(i)} \text{ are $10\times 10$ matrices which are the principal symbols of } v^{(i)}\}.
\end{gathered}
\eeq
Moreover, the principal symbol $\sigma(\mcu^{(4)})$ is non-vanishing on any relatively open subset of $\mcx$.
\end{prop}

\bpf The proof is divided into six steps. 

\textbf{Step 1:} We first observe that according to \eqref{eqmch1}--\eqref{eqmch5}, all the terms in $\mcu^{(4)}$ in Prop.\ \ref{intert} are linear combinations of the terms listed in \eqref{fourtho}. Using the formulas for the principal symbols \eqref{eqpsym}, we can easily see that the principal symbols $\sigma(\mcu^{(4)})$ is a polynomial function of the symbols $\sigma(v^{(i)}_{jk}), i = 1, 2, 3, 4, j, k = 0, 1, 2, 3.$ The coefficients are rational functions of $\zeta^{(i)}, i = 1, 2, 3, 4$ and $\zeta$. Thus we see that the principal symbols are real analytic functions defined on $\mcx$.

\textbf{Step 2:} We compute the symbols explicitly to show that they are non-vanishing. Since it is always possible to choose local coordinates so that the background metric $\hat g$ is Minkowski at $q_0$, we analyze the leading singularities in $\mcu^{(4)}$ when the background metric $(M, \hat g)$ is the Minkowski space-time $(\mbr^4, h)$, where $h = -(dx^0)^2 + \sum_{i = 1}^3 (dx^i)^2.$ In this case, the Einstein equations are simpler because $\hat \Gamma_{ij}^k$ and the derivatives of $\hat g_{ij}$ all vanish. 

Recall that $u = g - \hat g$ is the perturbation. The reduced Ricci tensor can be written as 
\beqq\label{eqricre}
\begin{gathered}
(\ric_{\hat g}(g))_{\mu\nu}  = - \ha g^{pq} \frac{\p^2 u_{\mu\nu}}{\p x^p \p x^q} + g^{ab}g_{ps}\Gamma^p_{\mu b}\Gamma^s_{\nu a} + \ha (g_{\nu l} \Gamma^l_{ab} g^{aq}g^{bd} \frac{\p u_{qd}}{\p x^\mu} + g_{\mu l} \Gamma^l_{ab} g^{aq} g^{bd} \frac{\p u_{qd}}{\p x^\nu}),
\end{gathered}
\eeqq 
where the Christoffel symbols can be simplified as 
\beq
\Gamma_{\alpha\beta}^\mu = \ha g^{\mu\la}(\frac{\p u_{\la\alpha}}{\p x^\beta} + \frac{\p u_{\la\beta}}{\p x^\alpha}- \frac{\p u_{\alpha\beta}}{\p x^\la}). 
\eeq
Next $g^{\alpha\beta}$ can be computed as following
\beq
\begin{gathered}
g^{-1} = (h + u)^{-1} = (\Id - h^{-1}u + (h^{-1}u)^2 - (h^{-1}u)^3 + \cdots) h^{-1}. 
\end{gathered}
\eeq 
Since $h = h^{-1}$, we obtain 
\beqq\label{eqgexpan}
g^{ab} = h^{ab} - \sum_{a, b = 0}^3 h^{aa}h^{bb}u_{ab} + \sum_{a, b, c = 0}^3 h^{aa}h^{bb}h^{cc}u_{ac}u_{cb} + \cdots. 
\eeqq
With these formulas, we can find $P_i, \hat H_i, i = 2, 3, 4$ in Prop.\ \ref{intert} explicitly. Recall that we assumed that there is no derivative of $g$ in $T$. The linearized reduced Einstein equations \eqref{perlin} in this case are
\beq
\begin{gathered}
\square_h u + Vu + W = -2f +  \tr_h(f)h,
\end{gathered}
\eeq
where $Vu$ is the linear term with $V$ a $10\times 10$ smooth matrix and $W$ a smooth matrix. We remark that the term $Vu$ comes from the stress-energy tensor term.

\textbf{Step 3:} We claim that for some choice of $\zeta^{(i)} \in L^*_{q_0}(M), i = 1, 2, 3, 4$ and some choice of the symbols of $v^{(i)}$ at $\zeta^{(i)}$, the principal symbol of $\mcu^{(4)}$ is non-vanishing at $q_0$. After that, we use the properties of real analytic functions to finish the proof. We need four light-like vectors $\zeta^{(i)}$ at $q_0$ such that their sum is also a light like vector at $q_0$.  We first choose light-like vectors $\tilde \zeta^{(i)} \in L^*_{q_0}M$ which are linearly independent as following
\ba
\tilde \zeta^{(1)}=(1,0,1,0),\ \ \tilde \zeta^{(2)}=(1,0,0,1), \ \ 
\tilde \zeta^{(3)}=(-1, -1, 0,0),\ \ \tilde \zeta^{(4)}=(1,-1,0,0).
\ea
Then we let $\alpha_1 = 1, \alpha_2 = -1, \alpha_4 = \rho^{10}$ with $\rho$ a large parameter, and solve for $\alpha_3$ such that $\zeta = \sum_{i = 1}^4 \alpha_i\tilde \zeta^{(i)}$ is light-like. In particular, we find that $
\alpha_3 = -\ha \rho^{-10}.$ Therefore, the light-like vectors $\zeta^{(i)} = \alpha_i \tilde \zeta^{(i)}$ at $q_0$ which we shall work with are
\beqq\label{zetas}
\left.\begin{array}{ll}
\zeta^{(1)}=(1,0,1,0),\ \ &\zeta^{(2)}= -(1,0,0,1), \\
\zeta^{(3)}=  \ha \rho^{-10} ( 1,  1, 0,0),\ \ &\zeta^{(4)}= \rho^{10}(1,-1,0,0).
\end{array}\right.
\eeqq
We compute 
\beq
\left.\begin{array}{ll}
h(\zeta^{(1)},\zeta^{(2)})= 1, \ \ & h(\zeta^{(2)},\zeta^{(3)})= \ha \rho^{-10}, \\ 
h(\zeta^{(1)},\zeta^{(3)})=-\ha \rho^{-10}, & h(\zeta^{(2)},\zeta^{(4)})=\rho^{10}, \\
h(\zeta^{(1)},\zeta^{(4)})=-\rho^{10}, & h(\zeta^{(3)},\zeta^{(4)})=-1.
\end{array}\right.
\eeq
Also, we have 
\beq
\left.\begin{array}{ll}
 |\zeta^{(1)} + \zeta^{(2)}+ \zeta^{(3)}|_h^2 = 2,\ \  &|\zeta^{(1)} + \zeta^{(2)}+ \zeta^{(4)}|_h^2= 2,\\
 |\zeta^{(1)} + \zeta^{(3)}+ \zeta^{(4)}|_h^2= -2\rho^{10}-2+ O(\rho^{-10}),\ \  &|\zeta^{(2)} + \zeta^{(3)}+ \zeta^{(4)}|_h^2 = 2\rho^{10} -2 + O(\rho^{-10}).
\end{array}\right.
\eeq
These terms will appear many times in the computation of the principal symbols. 

Now consider the choice of the principal symbols of $v^{(i)}, i = 1, 2, 3, 4$. Recall that $v^{(i)}$ should satisfy the linearized gauge conditions
\beq
-h^{an} \p_a v^{(i)}_{nj} + \ha h^{pq} \p_j v^{(i)}_{pq} = 0, \ \ j = 0, 1, 2, 3.
\eeq
This is also called the polarization conditions. Let $A^{(i)}$ be the principal symbols $\sigma(v^{(i)})$, then they should satisfy the microlocal linearized gauge conditions
\beq
-h^{an} \zeta^{(i)}_a A^{(i)}_{nj} + \ha h^{pq} \zeta^{(i)}_j A^{(i)}_{pq} = 0, \ \ j = 0, 1, 2, 3.
\eeq
For $\zeta^{(i)}, i = 1, 2, 3, 4$, we choose
\beq
A_{ab}^{(i)} = \zeta_a^{(i)}\zeta^{(i)}_b, \ \ a, b = 0, 1,2, 3.
\eeq
It is straightforward to check that 
\beq
h^{ab}A^{(i)}_{ab} = 0, \ \ h^{mn}\zeta^{(i)}_nA^{(i)}_{mk} = 0.
\eeq
Thus the symbols satisfy the microlocal linearized gauge conditions. Also, it is easy to see that 
\beq
A^{(1)} = O(1), \ \ A^{(2)} = O(1),\ \ A^{(3)} = O(\rho^{-20}),\ \ A^{(4)} = O(\rho^{20}).
\eeq

\textbf{Step 4:} With these choices, we count the order of the principal symbols of terms in $\mcu^{(4)}$ in Prop.\ \ref{intert} as $\rho\rightarrow \infty$ to find the leading order terms. It suffices to consider the symbols of $\mch = \sum_{i = 1}^5 \mch_{i}$ in Prop.\ \ref{intert} at $q_0$. Below, we use $\bullet$ to denote a generic index of a vector or a two tensor.
\begin{enumerate}[(i)]
\item  $\sigma(\mch_1)$: Observe that these terms have two derivatives in $v^{\bullet}$. Therefore, the symbol of $\mch_{1}$ are linear combinations of terms like
\beq
\zeta^{(a)}_{\bullet} \zeta^{(b)}_\bullet \prod_{i = 1}^4A^{(i)}_\bullet, \ \ a, b = 1, 2, 3, 4.
\eeq
The $\zeta^{(a)}_{\bullet}$ and $\zeta^{(b)}_\bullet$ come from the two derivatives in $\mch_1$. Notice that the order of $ \prod_{i = 1}^4A^{(i)}_\bullet$ is $O(1)$. To achieve the leading order, we would like $a = b = 4$ to get $\rho^{20}$. This happens in the term $P_4(v^{(i)}, v^{(j)}, v^{(k)}, v^{(4)})$.  

\item  $\sigma(\mch_2)$: The symbols of these terms are linear combinations of 
\beq
\zeta^{(a)}_{\bullet} \zeta^{(b)}_\bullet \frac{\zeta^{(c)}_\bullet \zeta^{(d)}_\bullet}{|\zeta^{(k)}+ \zeta^{(l)}|^2} \prod_{i = 1}^4A^{(i)}_\bullet, \ \ a, b, c, d = 1, 2, 3, 4,
\eeq
where at least one of $k, l$ is equal to  $c$ or $d.$ To maximize the order in $\rho$, we would like $a = b = c = d = 4$. This implies that one of $k, l$ is $4$. The leading order is achieved if the other one of $k, l$ is $3$, in which case the order is $\rho^{40}\rho^{0} = \rho^{40}$. The corresponding terms are 
\beq
P_3(v^{(i)}, v^{(j)}, \bfq(P_2(v^{(3)}, v^{(4)}))) \text{ where $i, j$ are $1, 2.$}
\eeq

\item  $\sigma(\mch_3)$: These are  linear combinations of 
\beq
\zeta^{(a)}_{\bullet} \zeta^{(b)}_\bullet \frac{\zeta^{(c)}_\bullet \zeta^{(d)}_\bullet}{ |\zeta^{(j)} + \zeta^{(k)} + \zeta^{(l)}|^2} \prod_{i = 1}^4A^{(i)}_\bullet, \ \ a, b, c, d = 1, 2, 3, 4,
\eeq
where at least one of $j, k, l$ is  equal to $c$ or $d.$ 
To maximize the order, we take $a=b=c=d=4$ which implies that one of $j, k, l$ is $4$. The leading order is achieved if the other two of $j, k, l$ are $1, 2$ and the order is $\rho^{40}\rho^0 = \rho^{40}$. The corresponding terms are 
\beq
P_2(v^{(3)}, \bfq(P_3(v^{(j)}, v^{(k)}, v^{(4)})))  \text{ where $j, k$ are $1, 2.$}
\eeq

\item  $\sigma(\mch_4)$: The symbols are linear combinations of 
\beq
\zeta^{(a)}_{\bullet} \zeta^{(b)}_\bullet  \frac{\zeta^{(c)}_\bullet \zeta^{(d)}_\bullet}{|\zeta^{(i)}+ \zeta^{(j)}|^2} \frac{\zeta^{(e)}_\bullet \zeta^{(f)}_\bullet}{|\zeta^{(k)}+ \zeta^{(l)}|^2} \prod_{i = 1}^4A^{(i)}_\bullet, \ \ a, b, c, d = 1, 2, 3, 4,
\eeq
in which at least one of $i, j$ is equal to  $c$ or $d$ and one of $k, l$ is equal to  $e$ or $f$. 
In this case, the leading order terms can be achieved when four of $a, b, c, d, e, f$ are equal to $4$
 and the leading order is $\rho^{40}$. The corresponding terms are 
 \beq
 P_2(\bfq(\hat G_2(v^{(i)}, v^{(j)})), \bfq(P_2(v^{(3)}, v^{(4)}))) \text{ where $i, j$ are $1, 2$}.
 \eeq
Recall that $\hat G_2 = P_2 + \hat H_2$ where $\hat H_2$ is the quadratic derivative semilinear terms.

\item  $\sigma(\mch^5)$: The symbols are linear combinations of 
\beqq\label{symmch5}
\zeta^{(a)}_{\bullet} \zeta^{(b)}_\bullet \frac{\zeta^{(c)}_\bullet \zeta^{(d)}_\bullet}{ |\zeta^{(j)} + \zeta^{(k)} + \zeta^{(l)}|^2} \frac{\zeta^{(e)}_\bullet \zeta^{(f)}_\bullet}{|\zeta^{(k)}+ \zeta^{(l)}|^2} \prod_{i = 1}^4A^{(i)}_\bullet, \ \ a, b, c, d = 1, 2, 3, 4,
\eeqq
where at least one of $k, l$ is equal to one of $e, f$.  Consider the term
\beq
 P_2(x, v^{(i)}, \bfq( P_2(x, v^{(j)}, \bfq(P_2(x,  v^{(k)},  v^{(4)}))))).
\eeq
For example, if $k = 3$ and $i, j$ are $1, 2$, one can check that the max order of the term is $\rho^{60}\rho^{-10} = \rho^{50}$. Unfortunately, as shown below, such terms cancel each other. 

We compute the symbol of 
\beqq\label{eqlead}
 \mci \doteq \sum_{(i, j, k)\in \mca(3)} \sigma(P_2(x, v^{(i)}, \bfq( P_2(x, v^{(j)}, \bfq(P_2(x,  v^{(k)},  v^{(4)})))))),
\eeqq
where $\mca(3)$ denote the set of permutations of $(1, 2, 3)$. Observe that each component of the above terms are of the type $\mcy_4$ in \eqref{fourtho} so we can write down the symbol easily. Using the metric expansion, we find that  
\beq
P_2(x, u, v) = (HuH)^{pq} \p_p \p_q v, 
\eeq
for $u, v \in C^\infty(M; \sym^2)$ which will be identified below  in local coordinates as $4\times 4$ matrix valued. Then we have
\beq
\begin{gathered}
 \sigma(P_2(x, v^{(i)}, \bfq( P_2(x, v^{(j)}, \bfq(P_2(x,  v^{(k)},  v^{(4)}))))))(q_0, \zeta)\\
  = (2\pi)^{-3} (HA^{(i)}H)^{pq} (\zeta^{(4)}_p+\zeta^{(j)}_p+\zeta^{(k)}_p) (\zeta^{(4)}_q+\zeta^{(j)}_q+\zeta^{(k)}_q) \\
  \cdot \frac{ (H A^{(j)} H)^{mn} (\zeta^{(4)}_m+\zeta^{(k)}_m) (\zeta^{(4)}_n+ \zeta^{(k)}_n)}{|\zeta^{(j)} + \zeta^{(k)} + \zeta^{(4)}|_h^2}\cdot \frac{(H A^{(k)} H)^{ab} \zeta^{(4)}_a \zeta^{(4)}_b A^{(4)}}{ |\zeta^{(k)} + \zeta^{(4)}|_h^2}.
\end{gathered}
\eeq
We also observe that
\beqq\label{eqhh2sym}
\begin{split}
(H A^{(k)} H)^{pq} \zeta^{(4)}_p \zeta^{(4)}_q A^{(4)}= h^{pp}\zeta^{(k)}_p \zeta^{(k)}_q h^{qq} \zeta^{(4)}_p \zeta^{(4)}_q A^{(4)}  = [h(\zeta^{(k)}, \zeta^{(4)})]^2 A^{(4)}. 
\end{split}
\eeqq
So we get 
\beq
\begin{gathered}
 \sigma(P_2(x, v^{(i)}, \bfq( P_2(x, v^{(j)}, \bfq(P_2(x,  v^{(k)},  v^{(4)}))))))(q_0, \zeta)\\
  = (2\pi)^{-3} \frac{[h(\zeta^{(i)}, \zeta^{(4)} +\zeta^{(j)} +\zeta^{(k)}) h(\zeta^{(j)}, \zeta^{(4)}  +\zeta^{(k)})h(\zeta^{(k)}, \zeta^{(4)})]^2}{|\zeta^{(i)} + \zeta^{(j)} + \zeta^{(k)}|_h^2 |\zeta^{(k)} + \zeta^{(4)}|_h^2}.
\end{gathered}
\eeq
Let's compute the symbol for the six permutations in $\mci$. First we have 
\beq
\begin{split}
 &\textbf{(a): }\sigma(P_2(x, v^{(1)}, \bfq( P_2(x, v^{(2)}, \bfq(P_2(x,  v^{(3)},  v^{(4)}))))))(q_0, \zeta)\\
&= (2\pi)^{-3} \cdot  \frac{[h(\zeta^{(1)}, \zeta^{(4)} +\zeta^{(2)})h(\zeta^{(2)}, \zeta^{(4)})h(\zeta^{(3)}, \zeta^{(4)})]^2}{|\zeta^{(2)} + \zeta^{(3)} + \zeta^{(4)}|^2 |\zeta^{(3)} + \zeta^{(4)}|^2}(1 + O(\rho^{-20})) A^{(4)}\\
& =  (2\pi)^{-3}  \frac{(\rho^{10} - 1)^2 \rho^{20}}{(2\rho^{10}-2)(-2)}(1 + O(\rho^{-20})) A^{(4)} =  (2\pi)^{-3} (-\frac 14) (\rho^{30}-\rho^{20})A^{(4)}  + O(\rho^{30}) .
\end{split}
\eeq
Similarly, we get 
 \beq
\begin{split}
 &\textbf{(b): }\sigma(P_2(x, v^{(2)}, \bfq( P_2(x, v^{(1)}, \bfq(P_2(x,  v^{(3)},  v^{(4)}))))))(q_0, \zeta)\\
&=  \frac{[h(\zeta^{(2)}, \zeta^{(4)} +\zeta^{(1)})h(\zeta^{(1)}, \zeta^{(4)})h(\zeta^{(3)}, \zeta^{(4)})]^2}{|\zeta^{(1)} + \zeta^{(3)} + \zeta^{(4)}|^2 |\zeta^{(3)} + \zeta^{(4)}|^2}(1 + O(\rho^{-20}))A^{(4)}\\
 &=  (2\pi)^{-3}  \frac{(\rho^{10} + 1)^2 \rho^{20}}{(-2\rho^{10}-2)(-2)}(1 + O(\rho^{-20})) A^{(4)} =  (2\pi)^{-3} (\frac 14) (\rho^{30} + \rho^{20})A^{(4)} + O(\rho^{30}) .
\end{split}
\eeq

\beq
\begin{split}
 &\textbf{(c): }\sigma(P_2(x, v^{(1)}, \bfq( P_2(x, v^{(3)}, \bfq(P_2(x,  v^{(2)},  v^{(4)}))))))(q_0, \zeta)\\
&= \frac{[h(\zeta^{(1)}, \zeta^{(4)})h(\zeta^{(3)}, \zeta^{(4)})h(\zeta^{(2)}, \zeta^{(4)})]^2}{|\zeta^{(3)} + \zeta^{(2)} + \zeta^{(4)}|^2 |\zeta^{(2)} + \zeta^{(4)}|^2}  (1 + O(\rho^{-20}))A^{(4)}\\
 &=  (2\pi)^{-3}  \frac{(\rho^{20}  \rho^{20}}{( 2\rho^{10}-2)( 2 \rho^{10})}(1 + O(\rho^{-20})) A^{(4)} =  (2\pi)^{-3} (\frac 14)  \rho^{20} A^{(4)} + O(\rho^{30}) .
\end{split}
\eeq

\beq
\begin{split}
 &\textbf{(d): }\sigma(P_2(x, v^{(3)}, \bfq( P_2(x, v^{(1)}, \bfq(P_2(x,  v^{(2)},  v^{(4)}))))))(q_0, \zeta)\\
&=  \frac{[h(\zeta^{(3)}, \zeta^{(4)})h(\zeta^{(1)}, \zeta^{(4)} + \zeta^{(2)})h(\zeta^{(2)}, \zeta^{(4)})]^2}{|\zeta^{(1)} + \zeta^{(2)} + \zeta^{(4)}|^2 |\zeta^{(2)} + \zeta^{(4)}|^2}(1 + O(\rho^{-20}))A^{(4)}\\
 &=  (2\pi)^{-3}  \frac{(-\rho^{10} + 1)^2  \rho^{20}}{ 2 ( 2 \rho^{10})}(1 + O(\rho^{-20})) A^{(4)} =  (2\pi)^{-3} (\frac 14)  (\rho^{30} - 2\rho^{20}) A^{(4)}  + O(\rho^{30}) .
\end{split}
\eeq

\beq
\begin{split}
 &\textbf{(e): }\sigma(P_2(x, v^{(2)}, \bfq( P_2(x, v^{(3)}, \bfq(P_2(x,  v^{(1)},  v^{(4)}))))))(q_0, \zeta)\\
&=  \frac{[h(\zeta^{(2)}, \zeta^{(4)})h(\zeta^{(3)}, \zeta^{(4)})h(\zeta^{(1)}, \zeta^{(4)})]^2}{|\zeta^{(3)} + \zeta^{(1)} + \zeta^{(4)}|^2 |\zeta^{(1)} + \zeta^{(4)}|^2}  (1 + O(\rho^{-20}))A^{(4)}\\
 &=  (2\pi)^{-3}  \frac{ \rho^{20}   \rho^{20}}{  ( -2 \rho^{10} - 2) (-2\rho^{10})}(1 + O(\rho^{-20})) A^{(4)} =  (2\pi)^{-3} (\frac 14)   \rho^{20} A^{(4)}  + O(\rho^{30}) .
\end{split}
\eeq

\beq
\begin{split}
 &\textbf{(f): }\sigma(P_2(x, v^{(3)}, \bfq( P_2(x, v^{(2)}, \bfq(P_2(x,  v^{(1)},  v^{(4)}))))))(q_0, \zeta)\\
&=   \frac{[h(\zeta^{(3)}, \zeta^{(4)})h(\zeta^{(2)}, \zeta^{(4)} + \zeta^{(1)})h(\zeta^{(1)}, \zeta^{(4)})]^2}{|\zeta^{(2)} + \zeta^{(1)} + \zeta^{(4)}|^2 |\zeta^{(1)} + \zeta^{(4)}|^2}(1 + O(\rho^{-20}))A^{(4)}\\
 &=  (2\pi)^{-3}  \frac{ (\rho^{10}+1)^2   \rho^{20}}{  2 (-2\rho^{10})}(1 + O(\rho^{-20})) A^{(4)} =  (2\pi)^{-3} (-\frac 14)  (\rho^{30} + 2\rho^{20}) A^{(4)}  + O(\rho^{30}) .
\end{split}
\eeq

Finally, summing up the above terms (a) to (f), we arrive at $\sigma(\mci)(q_0, \zeta)  = O(\rho^{30})$, which is lower than expected. 
\begin{remark}\label{symcan}
The cancellation of leading order terms is not a coincidence and it happens in later calculations. This phenomena perhaps is related to the null structure of the reduced Einstein equations. In a simpler setting of scalar wave equations with quadratic derivative nonlinear terms, it is showed in \cite{WZ}  that the null terms do not contribute to the leading order singularities in the interaction term, which is regarded as a manifestation of the smoothing effects. In view of \eqref{eqhh2sym}, we heuristically think of $P_2$ as a null form for our choice of $v^{(i)}$. (This is part of the so called generalized null condition, see \cite[Section XI.5]{Cb}). $P_3$ is similar in view of \eqref{eqhh3sym} below. In fact, as computed explicitly in the following, the symbol of the leading order terms in $\mch$ we look for involving only $P_2, P_3$ cancel each other. Fortunately, the quadratic derivative terms in the reduced Einstein equations are not all null forms, see e.g.\ \cite{LR}, and we  expect to (and do) get a non-trivial leading order term.  A thorough understanding of this phenomena could greatly simplify the symbol calculation, however, we leave it for future considerations. 
\end{remark}

We conclude that the leading order for the symbol of $\mch_5$ is at most $\rho^{40}.$ Thus we need to determine all possible $\rho^{40}$ order terms in $\sigma(\mch_5)$. The case we analyzed above is when all the indices  $a, b, c, d, e, f$ in \eqref{symmch5} equals to $4$ and the term $\mci$ include all the possible $\rho^{40}$ order terms. Next, consider the case when five of these indices are $4$. This implies that one of $\hat G_2$ in $\mch_5$ must be $\hat H_2$. Roughly speaking, there are three cases to consider depending on the position of $\hat H_2$. 
\begin{enumerate}
\item Suppose the first $\hat G_2$ in $\mch_5$ is $\hat H_2$. Then the leading order term must be among 
\beq
\hat H_2(x, v^{(i)}, \bfq( P_2(x, v^{(j)}, \bfq(P_2(x,  v^{(k)},  v^{(4)}))))).
\eeq
 We remark that the symmetric term  $\hat H_2(x, \bfq( P_2(x, v^{(j)}, \bfq(P_2(x,  v^{(k)},  v^{(4)})))), v^{(i)})$ should be considered as well. For simplicity, here and several places below, we just show one such term and leave the symmetric term  in the actual calculation, and we note that this is only necessary when $\hat H_i$ terms are present. If $i = 3$ or $j = 3$, we find that the order is at most $\rho^{30}$. So we must have $k = 3$ to get the $\rho^{40}$ order term. 
 
 \item Suppose the middle $\hat G_2$ in $\mch_5$ is $\hat H_2$. Then the leading order term must be among 
\beq
P_2(x, v^{(i)}, \bfq( \hat H_2(x, v^{(j)}, \bfq(P_2(x,  v^{(k)},  v^{(4)}))))).
\eeq
Here, again the symmetric term within $\hat H_2$ should be considered. If $i = 3$ or $k = 3$, it is possible to obtain $\rho^{40}$. However, if $j = 3$, the maximal order is $\rho^{30}$. 

\item Suppose the last $\hat G_2$ in $\mch_5$ is $\hat H_2$. 
Then the leading order term must be among 
\beq
\begin{gathered}
P_2(x, v^{(i)}, \bfq( P_2(x, v^{(j)}, \bfq(\hat H_2(x,  v^{(k)},  v^{(4)}))))) .
\end{gathered}
\eeq
We find that $i$ must be $3$ to achieve order $\rho^{40}$.  \\
\end{enumerate}

To get the $O(\rho^{40})$ terms, the last possibility is that  four indices of $a, b, c, d, e, f$ are $4$. If all the $\hat G_2$ in $\mch_5$ are $P_2$,  there is one possible term 
\beq
 P_2(x, v^{(3)}, \bfq( P_2(x, \bfq(P_2(x,  v^{(k)},  v^{(l)})), v^{(4)})))  \text{ where $k, l$ are $1, 2$}.
\eeq
If one of the $\hat G_2$ in $\mch_5$ is $\hat H_2$, we should consider 
\beq
P_2(x, v^{(i)}, \bfq( P_2(x,  \bfq(\hat H_2(x,  v^{(k)},  v^{(l)}))), v^{(4)})) \text{ where $k, l$ are $1, 2$}.
\eeq
At last, if two of  the $\hat G_2$ in $\mch_5$ are $\hat H_2$, again there are three cases to consider depending on the position of $P_2$: 

\begin{enumerate}
\item Suppose the first $\hat G_2$ in $\mch_5$ is $P_2$. Then the leading order term must be among 
\beq
P_2(x, v^{(i)}, \bfq( \hat H_2(x, v^{(j)}, \bfq(\hat H_2(x,  v^{(k)},  v^{(4)}))))).
\eeq
However, one can check that all the possible combinations are at most of order $\rho^{30}$. So these terms do not contribute to the leading order terms. 
 
 \item Suppose the middle $\hat G_2$ in $\mch_5$ is $P_2$. Then the leading order term must be among 
\beq
\hat H_2(x, v^{(i)}, \bfq( P_2(x, v^{(j)}, \bfq(\hat H_2(x,  v^{(k)},  v^{(4)}))))).
\eeq
Similarly, one find that no such terms are of the order $\rho^{40}$. 

\item Suppose the last $\hat G_2$ in $\mch_5$ is $P_2$. 
Then the leading order term must be among 
\beq
\hat H_2(x, v^{(i)}, \bfq( \hat H_2(x, v^{(j)}, \bfq(P_2(x,  v^{(k)},  v^{(4)}))))).
\eeq
Again, no such term is of order $\rho^{40}$.\\
\end{enumerate}

\end{enumerate}

To summarize, the leading order of the symbol of $\mch$ as $\rho\rightarrow \infty$ is $\rho^{40}$. This order could be achieved among the following terms (minus sign is added appropriately in accordance with $\mch_i$ in Prop.\ \ref{intert}):
\begin{enumerate}
\item $P_3(v^{(i)}, v^{(j)}, \bfq(P_2(v^{(3)}, v^{(4)}))) \text{ where $i, j$ are $1, 2.$} $
\item $P_2(v^{(3)}, \bfq(P_3(v^{(j)}, v^{(k)}, v^{(4)})))  \text{ where $j, k$ are $1, 2.$} $
\item $-P_2(\bfq(\hat G_2(v^{(i)}, v^{(j)})), \bfq(P_2(v^{(3)}, v^{(4)}))) \text{ where $i, j$ are $1, 2$}.$
\item term $-\mci$ in \eqref{eqlead}.
\item $-\hat H_2(x, v^{(i)}, \bfq( P_2(x, v^{(j)}, \bfq(P_2(x,  v^{(3)},  v^{(4)}))))) \text{ where $i, j$ are $1,2$}$.
\item $-P_2(x, v^{(i)}, \bfq( \hat H_2(x, v^{(j)}, \bfq(P_2(x,  v^{(k)},  v^{(4)}))))) \text{ where $i$ or $k$ is $3$.}$
\item $-P_2(x, v^{(3)}, \bfq( P_2(x, v^{(j)}, \bfq(\hat H_2(x,  v^{(k)},  v^{(4)}))))) \text{ where $j, k$ are $1,2$}$.
\item $- P_2(x, v^{(3)}, \bfq( P_2(x, \bfq(\hat G_2(x,  v^{(k)},  v^{(l)})), v^{(4)}))) \text{ where $k, l$ are $1,2$}$. 
\end{enumerate}

\textbf{Step 5:} We compute the symbols of the terms listed above. 

\textbf{(1) and (2): }For these two terms, we only need the expressions of $P_2$ which is already found, and $P_3$. For $u, v, w\in C^\infty(M; \sym^2)$, using the metric expansion \eqref{eqgexpan}, we find that 
\beq
P_3(u, v, w) = -(HuHvH)^{pq} \p_p \p_q w. 
\eeq
Therefore, using the symbol expressions in \eqref{fourtho} , we get
\beq
\begin{split}
&\sigma(P_3(v^{(i)}, v^{(j)}, \bfq(P_2(v^{(3)}, v^{(4)}))))(q_0, \zeta) \\
= &-(2\pi)^{-3}  (HA^{(i)}HA^{(j)}H)^{pq}\zeta^{(4)}_p \zeta^{(4)}_q \cdot \dfrac{(HA^{(3)}H)^{ab} \zeta^{(4)}_a \zeta^{(4)}_b A^{(4)}}{|\zeta^{(3)} + \zeta^{(4)}|^2}(1 + O(\rho^{-10}))\\
= &-(2\pi)^{-3}  (HA^{(i)}HA^{(j)}H)^{pq}\zeta^{(4)}_p \zeta^{(4)}_q \cdot  \ha h(\zeta^{(3)}, \zeta^{(4)}) A^{(4)} (1 + O(\rho^{-10})).
\end{split}
\eeq
Now we observe that 
\beqq\label{eqhh3sym}
\begin{gathered}
(HA^{(i)}HA^{(j)}H)^{pq}\zeta^{(4)}_p \zeta^{(4)}_q = h^{pp}\zeta^{(i)}_p\zeta^{(i)}_ah^{aa} \zeta^{(j)}_a \zeta^{(j)}_q h^{qq} \zeta^{(4)}_p \zeta^{(4)}_q \\
 = h(\zeta^{(i)}, \zeta^{(4)}) h(\zeta^{(j)}, \zeta^{(4)}) h(\zeta^{(i)}, \zeta^{(j)}). 
\end{gathered}
\eeqq
For $i, j = 1, 2$, we find that the above number is $-\rho^{20}$. Therefore, we complete the calculate of the symbol of the term in (1):
\beq
\begin{split}
&\sum_{i, j = 1, 2, i\neq j} \sigma(P_3(v^{(i)}, v^{(j)}, \bfq(P_2(v^{(3)}, v^{(4)}))))(q_0, \zeta) \\
= &-2 (2\pi)^{-3} (-\rho^{20}) \cdot  \ha (-1)A^{(4)} (1 + O(\rho^{-10})) = -(2\pi)^{-3}  \rho^{20} A^{(4)}  + O(\rho^{30}).
\end{split}
\eeq

For terms in (2), the calculation is similar. We have 
\beq
\begin{split}
& \sum_{j, k = 1, 2, j\neq k}\sigma(P_2(v^{(3)}, \bfq(P_3(v^{(j)}, v^{(k)}, v^{(4)}))) )(q_0, \zeta) \\
= &-(2\pi)^{-3}   \sum_{j, k = 1, 2, j\neq k}(HA^{(3)}H)^{pq} \zeta^{(4)}_p \zeta^{(4)}_q  \dfrac {(HA^{(j)}HA^{(k)}H)^{ab}\zeta^{(4)}_a \zeta^{(4)}_b A^{(4)}}{|\zeta^{(j)} + \zeta^{(k)} + \zeta^{(4)}|^2} (1 + O(\rho^{-10}))\\
= &-(2\pi)^{-3}  \sum_{j, k = 1, 2, j\neq k}[h(\zeta^{(3)}, \zeta^{(4)})]^2   \frac{h(\zeta^{(k)}, \zeta^{(4)}) h(\zeta^{(j)}, \zeta^{(4)}) h(\zeta^{(k)}, \zeta^{(j)})}{|\zeta^{(j)} + \zeta^{(k)} + \zeta^{(4)}|^2}A^{(4)} (1 + O(\rho^{-10}))\\
= & -2(2\pi)^{-3}  \frac{h(\zeta^{(1)}, \zeta^{(4)}) h(\zeta^{(2)}, \zeta^{(4)}) h(\zeta^{(1)}, \zeta^{(2)})}{|\zeta^{(1)} + \zeta^{(2)} + \zeta^{(4)}|^2}A^{(4)} (1 + O(\rho^{-10})) \\
 = &(2\pi)^{-3} \rho^{20} A^{(4)} + O(\rho^{30}). 
\end{split}
\eeq
We see that the sum of the symbols of (1) and (2) is of order $\rho^{30}$, so they do not contribute to $O(\rho^{40})$.  \\

\textbf{(3): }For this term, if the $\hat G_2$ is $P_2$, the calculation is simple. Actually, we have 
\beq
\begin{split}
&\sigma(-P_2(\bfq(P_2(v^{(i)}, v^{(j)})), \bfq(P_2(v^{(3)}, v^{(4)}))))(q_0, \zeta) \\
= &(2\pi)^{-3}  (H  \frac{(HA^{(i)}H)^{ab}\zeta^{(j)}_a\zeta^{(j)}_b A^{(j)}}{|\zeta^{(i)} + \zeta^{(j)}|^2} H)^{pq} \zeta^{(4)}_p \zeta^{(4)}_q  \frac {(HA^{(3)}H)^{cd}\zeta^{(4)}_c \zeta^{(4)}_d A^{(4)}}{|\zeta^{(3)} + \zeta^{(4)}|^2} (1 + O(\rho^{-10}))\\
= &(2\pi)^{-3} \cdot  h^2(\zeta^{(j)}, \zeta^{(4)}) \cdot \ha h(\zeta^{(i)}, \zeta^{(j)}) \cdot \ha h(\zeta^{(3)}, \zeta^{(4)}) A^{(4)} (1 + O(\rho^{-10})).
\end{split}
\eeq
When $i, j$ are $1, 2$, the symbol of the terms in (3) is 
\beq
\begin{split}
&\sum_{i, j = 1, 2, i\neq j} \sigma(-P_2(\bfq(P_2(v^{(i)}, v^{(j)})), \bfq(P_2(v^{(3)}, v^{(4)}))))(q_0, \zeta)\\
&=  (2\pi)^{-3} \frac {1}{4}   [h^2(\zeta^{(1)}, \zeta^{(4)}) + h^2(\zeta^{(2)}, \zeta^{(4)})]  h(\zeta^{(1)}, \zeta^{(2)})  h(\zeta^{(3)}, \zeta^{(4)}) A^{(4)} (1 + O(\rho^{-10})) \\
&=  (2\pi)^{-3}(-\ha )\rho^{20} A^{(4)} + O(\rho^{30}). 
\end{split}
\eeq

Next, we need $\hat H_2$ to compute the rest. Recall that the $\hat H_2$ term is the semilinear term in the (reduced) Einstein equations and has two derivatives. Using \eqref{eqricre}, we see that the semilinear terms in our formulation are 
\beq 
\begin{gathered}
H(u)  = 2 g^{ab}g_{ps}\Gamma^p_{\mu b}\Gamma^s_{\nu a} +  (g_{\nu l} \Gamma^l_{ab} g^{aq}g^{bd} \frac{\p u_{qd}}{\p x^\mu} + g_{\mu l} \Gamma^l_{ab} g^{aq} g^{bd} \frac{\p u_{qd}}{\p x^\nu}),
\end{gathered}
\eeq
in which the Christoffel symbols are
\beqq\label{eqmcg}
\Gamma_{\alpha\beta}^\mu = \ha g^{\mu\la}(\frac{\p u_{\la\alpha}}{\p x^\beta} + \frac{\p u_{\la\beta}}{\p x^\alpha}- \frac{\p u_{\alpha\beta}}{\p x^\la}) =  g^{\mu\la} \mcg_{\la\alpha\beta}(u)
\eeqq
where $\mcg(u)$ is defined by the above equation. Then we can determine the quadratic term $\hat H_2$ as
\beq 
\begin{gathered}
\hat H_{2, \mu\nu} (u, u)  = 2 h^{ab}h_{ps}h^{p\la}\mcg(u)_{\la \mu b}h^{s\gamma}\mcg(u)_{\gamma \nu a} +  (h_{\nu l} h^{l\la}\mcg(u)_{\la ab} h^{aq}h^{bd} \frac{\p u_{qd}}{\p x^\mu} + h_{\mu l} h^{l\gamma}\mcg(u)_{\gamma ab} h^{aq} h^{bd} \frac{\p u_{qd}}{\p x^\nu}).
\end{gathered}
\eeq

Using this symbol expression in \eqref{eqpsym}, we find that 
\beq
\begin{split}
&\sum_{i, j = 1, 2, i\neq j} \sigma(-P_2(\bfq(\hat H_2(v^{(i)}, v^{(j)})), \bfq(P_2(v^{(3)}, v^{(4)}))))(q_0, \zeta) \\
= & (2\pi)^{-3} \sum_{i, j = 1, 2, i\neq j}  (H C^{(ij)} H)^{pq} \zeta^{(4)}_p \zeta^{(4)}_q  \frac {(HA^{(3)}H)^{cd}\zeta^{(4)}_c \zeta^{(4)}_d A^{(4)}}{|\zeta^{(3)} + \zeta^{(4)}|^2} (1 + O(\rho^{-10})),
\end{split}
\eeq
where $C^{(ij)} = B^{(ij)} + B^{(ji)}$ is a matrix and 
\beq 
\begin{split}
B^{(ij)} \doteq  &  2 h^{ab}h_{ps}h^{p\la}\ha \zeta^{(i)}_\la\zeta^{(i)}_\mu \zeta^{(i)}_b h^{s\gamma} \ha \zeta^{(j)}_\gamma \zeta^{(j)}_\nu \zeta^{(j)}_a  
+ (h_{\nu l} h^{l\la}\ha \zeta^{(i)}_\la\zeta^{(i)}_a \zeta^{(i)}_b h^{aq}h^{bd} \zeta^{(j)}_q \zeta^{(j)}_d \zeta^{(j)}_\mu \\
&+ h_{\mu l} h^{l\gamma}\ha \zeta^{(i)}_\gamma \zeta^{(i)}_a\zeta^{(i)}_b h^{aq} h^{bd}\zeta^{(j)}_q \zeta^{(j)}_d \zeta^{(j)}_\nu)\\
=& \ha [h(\zeta^{(i)}, \zeta^{(j)})]^2\zeta^{(i)}_\mu \zeta^{(j)}_\nu + \ha [h(\zeta^{(i)}, \zeta^{(j)})]^2\zeta^{(i)}_\nu \zeta^{(j)}_\mu + \ha [h(\zeta^{(i)}\zeta^{(j)})]^2\zeta^{(i)}_\mu \zeta^{(j)}_\nu.
\end{split}
\eeq
%
In particular, $B^{(ij)}$ correspond to the term $\hat H^{(2)}(v^{(i)}, v^{(j)}).$ The minus sign is gone because we have six derivatives which gives $(\imath)^6 = -1$. It is easy to see that 
\beq 
\begin{split}
C^{(ij)}  = \frac{3}{2} [h(\zeta^{(i)}, \zeta^{(j)})]^2[\zeta^{(i)}_\mu \zeta^{(j)}_\nu + \zeta^{(i)}_\nu \zeta^{(j)}_\mu].
\end{split}
\eeq
For later reference, we define a matrix $A^{(ij)}$ for $1\leq i< j\leq 4$ so that $A^{(ij)}_{\mu\nu} = \zeta^{(i)}_\mu \zeta^{(j)}_\nu + \zeta^{(i)}_\nu \zeta^{(j)}_\mu. $
Using these pre-calculation, we  find that 
\beq
\begin{split}
&\sum_{i, j = 1, 2, i\neq j} \sigma(-P_2(\bfq(\hat H_2(v^{(i)}, v^{(j)})), \bfq(P_2(v^{(3)}, v^{(4)}))))(q_0, \zeta) \\
= &(2\pi)^{-3}  \frac{3}{2} h(\zeta^{(1)}, \zeta^{(2)})\cdot (HA^{(12)}H)^{pq} \zeta^{(4)}_p \zeta^{(4)}_q  \cdot \ha h(\zeta^{(3)}, \zeta^{(4)}) A^{(4)} (1 + O(\rho^{-10}))\\
=& (2\pi)^{-3}  \frac{3}{4} h(\zeta^{(1)}, \zeta^{(2)}) [2 h(\zeta^{(1)}, \zeta^{(4)})  h(\zeta^{(2)}, \zeta^{(4)})] h(\zeta^{(3)}, \zeta^{(4)}) A^{(4)} (1 + O(\rho^{-10}))\\
=& (2\pi)^{-3} \frac 32 \rho^{20}A^{(4)} + O(\rho^{30}). 
\end{split}
\eeq
This completes the calculation of the symbol in (3). \\

\textbf{(7) and (8): }These two terms can be calculated using the calculation result we have so far. For (7), we see that the symbol is 
\beq
\begin{split}
&\sum_{k, j = 1, 2, k\neq j}\sigma(-P_2(x, v^{(3)}, \bfq( P_2(x, v^{(j)}, \bfq(\hat H_2(x,  v^{(k)},  v^{(4)}) + \hat H_2(x,  v^{(4)},  v^{(k)}))))))(q_0, \zeta) \\
= &\sum_{k, j = 1, 2, k\neq j} (2\pi)^{-3}  (H  A^{(3)}H)^{cd} \zeta^{(4)}_c \zeta^{(4)}_d  \cdot   \dfrac{(HA^{(j)}H)^{pq} \zeta^{(4)}_p \zeta^{(4)}_q  }{|\zeta^{(j)} + \zeta^{(k)} + \zeta^{(4)}|^2} \cdot   \frac{\frac{3}{2} h^2(\zeta^{(k)}, \zeta^{(4)})}{2h(\zeta^{(k)}, \zeta^{(4)})} A^{(k4)}(1 + O(\rho^{-10}))\\
= & (2\pi)^{-3} \cdot  h^2(\zeta^{(3)}, \zeta^{(4)}) \cdot  \frac 38 [h^2(\zeta^{(1)}, \zeta^{(4)})h(\zeta^{(2)}, \zeta^{(4)}) A^{(24)} + h^2(\zeta^{(2)}, \zeta^{(4)})h(\zeta^{(1)}, \zeta^{(4)}) A^{(14)}]  (1 + O(\rho^{-10}))\\
=& \frac 38  (2\pi)^{-3} \rho^{30} [A^{(24)}-A^{(14)}] + O(\rho^{30}). 
\end{split}
\eeq

For (8),  we first take the $\hat G_2$ to be $P_2$. Then we get 
\beq
\begin{split}
&\sum_{k, l = 1, 2, k\neq l}\sigma(- P_2(x, v^{(3)}, \bfq( P_2(x, \bfq(P_2(x,  v^{(k)},  v^{(l)})), v^{(4)})))) \\
= &\sum_{k, l = 1, 2, k\neq l} (2\pi)^{-3}  (H  A^{(3)}H)^{cd} \zeta^{(4)}_c \zeta^{(4)}_d  (H\cdot \frac{(HA^{(k)}H)^{pq} \zeta^{(l)}_p \zeta^{(l)}_q A^{(l)}}{|\zeta^{(l)} + \zeta^{(k)}|^2}\cdot H)^{ab}  \zeta^{(4)}_a \zeta^{(4)}_b A^{(4)} (1 + O(\rho^{-10}))\\
= & (2\pi)^{-3} \cdot  h^2(\zeta^{(3)}, \zeta^{(4)}) \cdot \ha h(\zeta^{(1)}, \zeta^{(2)}) \cdot \ha [h^2(\zeta^{(1)}, \zeta^{(4)})+ h^2(\zeta^{(2)}, \zeta^{(4)})] A^{(4)} (1 + O(\rho^{-10}))\\
=& \ha  (2\pi)^{-3} \rho^{20} A^{(4)} + O(\rho^{30}). 
\end{split}
\eeq
When the $\hat G_2$ is $\hat H_2$ in (8), we use the calculation in (3) to get 
\beq
\begin{split}
&\sum_{k, l = 1, 2, k\neq l}\sigma(- P_2(x, v^{(3)}, \bfq( P_2(x, \bfq(\hat H_2(x,  v^{(k)},  v^{(l)})), v^{(4)}))))(q_0, \zeta) \\
= &\sum_{k, l = 1, 2, k\neq l} (2\pi)^{-3}  (H  A^{(3)}H)^{cd} \zeta^{(4)}_c \zeta^{(4)}_d  (H\cdot \frac{ \frac 32 [h(\zeta^{(k)}, \zeta^{(l)})]^2 A^{(kl)}}{|\zeta^{(l)} + \zeta^{(k)}|^2}\cdot H)^{ab}  \zeta^{(4)}_a \zeta^{(4)}_b A^{(4)} (1 + O(\rho^{-10}))\\
= & (2\pi)^{-3} \cdot  h^2(\zeta^{(3)}, \zeta^{(4)}) \cdot \frac 34 h(\zeta^{(1)}, \zeta^{(2)}) \cdot [2 h(\zeta^{(1)}, \zeta^{(4)}) h(\zeta^{(2)}, \zeta^{(4)})] A^{(4)} (1 + O(\rho^{-10}))\\
=&  (2\pi)^{-3}(-\frac 32) \rho^{20} A^{(4)} + O(\rho^{30}). 
\end{split}
\eeq

\textbf{(5): }For this term, we use the symbol expression \eqref{eqpsym}. To see the structure clearly, we first observe that in the symbol of (5), the matrix corresponding to $\bfq( P_2(x, v^{(j)}, \bfq(P_2(x,  v^{(3)},  v^{(4)}))))$ is given by 
\beq
\begin{split}
 D^{(j)} \doteq & \frac{(H A^{(j)} H)^{pq} \zeta^{(4)}_p \zeta^{(4)}_q }{|\zeta^{(j)} + \zeta^{(3)} + \zeta^{(4)}|^2} \cdot \frac{(H A^{(3)} H)^{ab} \zeta^{(4)}_a \zeta^{(4)}_b }{|\zeta^{(3)} + \zeta^{(4)}|^2} \cdot A^{(4)}(1+O(\rho^{-10}))\\
 =&\frac{h^2(\zeta^{(j)}, \zeta^{(4)}) }{|\zeta^{(j)} + \zeta^{(3)} + \zeta^{(4)}|^2} \cdot  \ha h(\zeta^{(3)}, \zeta^{(4)})  \cdot A^{(4)}(1+O(\rho^{-10})).
\end{split}
\eeq
This matrix is a constant multiple of $A^{(4)}$. So the symbol calculation  is similar those for (3), see in particular the matrix $B^{(ij)}$.
Therefore,  we get 
\beq
\begin{split}
&\sum_{i, j = 1, 2, i\neq j}[\sigma(-\hat H_2(x, v^{(i)}, \bfq( P_2(x, v^{(j)}, \bfq(P_2(x,  v^{(3)},  v^{(4)}))))))(q_0, \zeta) \\
&+ \sigma(-\hat H_2(x, \bfq( P_2(x, v^{(j)}, \bfq(P_2(x,  v^{(3)},  v^{(4)})))),  v^{(i)}))(q_0, \zeta) ]\\
= &\sum_{i, j = 1, 2, i\neq j} (2\pi)^{-3} \frac{h^2(\zeta^{(j)}, \zeta^{(4)}) }{|\zeta^{(j)} + \zeta^{(3)} + \zeta^{(4)}|^2} \cdot  \ha h(\zeta^{(3)}, \zeta^{(4)}) \cdot \frac 32 [h(\zeta^{(i)}, \zeta^{(4)})]^2 A^{(i4)}(1+O(\rho^{-10}))\\
=& (2\pi)^{-3} \rho^{40} (-\frac 34)[\frac{1}{|\zeta^{(2)} + \zeta^{(3)} + \zeta^{(4)}|^2}A^{(14)} + \frac{1}{|\zeta^{(1)} + \zeta^{(3)} + \zeta^{(4)}|^2}A^{(24)}](1+O(\rho^{-10}))\\
=&(2\pi)^{-3} \rho^{30} (-\frac 38)[ A^{(14)} - A^{(24)}](1+O(\rho^{-10})).\\
\end{split}
\eeq
 
 \textbf{(6): }Using \eqref{eqpsym}, we see that 
  \beq
\begin{split}
&\sum_{i, j = 1, 2, i\neq j}[\sigma(-P_2(x, v^{(i)}, \bfq( \hat H_2(x, v^{(j)}, \bfq(P_2(x,  v^{(k)},  v^{(4)})))))) \\
& + \sigma(-P_2(x, v^{(i)}, \bfq( \hat H_2(x, \bfq(P_2(x,  v^{(k)},  v^{(4)})), v^{(j)})))) ] \\
= &\sum_{i, j = 1, 2, i\neq j} (2\pi)^{-3}  (H  A^{(i)}H)^{cd} \zeta^{(4)}_c \zeta^{(4)}_d   \cdot  \frac{ E^{(kj)}}{|\zeta^{(j)} + \zeta^{(k)} + \zeta^{(4)}|^2}  (1 + O(\rho^{-10})),
\end{split}
\eeq
where $E^{(kj)} = \ha h(\zeta^{(k)}, \zeta^{(4)}) \frac 32  h^2 (\zeta^{(j)}, \zeta^{(4)}) A^{(j4)}$. Here, we used the same trick as in (5) by noticing that the matrix corresponding to $P_2(x,  v^{(k)},  v^{(4)})$ is $h^2(\zeta^{(k)}, \zeta^{(4)}) A^{(4)}$ which is a constant multiple of $A^{(4)}$.  Now we compute the symbol of (6) when $k = 3$: 
 \beq
\begin{split}
&\sum_{i, j = 1, 2, i\neq j}[\sigma(-P_2(x, v^{(i)}, \bfq( \hat H_2(x, v^{(j)}, \bfq(P_2(x,  v^{(3)},  v^{(4)})))))) \\
& + \sigma(-P_2(x, v^{(i)}, \bfq( \hat H_2(x, \bfq(P_2(x,  v^{(3)},  v^{(4)})), v^{(j)}))))] \\
= &\sum_{i, j = 1, 2, i\neq j} (2\pi)^{-3}  (H  A^{(i)}H)^{cd} \zeta^{(4)}_c \zeta^{(4)}_d   \cdot  \frac{ \ha h(\zeta^{(3)}, \zeta^{(4)}) \frac 32  h^2 (\zeta^{(j)}, \zeta^{(4)}) A^{(j4)}}{|\zeta^{(j)} + \zeta^{(3)} + \zeta^{(4)}|^2}  (1 + O(\rho^{-10}))\\
= & (2\pi)^{-3} \frac 34  h^2(\zeta^{(1)}, \zeta^{(4)}) h(\zeta^{(3)}, \zeta^{(4)})  h^2 (\zeta^{(2)}, \zeta^{(4)})  [\frac{1}{|\zeta^{(1)} + \zeta^{(3)} + \zeta^{(4)}|^2}A^{(14)} \\
&+ \frac{1}{|\zeta^{(2)} + \zeta^{(3)} + \zeta^{(4)}|^2}A^{(24)}] (1 + O(\rho^{-10})) 
=   (2\pi)^{-3}\frac 34 \rho^{30} [A^{(14)} - A^{(24)}] + O(\rho^{30}). 
\end{split}
\eeq
Finally, when $i = 3$ in (6), we get 
 \beq
\begin{split}
&\sum_{k, j = 1, 2, k\neq j}[\sigma(-P_2(x, v^{(3)}, \bfq( \hat H_2(x, v^{(j)}, \bfq(P_2(x,  v^{(k)},  v^{(4)}))))))\\
&+ \sigma(-P_2(x, v^{(3)}, \bfq( \hat H_2(x, \bfq(P_2(x,  v^{(k)},  v^{(4)})), v^{(j)}))))] \\
= &\sum_{k, j = 1, 2, k\neq j} (2\pi)^{-3}  (H  A^{(3)}H)^{cd} \zeta^{(4)}_c \zeta^{(4)}_d   \cdot  \frac{ \frac 32  \ha h(\zeta^{(k)}, \zeta^{(4)})  h^2 (\zeta^{(j)}, \zeta^{(4)}) A^{(j4)}}{|\zeta^{(j)} + \zeta^{(k)} + \zeta^{(4)}|^2}  (1 + O(\rho^{-10}))\\
= & (2\pi)^{-3} \frac 34  \frac{h^2(\zeta^{(3)}, \zeta^{(4)})} {|\zeta^{(1)} + \zeta^{(2)} + \zeta^{(4)}|^2}[h(\zeta^{(1)}, \zeta^{(4)})  h^2 (\zeta^{(2)}, \zeta^{(4)}) A^{(24)} + h(\zeta^{(2)}, \zeta^{(4)})  h^2 (\zeta^{(1)}, \zeta^{(4)}) A^{(14)}] (1 + O(\rho^{-10}))\\
=&  (2\pi)^{-3} \frac 38 \rho^{30} [A^{(14)} - A^{(24)}] + O(\rho^{30}). 
\end{split}
\eeq

To sum up, we showed that the symbol of the sum of (1)-(8) is given by 
\beq
\begin{gathered}
 (2\pi)^{-3} \frac 38 \rho^{30}  [A^{(14)} - A^{(24)}] = (2\pi)^{-3} \rho^{40} \begin{pmatrix}
 0& 0 &\frac 38 &-\frac 38\\[3pt]
 0 &0 &-\frac 38 &\frac 38 \\[3pt]
\frac 38 &-\frac 38&0 &0 \\[3pt]
-\frac 38& \frac 38& 0& 0 
 \end{pmatrix} + O(\rho^{30})
 \end{gathered}
\eeq
which is non-vanishing for $\rho$ large. 
 
\textbf{Step 6:} We showed that for the choice of our $\zeta^{(i)}$ and $A^{(i)}$, the principal symbol of $\mcu^{(4)}$ is non-vanishing. Since the symbol is a real analytic function, it cannot vanish on any relatively open subset of $\mcx$. This completes the proof of the proposition.
\epf

\section{The inverse problem in wave gauge}\label{sec-wave}
In this section, we prove Theorem \ref{main}  when the observation data is made in the wave gauge. From  the previous sections, we know that the wave gauge is convenient for mathematical analysis and we want to give a clear argument for the inverse problem in this case first, before we turn to the physical gauge in the next section. 

Let's recall that for $g'$ close to $\hat g$ in $C^4(M(\bt_0'))$, we let $f: (M(\bt_0'), g') \rightarrow (M(\bt_0), \hat g)$ be the wave map. Then $g = f_*g'$ is the image of the metric $g'$ in wave gauge. To emphasize the dependency of $f$ on $g'$, we also use the notation $f_{g'}$ for $f$. By the observation data in wave gauge, we mean  
\beqq\label{dataW}
\mcd_{W}(\delta; V) = \{(f_{g, *}g, f_{g, *}\psi, f_{g, *}\mcf): (g, \psi, \mcf) \in \mcd(\delta; V)\}
\eeqq
We shall prove
\begin{theorem}\label{main1}
Theorem \ref{main} holds if the data set $\mcd_F(\delta; \bullet)$ are replaced by $\mcd_W(\delta; \bullet)$. 
\end{theorem}

\bpf 
We divide the proof into three steps. 

\noindent \textbf{Step 1:}
We prove that the source-to-solution set $\mcd_W(\delta)$ determines the conformal class of the metric, namely there exists a diffeomorphism $\Psi: I(p^{(1)}_-, p^{(1)}_+)\rightarrow I(p^{(2)}_-, p^{(2)}_+)$ and $\gamma \in C^\infty(M)$ such that
\beq
\Psi^*\hat g^{(2)} = e^{2\gamma} \hat g^{(1)} \text{ in } I(p^{(1)}_-, p^{(1)}_+).
\eeq 

We follow the arguments  in Section 5 of \cite{KLU1}. It suffices to work with one copy of the Lorentzian manifold $(M, g)$ in wave gauge. We continue using the setup and notations in Section \ref{sec-wgauge}. The analysis of singularities in $\mcu^{(4)}$ we did in Prop.\ \ref{inter5} and Prop.\ \ref{symb}, immediately gave us the analogue of Prop.\ 3.3 and Prop.\ 3.4 of \cite{KLU1} (or similarly Theorem 3.3 and Prop.\ 3.4 in \cite{KLU}). In particular, for $s_0>0$ sufficiently small, we have: 
\begin{enumerate}
\item If $\bigcap_{j = 1}^4\gamma_{x^{(j)}, \theta^{(j)}}(\mbr_+) = \emptyset$ or $\bigcap_{j = 1}^4\gamma_{x^{(j)}, \theta^{(j)}}(\mbr_+) = q_0$ and the corresponding tangent vectors at $q_0$ are not linearly independent, then $\mcu^{(4)}$ is smooth in $M(\bt_0)\backslash \big(\mck\cup J_g(q_0)\big)$; 
\item If $\bigcap_{j = 1}^4\gamma_{x^{(j)}, \theta^{(j)}}(\mbr_+) = q_0$  and the corresponding tangent vectors at $q_0$ are linearly independent, then there exists $(x_j', \theta_j')$ sufficiently close to $(x^{(j)}, \theta^{(j)}), j = 1, 2, 3, 4$ such that $\bigcap_{j = 1}^4\gamma_{x'_j, \theta'_j}(\mbr_+) = q_0$. Moreover, by Prop.\ \ref{symb}, we can find $(x^{(j)}, \theta^{(j)})$ and $f^{(j)} \in  I^{\mu+1}(N^*Y_j; \sym^2)$ such that in $M(\bt_0)\backslash \mck$,  we have $\mcu^{(4)}  \in  I^{4\mu + \frac{3}{2}}(\La^{\hat g}_{q_0}; \sym^2)$ as well as $\sigma(\mcu^{(4)}) (y, \eta) \neq 0$ for any given $(y, \eta)\in \La^{\hat g}_{q_0}$.
\end{enumerate}

From the above two points, we can determine the set $\La^{\hat g}_{q_0}\backslash \mck$ using the data in the source-to-solution set $\mcd_W(\delta; V)$. As $s_0 \rightarrow 0$, the Hausdorff measure of $\mck$ in $\La^{\hat g}_{q_0}$ goes to zero so we determined the set $\La^{\hat g}_{q_0}$ up to a measure zero set. For  $V\subset M$ considered above, the light observation set of $q\in M$ is defined as $\mcp_V(q) = \mcl_q^+\cap V$, which is the projection of $\La^{\hat g}_{q_0}$ to $V$. We conclude that $\mcp_V(I(p_-, p_+))$ is determined by  $\mcd_W(\delta)$. The problem is reduced to the inverse problem with passive measurements studied in \cite{KLU2} and Theorem 2.5 of \cite{KLU1} tells that the differential structure of $I(p_-, p_+)$ and the conformal class of the metric can be determined.
\begin{remark}
For the vacuum background i.e.\ $\hat T^{(i)}(\hat g^{(i)}, \hat \psi^{(i)}) = 0, i = 1, 2$, one can complete the proof of the theorem  by directly applying Theorem 1.2 and Corollary 1.3 of \cite{KLU2}. This argument was used in \cite{LUW1} to determine  vacuum space-times. 
\end{remark}
\begin{remark}
The proof in Step 1 is simpler than those in \cite{KLU1} because of the no conjugate point assumption. When there are conjugate points, one can modify the proof by following the arguments in \cite{KLU1} line by line. In particular, one should replace the light observation set by the so called earliest light observation set. 
\end{remark}

\noindent \textbf{Step 2:} Without loss of generality, from now on we assume that  
$\hat g^{(1)} = e^{2\gamma}\hat g^{(2)}$ for some $\gamma\in C^\infty(M)$. 
From the data set, we know that the metrics $\hat g^{(1)} = \hat g^{(2)}$ on $V$ so that  $e^\gamma = 1$ on $V$. 


To further determine the conformal factor, we follow the argument in \cite{LUW} to make use of the symbols of the interaction terms. For this purpose, we need  the relation of principal symbols of $\mcu^{(4)}$ for two conformal metrics.  Let $\bfq_{\hat g^{(1)}}, \bfq_{\hat g^{(2)}}$ be the causal inverse of $\bold{P}_{\hat g^{(1)}}, \bold{P}_{\hat g^{(2)}}$ respectively. The Lagrangians $\La_{\hat g^{(1)}} = \La_{\hat g^{(2)}}$. We claim that the principal symbols of $\bfq_{\hat g^{(1)}}, \bfq_{\hat g^{(2)}} \in I^{-2}(N^*\diag\backslash \La_{\hat g^{(1)}})$ satisfy
\beq
\begin{gathered}
\sigma(\bfq_{\hat g^{(1)}}) = e^{2\gamma} \sigma(\bfq_{\hat g^{(2)}}),
\end{gathered}
\eeq
and  their principal symbols in $I^{-\frac{3}{2}}(\La_{\hat g^{(1)}} \backslash N^*\diag ; \bold{B})$ satisfy
\beq
\sigma(\bfq_{\hat g^{(1)}})(x, \xi, y, \eta) = e^{-\gamma(x)}\sigma(\bfq_{\hat g^{(2)}})(x, \xi, y, \eta) e^{3\gamma(y)},
\eeq
for $(x, \xi), (y, \eta)$ on the same bi-characteristics on $\La_{\hat g^{(1)}}$.

This result was proved in \cite[Prop.\ 4.6]{LUW} for the causal inverse of $\square_{\hat g}$ with no first order terms.  The same argument applies here but the first order term $V(x, \p)$ in \eqref{lineq} contribute to the subprincipal symbol of $\sigma(\bold{P}_{\hat g})$ so we only need to check that it transforms properly under conformal transformations.   

Recall see e.g.\ \cite[Appendix A.3.2]{Ho} that for any two Lorentzian metrics $g^{(1)}, g^{(2)}$ with $g^{(1)} = e^{2\gamma}g^{(2)}$, their corresponding Ricci curvatures satisfy 
\beq
\ric_{\mu\nu}(g^{(1)}) = \ric_{\mu\nu}(g^{(2)}) - 2[\nabla^{(2)}_\mu \p_\nu \gamma -\p_\mu \gamma \p_\nu \gamma] + (\lap \gamma - 2\|\nabla^{(2)} \gamma\|^2)g^{(2)}_{\mu\nu},
\eeq
where $\|\cdot\|$ stands for the standard Euclidean norm. 
So the difference has no derivative in the metric. From \eqref{eqredric}, we know that in corresponding wave gauges, the reduced Ricci curvature also satisfies the above equation i.e.\ 
\beq
\ric_{\hat g^{(1)}, \mu\nu}(g^{(1)}) = \ric_{\hat g^{(2)}, \mu\nu}(g^{(2)}) - 2[\nabla^{(2)}_\mu \p_\nu \gamma -\p_\mu \gamma \p_\nu \gamma] + (\lap \gamma - 2\|\nabla^{(2)} \gamma\|^2)g^{(2)}_{\mu\nu}.
\eeq
If we let $u^{(i)} = g^{(i)}-\hat g^{(i)}$ then $u^{(1)} = e^{2\gamma}u^{(2)}$. Let $V_{\mu\nu}^{(i)}(x, \p)(u^{(i)})$ be the first order terms in $\ric_{\hat g^{(i)}}(g^{(i)}), i = 1, 2$,  we obtain that
\beq
 V_{\mu\nu}^{(2)}(x, \p)(u^{(2)}) = V_{\mu\nu}^{(1)}(x, \p)(u^{(1)}) = V_{\mu\nu}^{(1)}(x, \p)(e^{2\gamma}u^{(2)}).
\eeq
Now the proof of Prop.\ 4.6 of \cite{LUW} especially equation (4.3) can be continued to obtain the conclusion. This proves our claim.\\

\noindent \textbf{Step 3:} 
Now we consider how the principal symbols of $\mcu^{(4)}$ are related for two conformal metrics. We first consider $\mch$ defined in Prop.\ \ref{intert}. Let  ${}^{(1)}\mch$ and ${}^{(2)}\mch$ be the defined with respect to two Einstein equation \eqref{pernon} where $\hat g$ is replaced by $\hat g^{(1)}$ and $\hat g^{(2)}$. In particular, they are polynomial functions of the linearized metric perturbation. We are interested in how the coefficients transform and we claim that their principal symbols on $\La_{q_0} \backslash \Theta$ satisfy the relation
\beq
\sigma( {}^{(1)}\mch) = e^{-8\gamma}\sigma({}^{(2)}\mch).
\eeq
Note this is for the principal symbol of the full term $\mch$ not just the leading term for some $\rho$ large. 
Since the coefficients of the terms $\hat G_i, i = 2, 3, 4$ depends on $\hat g$, we need to find them to see how they transform between two conformal metrics. We recall from the proof of Prop.\ \ref{intert} that $\hat G_i = P_i + \hat H_i, i = 1, 2, 3$ in which $u = g - \hat g$ is the metric perturbation,  $P_i$ are given by \eqref{eqpi} and $\hat H_i$ come from the semilinear term $Q_{\mu\nu}$ in \eqref{ricre} with two derivatives. 
What is important below is  the coefficients of $u$ in terms of $\hat g$, not the exact form of $\hat G_i$. We use $\mcg_{\la\alpha\beta}(u)$ introduced in  \eqref{eqmcg}. 
We start from $\hat H_2$ and see that it is a summation of terms like 
\beq
\hat g^\bullet \hat g_\bullet \hat g^\bullet \mcg_\bullet(u)\hat g^\bullet \mcg_\bullet(u) \text{ and } \hat g_\bullet \hat g^\bullet \mcg_\bullet (u) \hat g^\bullet \hat g^\bullet \p_\bullet u_{\bullet},
\eeq
where $\bullet$ denotes a generic index of $\hat g$.  Therefore, for the two conformal metrics, we have $\hat H^{(1)}_2 = e^{-4\gamma}\hat H^{(2)}_2$ for the same $u$. Next, $\hat H_3$ is a summation of terms like
 \beq
 \hat g^\bullet \hat g^\bullet u_\bullet \hat g_\bullet \hat g^\bullet \mcg_\bullet(u)\hat g^\bullet \mcg_\bullet(u),  \ \
 \hat g^\bullet u_\bullet \hat g^\bullet \mcg_\bullet(u)\hat g^\bullet \mcg_\bullet(u),
 \eeq
 and
 \beq
 u_\bullet \hat g^\bullet \mcg_\bullet (u) \hat g^\bullet \hat g^\bullet \p_\bullet u_{\bullet},  \ \
 \hat g_\bullet \hat g^\bullet\hat g^\bullet u_\bullet \mcg_\bullet (u) \hat g^\bullet \hat g^\bullet \p_\bullet u_{\bullet}.
 \eeq
 Therefore,  $\hat H^{(1)}_3 = e^{-6\gamma}\hat H^{(2)}_3$. Similarly, one can find that  $\hat H_4$ are summations of 
  \beq
 \hat g^\bullet \hat g^\bullet \hat g^\bullet u_\bullet  u_\bullet \hat g_\bullet \hat g^\bullet \mcg_\bullet(u)\hat g^\bullet \mcg_\bullet(u),  \ \
\hat g^\bullet \hat g^\bullet u_\bullet u_\bullet \hat g^\bullet \mcg_\bullet(u)\hat g^\bullet \mcg_\bullet(u),
 \eeq
 and
 \beq
 u_\bullet \hat g^\bullet \hat g^\bullet  u_\bullet \mcg_\bullet (u) \hat g^\bullet \hat g^\bullet \p_\bullet u_{\bullet},  \ \
 \hat g_\bullet \hat g^\bullet\hat g^\bullet \hat g^\bullet  u_\bullet u_\bullet \mcg_\bullet (u) \hat g^\bullet \hat g^\bullet \p_\bullet u_{\bullet}.
 \eeq
Thus $\hat H^{(1)}_4 = e^{-8\gamma}\hat H^{(2)}_4$. By similar argument, we find that 
\beq
P^{(1)}_2 = e^{-4\gamma} P^{(2)}_2, \ \ P^{(1)}_3 = e^{-6\gamma} P^{(2)}_3, \ \ P^{(1)}_4 = e^{-8\gamma} P^{(2)}_4.
\eeq 
Finally,  it is straightforward to see that the transformation of $\mch_i, i = 1, \cdots, 5$ are ${}^{(1)}\mch_i = e^{-8\gamma} \big({}^{(2)}\mch\big)$. Thus, we have proved the claim.\\


Now we compare the symbols. From Prop.\ \ref{intert}, we know that 
\beq
\sigma(\mcu^{(i)})(q, \eta) = \sigma(\bfq_{\hat g^{(i)}})(q, \eta, q_0, \zeta) \sigma(\mch^{(i)})(q_0, \zeta), \ \ i = 1, 2,
\eeq
where $(q, \eta)$ is jointed to $(q_0, \zeta)$ by bi-characteristics. Let $ = \bfq_{\hat g^{(i)}}({}^{(i)}f), i = 1, 2$, from Step 2, we know that 
\beq
\sigma({}^{(1)}v) = e^{-\gamma(q_0)} \sigma({}^{(2)}v).
\eeq
Here, we used the fact that $e^{\gamma(q)} = 1$ on $V$. Therefore, we get
\beq
\sigma({}^{(1)}\mch) = e^{-4\gamma(q_0)} e^{-8\gamma(q_0)} \sigma({}^{(2)}\mch).
\eeq
Finally, we have 
\beq
\sigma(\mcu^{(1)})(q, \eta) =  e^{-9\gamma(q_0)}\sigma(\mcu^{(2)})(q, \eta).
\eeq
Since the data set $\mcd_W(\delta)$ is the same, we know that the symbols should be the same. So we get $e^{\gamma(q_0)} = 1$. Since this is true for all $q_0\in I(p_-, p_+)$, the proof is finished.   
\epf

\section{The inverse problem in Fermi coordinates}\label{sec-fermi}
\begin{figure}[htbp]
\centering
\includegraphics[scale=0.7]{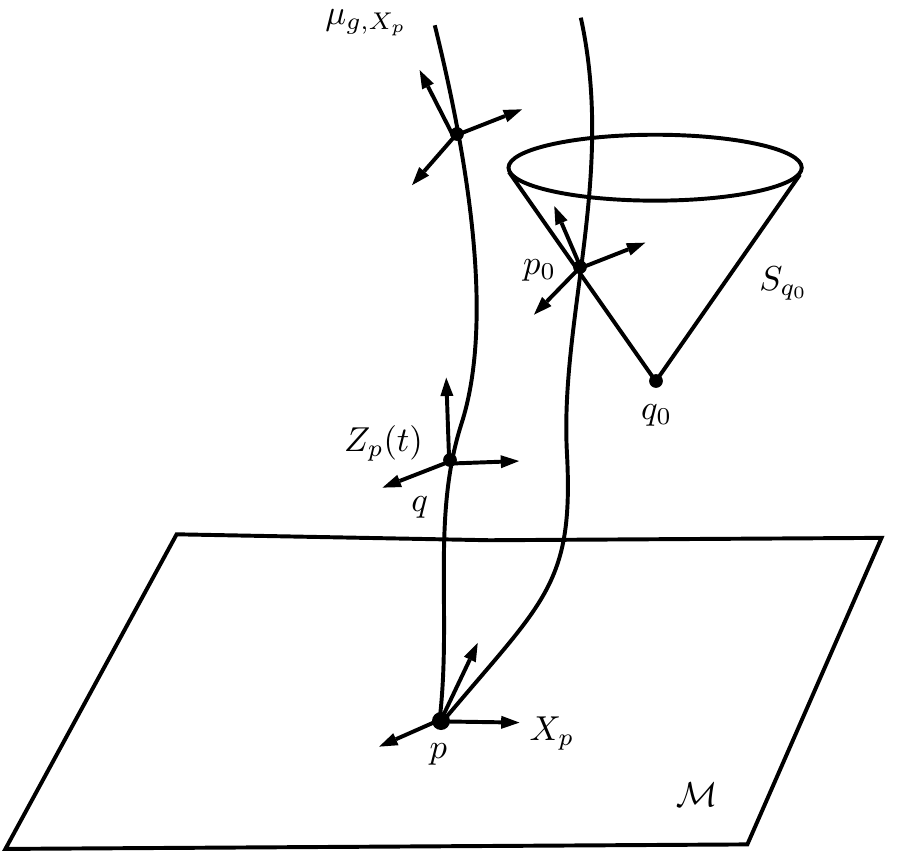}
\caption{Fermi coordinates along time-like geodesics $\mu_{g, X_p}.$ $q_0$ is a point source of gravitational waves. $S_{q_0}$ represents the future pointing light cone at $q_0$, where the gravitational wave is propagating. When the geodesic intersect $S_{q_0}$ transversally at $p_0$, the gravitational wave ( the singularity of the metric) is detected at $p_0.$}
\label{fermipic}
\end{figure}

We complete the proof of Theorem \ref{main} in this section. We begin with some  physical considerations, that is how  experiments are designed to detect gravitational perturbations. We follow the nice presentation in \cite[Chapter 37]{MTW} closely. Roughly speaking, one considers detecting gravitational waves in the {\em proper reference frame} attached to the center of mass of a mechanical detector (e.g.\ a vibrating bar). To construct the frame, let $\tau$ be the proper time measured by the observer's clock (attached to the detector). Let $\mu(\tau), \tau\in \mbr$ be the world line of the observer, which is a time-like geodesic. Then one chooses a orthonormal frame $\vec e_a, a = 0, 1, 2, 3$ at $\tau = 0$ such that $\vec e_0 = \dot \mu(0)$. Finally the frame changes from point to point along the world line following the parallel transport, see Figure \ref{fermipic}. The details of this construction can be found in \cite[Section 13.6]{MTW}. Here, we are interested in  non-accelerated observers hence the proper reference frame is just the Fermi coordinates associated with the time-like geodesic. This is the observation coordinate we study in this section. We remark that in this set up, one can think that the gravitational disturbance is detected at the location of the detector, in other words, on the time-like geodesic. 

Now we define Fermi coordinates for our problem more rigorously. For any $\hat p\in \{0\}\times \mcm$, we choose a basis $X^i_{\hat p} \in T_{\hat p} M, i = 0, 1,2, 3$. We assume that $X^0_{\hat p}$ is time-like. Let $\mu_{g}(t), t\in \mbr$ be the time-like geodesic such that $\mu_{g}(0) = \hat p, \dot \mu_{g}(0) = X^0_{\hat p}$. To emphasize the dependency of $\mu_g$ on the frame $X_{\hat p}$, we denote it by $\mu_{g, X_{\hat p}}.$ Then we parallel translate the basis $X_{\hat p}^i, i = 0, 1, 2, 3$ to  get vector fields $Z^i(t)$ (a frame) along $\mu_{g, X_{\hat p}}(t), t\in \mbr$. Let $q = \mu_{g, X_{\hat p}}(z_0)$ be any point on the geodesic. For any $q_0\in M$, the Fermi coordinate along $\mu_{g, X_{\hat p}}$ are given by $( (z_i)_{i = 0}^3\in \mbr^4)$ such that
\beq
q_0 = \Phi_{g, X_{\hat p}}(z_0, z_1, z_2, z_3) = \exp_{g, q} \bigg(\sum_{i = 1}^3 z_i Z^i_q\bigg).
\eeq
See Figure \ref{fermipic}. 
Here, we benefited from the no conjugate point assumption so that the exponential map is defined for $\mbr^4$. We emphasize that to specify the Fermi coordinates, we need to take $p\in \{0\}\times \mcm$ and choose frames $X^i_p, i = 1, 2, 3,4$. However, we remark that if the metric $g$ is known (which is the case on $V$ for the inverse problem we have been  considering), the Fermi coordinates associated with different geodesics can be obtained from each other. 

Now we let $\hat \zeta = X^0_{\hat p}$. For $\delta > 0$ and small, we let  
\beq
\bar B_{\delta}(\hat p, \hat \zeta) = \{(p, \zeta) \in T^*M: \text{dist} \big((p, \zeta) - (\hat p, \hat \zeta)\big) < \delta,\ \ p\in \{0\} \times \mcm\},
\eeq
where the distance is defined using the Sasaki metric on $T^*M$. We shall take and fix vector fields $X^i \in T\mcm, i = 1, 2, 3$ such that $X^i(\hat p) = X^i_{\hat p}$. Then for each $p \in \bar B_{\delta}(\hat p, \hat \zeta)$, we get a frame $(\zeta, X^i(p))$ at $p$ which can be used to define the Fermi coordinates along the time-like geodesic from $(p, \zeta)$. For simplicity, we use $\zeta$ to indicate the frame as $X^i$ are fixed on $\mcm.$ For example, these geodesics will be denoted by $\mu_{g, \zeta}$.  
 The Fermi coordinates for $(p, \zeta)\in \bar B_{\delta}(\hat p, \hat \zeta)$ give a coordinate system for a neighborhood of $\hat \mu$.  We can formulate the observation set in these Fermi coordinate as
\beqq\label{dataF}
\begin{split}
 \mcd_F(\delta; V) = \{ ( \Phi_{g, \zeta}^*g, \Phi_{g, \zeta}^*\psi, \Phi_{g, \zeta}^*\mcf): (g, \psi, \mcf) \in \mcd(\delta; V), \ \ (p, \zeta) \in \bar B_{\delta}(\hat p, \hat \zeta) \}.
\end{split}
\eeqq
Interpreted more physically, what we just described is that many detectors are placed near the observer at time $t = 0$ and they are set to move along different directions. The data is the gravitational waves detected by these detectors, see Figure \ref{fermipic}. \\

For the inverse problem we consider, the main issue now is whether we can still observe in these coordinates the waves (singularities) we produced in wave gauge. Let's make the situation more clear. Let $\eps = (\eps_1, \eps_2, \eps_3, \eps_4)$ and  $g_\eps$ be a family of metrics on $M$. Here, $g_0 = \hat g$ is the background metric. Let $f_{\eps}$ be the $(g_\eps, g_0)$ wave map and $\tilde g_\eps = f_{\eps, *} g_\eps$ be the image of $g_\eps$ in wave gauge. We assume that $f_\eps = \Id$ on $\{0\}\times \mcm.$ Notice that when $\eps = 0$, $f_0 = \Id$. Now we consider two Fermi coordinates. Let $\Phi_{g_\eps, X_p}$ be the Fermi coordinates in $(M, g_\eps)$ at $p$ with a frame $X_p$ in $T_p^*M$. Let $Z^i, i = 0, 1, 2, 3$ be the frame along the geodesic from $p$.  For Fermi coordinates in wave gauge, we take $\tilde Z^i_\eps = f_{\eps, *}Z^i, i = 0, 1, 2, 3$ to be the frame associated with $\mu_{\tilde g_\eps, \tilde X_p}(t) = f_\eps(\mu_{g_\eps, X_p}(t)), t\in \mbr$, where $\tilde X_p = f_{\eps, *}(X_p)$. If we define the Fermi coordinates in $(M, \tilde g_\eps)$ using these frames, we have the following commutative diagram 
\begin{center}
\begin{tikzpicture}
  \matrix (m) [matrix of math nodes,row sep=3em,column sep=4em,minimum width=2em]
  {
     (M, g_\eps) &  (M, \tilde g_\eps) \\
     (\mbr^4, Z^i) & (\mbr^4, \tilde Z^i_\eps) \\};
  \path[-stealth]
    (m-2-1) edge node [left] {$\Phi_{g_\eps, X_p}$} (m-1-1)
    (m-1-1) edge node [above] {$f_\eps$} (m-1-2)
    (m-2-1) edge node [below] {$\Id$} (m-2-2)
    (m-2-2) edge node [right] {$\Phi_{\tilde g_\eps, \tilde X_p}$} (m-1-2);
\end{tikzpicture}
\end{center}
Here, we emphasized which frame of $\mbr^4$ we are referring to because only under the above two frames, the map $f_{\eps, *}$ becomes the identity map $\Id$. If we pull back the metrics to Fermi coordinates, we have 
\beq
\Phi_{g_\eps, X_p}^*g_\eps = \Phi_{\tilde g_\eps, \tilde X_p}^*\tilde g_\eps.
\eeq
In particular, the left hand side is what we can observe in the data set $\mcd_F(\delta)$. So it suffices to analyze the singularities of the right hand side. Let's first clarify what singularities we should consider. In Section \ref{sec-wgauge}, we studied the singularities in the interaction term 
\beq
\mcu^{(4)} = \p_\eps \tilde g_\eps|_{\eps = 0} = \p_{\eps_1} \p_{\eps_2}\p_{\eps_3}\p_{\eps_4} \tilde g_\eps|_{\eps = 0} 
\eeq
in wave gauge. By choosing four distorted plane waves, this is a Lagrangian distribution in $I^{4\mu + \frac 32}(\La^{\hat g}_{q_0}\backslash \Theta; \sym^2)$, where $q_0$ is the point where the waves interact. For each point $q = \mu_{\tilde g_\eps, \tilde X_p}(t)$ on the geodesic, we will consider singularities in the Fermi coordinates in wave gauge
 \beqq\label{eqmcu1}
 \begin{split}
 \p_\eps \Phi_{\tilde g_\eps, \tilde X_p}^*\tilde g_\eps|_{\eps = 0}= \p_\eps \tilde g_\eps(\Phi_{\tilde g_\eps, \tilde X_p})(D\Phi_{\tilde g_\eps, \tilde X_p}, D\Phi_{\tilde g_\eps, \tilde X_p})|_{\eps = 0}.
 \end{split}
 \eeqq
 Here, $D(\bullet)$ stands for the Jacobian of $\bullet.$ Notice that at $\eps = 0$, we have $f_0= \Id$ and $\tilde g_0=  g_0$. Taking the derivatives in $\eps$, we have
 \beq
 \p_\eps \Phi_{\tilde g_\eps, \tilde X_p}^*\tilde g_\eps|_{\eps = 0}    = I_1 + I_2 + I_3 + I_4, 
 \eeq
 where 
 \beqq\label{eqsp4}
 \begin{split}
I_1 &= \p_\eps \tilde g_\eps|_{\eps=0}(\Phi_{g_0, X_p})(D\Phi_{g_0, X_p}, D\Phi_{g_0, X_p}) = \Phi_{g_0, X_p}^*(\p_\eps \tilde g_\eps|_{\eps = 0}), \\
 I_2 &=  \sum_{j = 0}^3 \p_jg_0(\Phi_{g_0, X_p})(D\Phi_{g_0, X_p}, D\Phi_{g_0, X_p}) \p_\eps \Phi^j_{\tilde g_\eps, \tilde X_p}|_{\eps = 0},\\
 I_3 &=  2 g_0(\Phi_{g_0, X_p})(\p_\eps D\Phi_{\tilde g_\eps, \tilde X_p}|_{\eps = 0}, D\Phi_{g_0, X_p}).
 \end{split}
 \eeqq
 These three terms are obtained when four derivatives in $\eps_i$ fall on one mapping in \eqref{eqmcu1}. The rest of the derivatives are collected in $I_4$, which is a summation of terms like 
 \beq
 (\mcv \tilde g_\eps)(\Phi_{\tilde g_\eps, \tilde X_p})(\mcv (D\Phi_{\tilde g_\eps, \tilde X_p}), D\Phi_{\tilde g_\eps, \tilde X_p}) (\mcv \Phi^j_{\tilde g_\eps, \tilde X_p})|_{\eps = 0},
 \eeq
 where $\mcv = \p^\alpha_\eps$ with $\alpha\in \mathbb{Z}^4, |\alpha|\leq 3.$ The $I_1$ term is just the pull back of the singularities in wave gauge by a smooth map. This is the singularity which we can construct in wave gauge and we hope to observe in the data set. Eventually, we will see that $I_4$ is smooth thus can be ignored.  The obstruction is that $I_2$ and $I_3$ have the derivative of $\Phi_{\tilde g_\eps, \tilde X_p}$ which may contain other singularities. These are the singularities caused by the coordinate change. We need to distinguish them from the singularities in $I_1$. The key point below is that we shall only consider  singularities along the geodesic, which also agrees with the physical explanation that the gravitational waves are supposed to be detected at the detector.

Let's start with the term $I_3$. By definition, for any point $q$ on the geodesic, the map $\Phi_{\tilde g_\eps, \tilde X_p}$ is just the normal coordinates based on $q$. Therefore, the Jacobian $D\Phi_{\tilde g_\eps, \tilde X_p} = \Id$ at $T_0(T_qM)$, see e.g.\ \cite[Lemma 5.10]{Lee}. So, if we restrict $I_3$ to the associated geodesics, we get $I_3 = 0$ for all $\eps$.  When $\eps = 0$, the geodesic is $\mu_{g_0, X_p}$. So it suffices to analyze the singularities of $I_1, I_2$ only along the associated geodesics.

Now we consider the restriction of conormal distributions to a time-like geodesic. 
\begin{lemma}\label{lemres}
Let $S$ be a co-dimension one submanifold of $M$ and $\gamma(t), t\in \mbr$ be a smooth curve on $M$. Suppose that $\gamma$ intersects $S$ transversally only at a point $p_0$. Then for any $f \in I^m (N^*S), m\in \mbr$, the restriction of $f$ to $\gamma$, denoted by $\tilde f = f|_{\gamma}$, is in $I^{m+\frac{3}{4}}(\gamma; N^*\{p_0\})$. This means that $\tilde f$ is a distribution on $\gamma$ and conormal to $p_0$. Furthermore, if the principal symbol $\sigma(f) \neq 0$ at $p_0$, then the principal symbol $\sigma(\tilde f)$ is also non-vanishing at $p_0$.
\end{lemma}

\bpf
Since $\mu$ intersect $S$ at $p_0$ transversally, we can choose local coordinates $x = (x^0, x^1, x^2, x^3)$ near $p_0$ such that locally $S = \{x^0 = 0\}$ and $\gamma = \{x^1 = x^2 = x^3 = 0\}$ so $p_0 = 0$. Since $f$ is a conormal distribution, we can write it as an oscillatory integral 
\beq
f(x^0, x^1, x^2, x^3) = \int_\mbr e^{i x^0 \xi_0}a(x^0, x^1, x^2, x^3, \xi_0)d \xi_0,
\eeq
where $a \in S^{m+\frac{n}{4} - \frac{k}{2}}(\mbr^4\times \mbr)$ with $n = 4, k = 1$. The restriction of $f$ to $\mu$ is 
\beq
\tilde f(x_0) = f(x^0, 0, 0, 0) = \int_\mbr e^{i x^0 \xi_0}a(x^0, 0, 0, 0, \xi_0)d \xi_0,
\eeq
so that $a \in S^{\tilde m + \frac{1}{4} - \frac{1}{2}}(\mbr\times \mbr)$. Thus we find that $\tilde m = m + \frac 34.$ So $\tilde f\in I^{m + \frac 34}(\gamma; N^*\{p_0\})$. The conclusion about the principal symbols are easy to see from these formulas. 
\epf

Let $S_{q_0} = \pi(\La^{\hat g}_{q_0})$ be the projection of $\La^{\hat g}_{q_0}$ to $M$. So $S_{q_0}$ is the future pointing light-cone at $q_0$ and  a co-dimension one submanifold. Here, for simplicity, we ignored the set $\Theta$ and its projection $\mck$ as it does not matter for the analysis below. Suppose the geodesic $\mu_{\hat g, X_p}$ intersect $S_{q_0}$ transversally at $p_0$. Applying the above result, we see that $I_1(t) \in  I^{\mu + \frac{3}{4}}(\mu_{\hat g, X_p}; N^*\{p_0\})$ and we know that the singularity is non-vanishing. See Figure \ref{fermipic}.

Next we  use this lemma to analyze singularities of $I_2$ along $\mu_{g_0, X_p}(t)$. 
\begin{lemma}
Assume that $g_\eps$ is such that $\p^\alpha_\eps \tilde g_\eps |_{\eps = 0}$ is smooth for $|\alpha|< 4$ and $\p_\eps\tilde g_\eps|_{\eps = 0}$ is a conormal distribution to the light cone at $S_{q_0}$ at $q_0$. Assume that $\mu_{\hat g, X_p}(t)$ intersects $S_{q_0}$ transversally at $p_0$. Then the restriction of $I_2(t)$ to the geodesic $\mu_{\hat g, X_p}(t), t>0$ is  in $I^{\mu -1 + \frac{3}{4}}(\mu_{\hat g, X_p}; N^*\{p_0\}). $
\end{lemma}
\bpf
Observe that $\Phi_{\tilde g_\eps, \tilde X_p}$ along $\mu_{\tilde g_\eps, \tilde X_p}$ is just the geodesic itself. More precisely, let $\gamma_\eps(t) = \mu_{\tilde g_\eps, \tilde X_p}(t), t>0$. Then $\gamma_\eps$ satisfies the geodesic equation
\beq
\begin{gathered}
\ddot\gamma^k_\eps(t) + \Gamma^k_{\eps, ij}(\gamma_\eps(t)) \dot \gamma_\eps^i(t) \dot \gamma_\eps^j(t) = 0\\
\gamma_\eps(0) = p, \ \ \dot \gamma_\eps(0) = X^0_p.
\end{gathered}
\eeq
Here, we write down the equation in some local coordinates. This is a second order nonlinear ODE system. Since the Christoffel symbol $\Gamma_{\eps, ij}^k$ of $\tilde g_\eps$ may be singular, the solution may develop singularities. Since we only need the singularities of $\p_\eps \gamma_\eps(t)|_{\eps = 0}$, we resort to asymptotic analysis. 

We consider the asymptotic expansion of the Christoffel symbol as $\eps\rightarrow 0.$ First of all, the order $\eps^0$ term is given by $G_{ij}^k(t) = \Gamma_{0, ij}^k(\gamma_0(t))$. Next, the $\eps_a, a = 1, 2, 3, 4$ term is of the form 
\beq
G^k_{a; ij}(t) = G^k_{a; ij}(t, \p^\alpha_\eps \p^k \tilde g_\eps|_{\eps = 0}), \ \ |\alpha| \leq 1, k \leq 2, 
\eeq
where the function is smooth in its arguments. By our assumption on $\tilde g_\eps$, $G_1(t)$ is smooth on $M$. Now consider the $\eps_a\eps_b, 1\leq a< b\leq 4$ terms which can be written as
\beq
G^k_{ab; ij}(t) = G^k_{ab; ij}(t, \p^\alpha_\eps \p^k \tilde g_\eps|_{\eps = 0}), \ \ |\alpha| \leq 2, k\leq 3,
\eeq
and the functions are smooth in its argument so they are also smooth on $M$. The next term i.e.\ $\eps_a\eps_b\eps_c, 1\leq a< b< c \leq 4$ can be written similarly as
\beq
G^k_{abc; ij}(t) = G^k_{abc; ij}(t, \p^\alpha_\eps \p^k \tilde g_\eps|_{\eps = 0}), \ \ |\alpha| \leq 3, k\leq 4,
\eeq
and they are smooth on $M$. The order $\eps_1\eps_2\eps_3\eps_4$ term is given by 
\beq
G_{1234; ij}^k(t) = G^k_{1234; ij}(t, \p^\alpha_\eps \p^k \tilde g_\eps|_{\eps = 0}), \ \ |\alpha| \leq 4, k\leq 5.
\eeq
By the assumption $\p_\eps \tilde g_\eps|_{\eps = 0}\in I^\mu(N^*S_{q_0})$, we know that this term is in $I^{\mu+1}(N^*S_{q_0})$. We remark that it is possible that the singularities may cancel each other so the term is smoother. But the worst case in terms of strength of singularities is that the term belongs to $I^{\mu+1}(N^*S_{q_0})$.

Now consider the asymptotic expansion of $\gamma_\eps(t)$. We write
\beq
\gamma_\eps(t) = \gamma_0(t) + \sum_{a = 1}^4 \eps_a A_a(t) + \sum_{1 \leq a< b \leq 4}\eps_a\eps_b A_{ab}(t) + \sum_{1\leq a< b< c \leq 4} \eps_a\eps_b\eps_c A_{abc}(t) + \eps_1\eps_2\eps_3\eps_4 A(t) 
\eeq
modulus terms in $\bigcup_{i = 1}^4 O(\eps_i^2).$ 
The equations for $\gamma_0$  is 
\beq
\begin{gathered}
\ddot\gamma^k_0(t) + \Gamma^k_{0, ij}(\gamma_0(t)) \dot \gamma_0^i(t) \dot \gamma_0^j(t) = 0\\
\gamma_0(0) = p, \ \ \dot \gamma_0(0) = X^0_p,
\end{gathered}
\eeq
which is just the geodesic equation for $g_0$. The order $\eps_a$ terms in the equation gives the equation for $A_a(t), a = 1, 2, 3, 4,$
\beq
\begin{gathered}
\ddot A^k_a(t) + G_{a; ij}^k(t) \dot \gamma_0^i(t) \dot \gamma_0^j(t) +2\Gamma^k_{0, ij}(\gamma_0(t)) \dot A_a^i(t) \dot \gamma_0^j(t) = 0\\
A_a(0) = 0, \ \ \dot A_a(0) = 0.
\end{gathered}
\eeq
This is a second order linear ODE for $A_a(t)$ and since the coefficients are all smooth, we conclude that $A_a(t)$ are smooth in $t$. This procedure can be continued to conclude that $A_{ab}(t), A_{abc}(t)$ are all smooth in $t$. Finally we reached the equation for $A(t)$ which comes from the order $\eps_1\eps_2\eps_3\eps_4$ term of the geodesic equation. The equation has many terms but schematically, the equation is of the form
\beq
\begin{gathered}
\ddot A^k(t) + G_{1234; ij}^k(t) \dot \gamma_0^i(t) \dot \gamma_0^j(t) +F = 0\\
A(0) = 0, \ \ \dot A(0) = 0,
\end{gathered}
\eeq
where $F$ denotes a smooth function which depends on all the previous asymptotic terms. From Lemma \ref{lemres}, we know that the coefficients $G_{1234; ij}^k(t) \in I^{\mu + 1 + \frac{3}{4}}(\mu_{g_0, X_p}; N^*\{p_0\})$ along the geodesic. Solving the above equation (integrating in $t$ twice), we get $A^{k}(t) \in I^{\mu -1 + \frac{3}{4}}(\mu_{g_0, X_p}; N^*\{p_0\})$. Therefore, $I_2(t)$ term along the geodesic $\gamma_0$ is at most $I^{\mu -1 + \frac{3}{4}}(\mu_{g_0, X_p}; N^*\{p_0\}). $
\epf

From the proof above, we see that the term $I_4$ in \eqref{eqsp4} is smooth along the geodesic. We conclude that the most singular term among $I_a, a = 1, 2, 3$ is $I_1$ and we know that the singularity is non-vanishing, under the assumption that the geodesic $\mu_{\hat g, X_p}$ intersects the singular support of $\mcu^{(4)}$ transversally. However, this is always possible by our construction of Fermi coordinates.  So we proved that one can observe the singularity in Fermi coordinate.  This means that $\mcd_F(\delta; V)$ determines the light observation set. We follow the argument in Section \ref{sec-wave} to prove the determination of the conformal class of the metric, and by considering the symbol on the geodesics, we can determine the conformal factor as in Section \ref{sec-wave}. This complete the proof of Theorem \ref{main}.

\section{The Einstein-scalar field equations}\label{sec-einscal}
We consider the Einstein-scalar field equations studied in \cite{KLU1}. We show that this model fits in our scheme. Moreover, we improve the result of \cite{KLU1} to show that one can determine the metric in $I_{\hat g}(p_-, p_+)$, not just the conformal class. 

Consider  the stress-energy tensor given by scalar fields $\phi_l, l = 1, 2, \cdots, L$, 
\beqq\label{einmatsource}
\begin{gathered}
T^{scalar}_{jk}(g, \phi) = \sum_{l = 1}^L (\p_j \phi_l \p_k \phi_l - \ha g_{jk} g^{pq}\p_q\phi_l \p_q \phi_l) - V(\phi) g_{jk},
\end{gathered} 
\eeqq
where $V$ is a smooth function.  Assume further that the fields $\phi_l$ satisfy the wave equations
\beq
\square_g\phi_l - V_l'(\phi) = \mcf_l^\Psi,  \ \ l = 1,2, \cdots, L,
\eeq
where $V_l'(s) = \p_{s_l} V(s)$ and $\mcf_l^\Psi$ are sources which generates the scalar fields on $(M, g)$.  We take 
\beqq\label{scasour}
T_{sour} = T^{scalar}(g, \phi) + \mcf 
\eeqq
as the stress-energy tensor  in \eqref{einsour0} where $\mcf$ is another source perturbation term. Let $\overline\mcf = (\mcf, \mcf^\Psi)$. The system \eqref{einsour00} becomes the coupled system  
\beqq\label{einmat}
\begin{gathered}
\ein(g) = T^{scalar}(g, \phi) + \mcf \text{ in } M(\bt_0)\\
\square_g\phi_l - V_l'(\phi) = \mcf_l^\Psi,  \ \ l = 1,2, \cdots, L,\\
g = \hat g, \ \ \phi_l  = \hat \phi_l, \text{ in } M(\bt_0)\backslash J_g^+(\supp (\overline\mcf) ).
\end{gathered}
\eeqq
The background fields $(\hat g, \hat \phi)$ are solutions to \eqref{einmat} with $\overline\mcf = 0$.  Again, once we determine the set of $\overline\mcf$ so that the system \eqref{einmat} is well-posed, we can take $T_{sour}$ in \eqref{scasour} as the source which is defined using the fields $\phi_l$ and $\mcf$. Then we can just work with the Einstein equations, ignoring the equations for the fields $\phi_l$. 

For this concrete model, the data set and microlocal linearization stability condition can be described more explicitly. 
 The data set is the source-to-solution set as in equation (11) of \cite{KLU1}.  For $\delta > 0$, we let
\beqq\label{datasca}
\begin{gathered}
 \mcd^{sca}(\delta; V) \doteq \{(g, \phi, \overline\mcf):  \mcf \in C^{4}(M(\bt_0)), \ \ \|\mcf\|_{C^{4}(M(\bt_0))}< \delta, \ \ \supp (\mcf) \subset V\\
\text{ such that there is a unique solution $(g, \phi)\in C^{4}(M(\bt_0))$ to \eqref{einmat}}\}.
\end{gathered}
\eeqq
We denote the observation in Fermi coordinates by $\mcd^{sca}_F(\delta; V).$ 
The microlocal linearization stability conditions can be reformulated as in \cite{KLU1}. Since $T_{sour}$ is given by \eqref{scasour} and the fields $\phi_l$ satisfy the wave equations, one can derive the conservation law $\div_g T_{sour} = 0$ more explicitly as 
\beq
\ha g^{pk} \nabla_p \mcf_{jk} + \sum_{l = 1}^L \mcf^\Psi_l \p_j \phi_l = 0, \ \ j = 0, 1, 2, 3,
\eeq
see equation (159) of \cite{KLU1}. Let $\overline\mcf_\eps = \eps (f, f^\Psi)$ depending on a small parameter $\eps$. We have the linearized conservation law
\beq
\ha \hat g^{pk} \hat \nabla_p f_{jk} + \sum_{l = 1}^L f^\Psi_l \p_j \hat \phi_l = 0, \ \ j = 0, 1, 2, 3,
\eeq
see equation (160) of \cite{KLU1}. In this case, the linearized conservation law does not depend on the linearization of $g$ and $\phi$. It is easy to see condition (2) of Def.\ \ref{mlstab} is satisfied. Now we can state the microlocal linearization stability condition for the system \eqref{einmat}.

\begin{definition}\label{mlssca}
Let $Y$ be a two dimensional submanifold of a space-like surface in $M$ and $Y\subset V\subset M$. Let $(y, \eta)\in N^*Y$ with $\eta$ a light-like co-vector.  Let $A\in \mbr^{10}$ satisfy
\beq
\hat g^{\beta \alpha}(y)  \eta_\beta A_\alpha = 0
\eeq
and $B\in \mbr^L$.  We say that the Einstein-scalar field equation \eqref{einmat} satisfies the microlocal linearization stability condition if we can find a family $(g_\eps, \phi_\eps, \overline \mcf_\eps)\in \mcd^{sca}(\delta; V)$ such that 
 \begin{enumerate}
\item  $\overline \mcf_\eps = (\mcf_\eps, \mcf^\Psi_\eps)$ are compactly supported in a conic neighborhood $W$ of $(y, \eta)$ in $N^*V$
\item   $ f = \p_\eps \mcf_\eps|_{\eps = 0} , f^\Psi = \p_\eps \mcf^\Psi_\eps|_{\eps = 0}$ are in $I^{\mu+ 1}(N^*Y), \mu < -17$  such that 
\beq
\sigma(f)(y, \eta) = A, \ \ \sigma(f^\Psi)(y, \eta) = B.
\eeq
\end{enumerate}
\end{definition}  
Examples when this condition holds are illustrated  in \cite[Appendix C]{KLU1}. Now we state and prove the result for Einstein-scalar field equations, which improves the result in \cite{KLU1}.
\begin{theorem}\label{mainsca}
Let $(M, \hat g^{(i)}), i = 1, 2$ be two four-dimensional smooth, time oriented, globally hyperbolic Lorentzian manifolds such that  $\hat g^{(i)}$ satisfy assumption (A1).  
Let $\mu_{\hat g^{(i)}}(t), t\in [-1, 1]$ be time-like geodesics on $(M, \hat g^{(i)})$. Let $V$ be open relatively compact neighborhoods of $\mu_{\hat g^{(i)}}([-1, 1])$ and  $V\subset M^{(i)}(\bt_0)$ for some $\bt_0>0$.  Let $-1<s_-<s_+ < 1$, $p_\pm^{(i)} = \mu_{\hat g^{(i)}}(s_\pm)$ and $p^{(1)}_\pm = p^{(2)}_\pm$. 
For  $\delta > 0$, we consider solutions $(g^{(i)}, \phi^{(i)})$ to the Einstein equations
\beq 
\begin{gathered}
\ein(g^{(i)}) =   {T}^{scalar, (i)} +  \mcf^{(i)}, \text{ in } M^{(i)}(\bt_0),\\
\square_{g^{(i)}}\phi_l^{(i)}- V_l'(\phi^{(i)}) =  \mcf_l^{\Psi, (i)},  \ \ l = 1,2, \cdots, L,\\
g^{(i)} = \hat g^{(i)}, \ \ \phi_l^{(i)}  = \hat \phi_l^{(i)}, \text{ in } M^{(i)}(\bt_0)\backslash J_{g^{(i)}}^+(\supp (\overline\mcf^{(i)}) ),
\end{gathered}
\eeq 
where $(g^{(i)}, \phi^{(i)}, \overline \mcf^{(i)}) \in \mcd^{sca, (i)}(\delta; V)$ defined as \eqref{datasca}. Suppose that 
 \begin{enumerate}
 \item  the microlocal linearization stability condition Definition \ref{mlssca} holds for $\mcd^{sca, (i)}(\delta; V)$.
 \item the source-to-solution maps satisfy 
\[
 \mcd^{sca, (1)}_F(\delta; V) = \mcd^{sca, (2)}_F(\delta; V).
\]
\end{enumerate}
Then there is a diffeomorphism $\Psi: I_{\hat g^{(1)}}(p^{(1)}_-, p^{(1)}_+)\rightarrow I_{\hat g^{(2)}}(p^{(2)}_-, p^{(2)}_+)$ such that  $\Psi^*\hat g^{(2)} = \hat g^{(1)}$ in $I_{\hat g^{(1)}}(p^{(1)}_-, p^{(1)}_+)$. 
\end{theorem}

\bpf
We reduce the problem into the framework of Theorem \ref{main}. The only problem is that $T(g, \phi)$ depends on the derivatives of the scalar fields $\phi$.  

The local well-posedness of the system \eqref{einmat} in wave gauge was proved in \cite{KLU1}, see Section 3.1.1 and Appendix B there.
We  find that 
\beq
\rho_{jk}(g, T^{scalar}) = \sum_{l= 1}^L \p_j \phi_l \p_k \phi_l - 3V(\phi) g_{jk}, \ \ j, k = 0, 1, 2, 3.
\eeq
So the reduced Einstein-scalar field equation is 
\beq
(\ric_{\hat g}(g))_{jk} =  \sum_{l= 1}^L \p_j \phi_l \p_k \phi_l - 3V(\phi) g_{jk} + \mcf_{jk} - \ha \tr_g(\mcf) g_{jk}, \ \ j, k = 0, 1, 2, 3.
\eeq
The wave equation for $\phi$ can be written as
\beq
 -g^{\alpha\beta}  \p_\alpha\p_\beta \phi_l + \hat \Gamma^\alpha \p_\alpha \phi_l + V_l'(\phi) = \mcf^\Psi_l, \ \ l = 1, 2, \cdots, L.
\eeq
Now we let $u = (g - \hat g, \phi - \hat \phi)$ be the perturbed fields and we label $u = (u_l)_{l = 1}^{10+ L}.$ The system \eqref{einmat} can be written as 
\beqq\label{einmat1}
\begin{gathered}
-g^{jk}(x, u) \p_j\p_k u + H(x, u, \p u) = F, \ \ x\in M(\bt_0), \\
u = 0 \text{ in } M(\bt_0)\backslash J_{g}^+(\supp (F)).
\end{gathered}
\eeqq
See equation (26) of \cite{KLU1}. In particular, if $F$ is supported in a compact set and $\|F\|_{E^{m_0}(\bt_0)} < c_0$ sufficiently small with $m_0\geq 4$ even, there is a unique solution $u$ to \eqref{einmat1} and $\|u\|_{E^{m_0}(\bt_0)}\leq C_1\|F\|_{E^{m_0}(\bt_0)}$. 
If we take $\mcf = \eps f, \mcf^\Psi = \eps f^\Psi$ with $\eps > 0$ a small parameter, the solution $\phi$ depends on $\eps$ as well. Therefore, the tensor $T_{sour}(x, g, \phi)$ may contribute terms in the asymptotic expansion of $u$ in $\eps$.  Notice that this term only appears in the Einstein equations and the term has no derivative in $g$. 

Let $\dot g = \p_\eps g|_{\eps = 0}, \dot \phi = \p_\eps \phi|_{\eps = 0}$ be the linearized fields. Consider the linearized system \eqref{einmat1}, which can be written as
\beqq\label{perlin1}
\begin{gathered}
\square_{\hat g} \dot g + \mcb^1 (x, \dot g, \p \dot g, \dot \phi, \p \dot \phi) = \rho (\hat g, f),\text{ in } M(\bt_0),\\
\square_{\hat g} \dot \phi + \mcb^2 (x, \dot g,  \dot \phi) = f^\Psi ,\\ 
\dot g = 0, \ \ \dot \phi = 0, \text{ in } M(\bt_0)\backslash J_{\hat g}(\supp(\bar f)),
\end{gathered}
\eeqq
where $\bar f = (f, f^\Psi)$. As before, we can write this in matrix form. Let $v = (\dot g, \dot \phi)$ and we label as $v = (v_l)_{l = 1}^{10 + L}$. Then we can write \eqref{perlin1} as $\bold{P} v = (\rho(\hat g, f), f^\Psi)$. 
Observe that the linearized Einstein equation is not decoupled from the system. This would cause difficulties as singularities in the metric components $\dot g$ and scalar fields $\phi$ may mix together as they propagate. We shall fix the problem by working with a subset of the source.  Essentially the idea is to choose the source $f$ so that the linearized scalar fields $\dot \phi$ becomes smoother, in particular, their leading singularities vanish. This essentially put the problem into our framework. 

The causal inverse of $\bold{P}$ denoted by $\bold{Q}$ can be written as
\beq
\begin{pmatrix}
Q \bold{I}_{10} + \bold{U}_{11} & \bold{U}_{12}\\[3pt]
\bold{U}_{21} & Q\bold{I}_L + \bold{U}_{22}
\end{pmatrix}
\eeq
where $Q$ denotes the causal inverse of $\square_g$ and $\bold{U}_\bullet$ are some matrices. For $(x, \xi), (y, \eta)$ on the same null bi-characteristics, as in Prop.\ \ref{distor}, the principal symbol $\bold{R} \doteq \sigma(\bold{Q})(y, \eta, x, \xi)$ is linear invertible on $\mbr^{10 + L}$. Because of the gauge condition and conservation law, we know from the discussion in the end of Section \ref{sec-muls} that this is indeed a bijection on $\mbr^{6+L}$. We'd like to consider the symbol space of $\bar f$ at $(x, \xi)$ so that the range of $\bold{R}$ is $\bold{\mcx}_{x,\xi}$ consisting of vectors in $\mbr^{6+ L}$ whose last $L$ components are zero. In other words, we consider the subspace of the symbol of the sources so that the scalar field component are zero. We let 
\beq
\bold{\mcy}_{y, \eta} = \bold{R}^{-1}(\bold{\mcx}_{x, \xi})
\eeq
then the restriction $\bold{R}: \bold{\mcx}_{x, \xi}\rightarrow  \bold{\mcy}_{y, \eta}$ is bijection on the $6$-dimensional vector spaces. We remark that when restricted to $\mcx_{x, \xi}$, the map $\bold{R}$ is the same $R = \sigma(\bold{Q}_{\hat g})(y, \eta, x, \xi)$ as what we used for the Einstein equation before, because the $\mcb^2$ term  in the linearized field equation  \eqref{perlin1} does not contribute to the principal symbols. Also, the Einstein equations in \eqref{einmat1} has the same principal term and subprincipal term as the original Einstein system. In particular, the stress-energy term does not contribute to this term (it does not have derivative in the metric components). 

Since the leading order singularities in $\dot \phi$ vanishes, when we repeat the calculation in Prop.\ \ref{intert}, we get the same interaction term. Then Prop.\ \ref{symb} holds for system  \eqref{einmat1}. Finally, the rest of the argument in Section \ref{sec-wave} and \ref{sec-fermi} can be continued line by line to complete the proof of Theorem \ref{mainsca}.   
 \epf

\section{The Einstein-Maxwell equations}\label{sec-einmax}
Let's consider  the  electromagnetic fields. In the covariant formulation, the electromagnetic field can be described as a two form $F= (F_{jk})$. The stress-energy tensor for $F$ is 
\beq
 T^{em}_{jk} = F_j^\gamma F_{k\gamma} - \frac{1}{4} F^{\gamma \delta}F_{\gamma \delta} g_{jk} , \ \ j, k = 0, 1, 2, 3.
\eeq
This can be further generalized to the Yang-Mills equations.  The source for the electromagnetic field is the (four) electric current $\mcj$, a vector field on $M$. The Maxwell equations for $F$ on $(M, g)$ with source $\mcj$ are given by 
\beqq\label{maxwell}
dF = 0, \ \ \delta_g F = \mcj^\flat,
\eeqq
where the codifferential $\delta_g$ is the dual of the exterior differential $d$ with respect to  $g$ and $\mcj^\flat$ denotes the one form obtained from $\mcj$ by lowering the index using $g$. In simply connected domains, one can find a one form $\phi$ by Poincar\'e lemma such that $F = d\phi$, and $\phi$ is called the electromagnetic potential. The Maxwell equations can be written as
\beq
\delta_g d\phi = \mcj^\flat,
\eeq
while the first set of Maxwell equations is automatically satisfied. The source $\mcj$ is subject to the conservation law $\div_g \mcj = 0$. There is a gauge choice for $\phi$. In Lorentz gauge $\delta_g \phi = 0$, the Maxwell equations for $\phi$ are reduced to wave equations $\square_g \phi = \mcj^\flat$. We refer to \cite{Cb} and \cite[Section 18.5]{Tay} for more details of the background of this problem.   
 
We consider the following coupled Einstein-Maxwell system  
\beqq\label{einmax}
\begin{gathered}
\ein(g) = T^{em}(g, \phi) + \mcf \text{ in } M(\bt_0)\\
\delta_g d\phi = \mcj^\flat,\\
g = \hat g, \ \ \phi  = \hat \phi, \text{ in } M(\bt_0)\backslash J_g^+(\supp (\mcf)\cup \supp(J)).
\end{gathered}
\eeqq
Here, we denoted by $T_{em}(g, \phi)$ the stress-energy tensor associated with $F = d\phi$. The background fields $(\hat g, \hat \phi)$ are solutions to \eqref{einmat} with $\mcf = \mcj = 0$.  Also, the $g, \phi$ above shall be regarded as equivalent classes or representatives in a chosen gauge. As before, we let $T_{sour} =  T^{em}(g, \phi) + \mcf$. 
\begin{remark}
In \cite{KLU1}, the Einstein-Maxwell equations are of the following form
\beqq\label{eqeinmaxsour}
\begin{gathered}
\ein(g) = T^{em}(g, \phi) + T^{inter}\\
\delta_g d\phi = \mcj^\flat ,
\end{gathered}
\eeqq
where $T^{inter} = - \ha (\mcj^\mu\phi_\mu)g$ is the stress-energy tensor for the interaction term. This model is derived from the Lagrangian formulation in which only electromagnetic source perturbations are considered. So in \eqref{eqeinmaxsour}, the only source perturbation is $\mcj$, while for \eqref{einmax} studied here, the source perturbations are $\mcf$ and $\mcj$. For system \eqref{eqeinmaxsour}, the conservation law is reduced to $\div_g \mcj = 0$, while for \eqref{einmax}, we need the conservation laws $\div_g T_{sour} = 0$ and $\div_g \mcj = 0$ as discussed below.
\end{remark}

The local well-posedness of \eqref{einmax} can be established as before.  For $\delta > 0$, we let
\beqq\label{dataem}
\begin{gathered}
 \mcd^{em}(\delta; V) \doteq \{(g, \phi, \mcf, \mcj):  \mcf, \mcj \in C^{4}(M(\bt_0)), \ \ \|\mcf\|_{C^{4}(M(\bt_0))}, \|\mcj\|_{C^{4}(M(\bt_0))} < \delta, \\
  \supp (\mcf), \supp(\mcj) \subset V\\
\text{ such that there is a unique solution $(g, \phi)\in C^{4}(M(\bt_0))$ to \eqref{einmax}}\}.
\end{gathered}
\eeqq
We denote the observation in Fermi coordinates by $\mcd^{em}_F(\delta; V).$ 
The microlocal linearization stability conditions can be reformulated as follows. From \cite[Theorem 6.2, Chapter III]{Cb}, we know that 
\beq
(\div_g T^{em})^\beta = \nabla_\alpha T^{em, \alpha\beta} = \mcj^\la F^\beta_{\quad\la} = \mcj^\la g^{\beta\alpha}F_{\alpha\la}, \quad \beta = 0, 1, 2, 3, 
\eeq
is the Lorentz force where $\nabla$ denotes the covariant derivative. We write this in terms of the potential $\phi$, so the conservation law $\div_g T_{sour} = 0$ is 
\beq
\mcj^\la g^{\beta\alpha}(\p_\alpha \phi_\la - \p_\la\phi_\alpha) + (\div_g \mcf)^\beta = 0, \ \ \beta = 0, 1, 2, 3.
\eeq
Let $(\mcf, \mcj)  = \eps (f, J)$ depending on a small parameter $\eps$. We have the linearized conservation law
\beq
J^\la \hat g^{\beta\alpha}(\p_\alpha \hat \phi_\la - \p_\la \hat \phi_\alpha) + (\div_{\hat g} f)^\beta = 0, \ \ \beta = 0, 1, 2, 3.
 \eeq
We still need to consider the conservation law for the Maxwell equations $\div_g \mcj = 0$. The linearization gives $\div_{\hat g} J = 0$ which in local coordinates reads
\beq
\p_\alpha((-\det \hat g)^{\ha} J^\alpha) = 0.
\eeq
We see that the conservation laws holds where $\mcf = \mcj = 0.$ 
Now we state the microlocal linearization stability condition for the system \eqref{einmax}.

\begin{definition}\label{mlsem}
Let $Y$ be a two-dimensional submanifold in a space-like hypersurface of $M$ and $Y\subset V\subset M$. Let $(y, \eta)\in N^*Y$ with $\eta$ a light-like co-vector.  Let $A\in \mbr^{4}\times \mbr^4, B\in \mbr^4$ satisfy
\beq
\hat g^{\beta \alpha}(y)  \eta_\beta A_\alpha = 0, \quad \eta_\alpha B_\alpha = 0. 
\eeq
We say that the Einstein-Maxwell equation \eqref{einmax} satisfies the microlocal linearization stability condition if we can find a family $(g_\eps, \phi_\eps, \mcf_\eps, \mcj_\eps)\in \mcd^{em}(\delta; V)$ such that 
 \begin{enumerate}
\item  $ (\mcf_\eps, \mcj_\eps)$ are compactly supported in a conic neighborhood $W$ of $(y, \eta)$ in $N^*V$
\item   $ f= \p_\eps \mcf_\eps|_{\eps = 0} , J = \p_\eps \mci_\eps|_{\eps = 0}$ are in $I^{\mu+ 1}(N^*Y), \mu < -17$  such that 
\beq
\sigma(f)(y, \eta) = A, \ \ \sigma(J)(y, \eta) = B.
\eeq
\end{enumerate}
\end{definition}  
One can find examples for which the above condition holds in the following way. We consider $\mcf = T^{scalar}(g, \phi) +  \mcf^1$ in which $\mcf^1$ is a compactly supported $2$-tensor and $\phi \in \bold{B}^L$ are scalar fields satisfying the wave equation as before. In other words, we assume that the perturbation tensor $\mcf$ consists of a scalar field and another perturbation $\mcf^1$ and one can think of further coupling \eqref{einmax} with the scalar field equations. From the construction in \cite{KLU1} for the Einstein-scalar field equation, we see that the microlocal linearization holds for the coupled system, which further implies that the condition holds for \eqref{einmax} with the tensor $\mcf$. We remark that in this case, $\mcf$ is no longer compactly supported in $V$ (although $\mcf^1$ is). So the additional assumptions in Definition \ref{mlstab} do not hold. However, this  will not affect the analysis of singularities in Section \ref{sec-wgauge}, hence is not an issue here. For convenience, we still stated Definition \ref{mlsem} with the compact support assumption.\\

Now we state and prove the result for Einstein-Maxwell equations.
\begin{theorem}\label{mainem}
Let $(M, \hat g^{(i)}), i = 1, 2$ be two four-dimensional  simply connected smooth time oriented globally hyperbolic Lorentzian manifolds such that  $\hat g^{(i)}$ satisfy assumptions (A1).  
Let $\mu_{\hat g^{(i)}}(t), t\in [-1, 1]$ be time-like geodesics on $(M, \hat g^{(i)})$. Let $V$ be open relatively compact neighborhoods of $\mu_{\hat g^{(i)}}([-1, 1])$ and  $V\subset M^{(i)}(\bt_0)$ for some $\bt_0>0$.  Let $-1<s_-<s_+ < 1$, $p_\pm^{(i)} = \mu_{\hat g^{(i)}}(s_\pm)$ and $p^{(1)}_\pm = p^{(2)}_\pm$. 
For  $\delta > 0$, we consider solutions $(g^{(i)}, \phi^{(i)})$ to the Einstein equations
\beq 
\begin{gathered}
\ein(g^{(i)}) =   {T}^{em, (i)} +  \mcf^{(i)}, \text{ in } M^{(i)}(\bt_0),\\
\delta_{g^{(i)}}d\phi^{(i)} =  \mcj^{(i)},\\
g^{(i)} = \hat g^{(i)}, \ \ \phi^{(i)}  = \hat \phi^{(i)}, \text{ in } M^{(i)}(\bt_0)\backslash J_{g^{(i)}}^+(\supp (\mcf^{(i)})\cup \supp(\mcj^{(i)})),
\end{gathered}
\eeq 
which is well-posed on the data set $\mcd^{em, (i)}(\delta; V)$ defined as \eqref{dataem}. Suppose that 
 \begin{enumerate}
 \item  the microlocal linearization stability condition Definition \ref{mlsem} holds for $\mcd^{em, (i)}(\delta; V)$.
 \item the source-to-solution maps satisfy 
\[
 \mcd^{em, (1)}_F(\delta; V) = \mcd^{em, (2)}_F(\delta; V).
\]
\end{enumerate}
Then there is a diffeomorphism $\Psi: I_{\hat g^{(1)}}(p^{(1)}_-, p^{(1)}_+)\rightarrow I_{\hat g^{(2)}}(p^{(2)}_-, p^{(2)}_+)$ such that  $\Psi^*\hat g^{(2)} = \hat g^{(1)}$ in $I_{\hat g^{(1)}}(p^{(1)}_-, p^{(1)}_+)$. 
\end{theorem}

\bpf
Consider \eqref{einmax} in the wave and Lorentz gauge, that is 
\beqq\label{einmax1}
\begin{gathered}
\left\{\begin{array}{c}
\ric_{\widehat g}(g) = T_{em}(g, \phi) + \rho(g, \mcf)\\[3pt]
\square_g \phi =  \mcj^\flat \\
\end{array}\right. \text{ in } M(T_0),\\[2pt]
g = \widehat g, \ \ \phi  = \hat \phi,  \text{ in } M(T_0)\backslash J_g^+(\supp(\mcj)\cup\supp(\mcf)).
\end{gathered}
\eeqq
Note that $T_{em}(g, \phi)$ is a polynomial involving $g_{ij}, g^{ij}$ and first derivatives of $\phi$. Consider the perturbed fields
\beqq\label{perfield}
\vec w \doteq (u, \phi) = (g, \phi) - (\widehat g, \hat \phi). 
\eeqq
Using \eqref{ricre} and the expression of $\square_g \phi$, we see that equations \eqref{einmax1} in terms of $\vec w$ can be reduced to 
\beqq\label{pernonem}
\begin{gathered}
\left\{\begin{array}{c}
\square_g u    + \mca (x, \vec w, \p \vec w) =  G(x, \vec w, \mcf),\\[3pt]
\square_g \phi  =  \mcj,
\end{array}\right. \text{ in } M(T_0),\\[3pt]
\vec w = 0, \text{ in } M(\bold{t_0})\backslash J_g^+(\supp(\mcj)\cup\supp(\mcf)),
\end{gathered}
\eeqq
where $\mu, \nu, \beta = 0, 1, 2, 3$, $\mca$ is smooth functions of $\vec w, \p \vec w$ and $G(\bullet)$ is a smooth function in its arguments.  The well-posedness can be established as in Proposition \ref{stabest}. Let $v = (\dot g, \dot \phi)$ be the linearized fields and the linearization of \eqref{pernonem} is of the form 
\beq 
\begin{gathered}
\left\{\begin{array}{c}
\square_{\hat g} \dot g   + \mcb_1(x, v, \p v) =  G(x, v, \mcf),\\[3pt]
\square_{\hat g} \dot \phi  + \mcb_2 (x, v) =  \mcj.
\end{array}\right. 
\end{gathered}
\eeq 
We write this system in the matrix form as $\bold{P} v = F$ and let $\bold{Q}$ be the causal inverse. We need to discuss the symbol mapping. Recall that we work in the wave and Lorentz gauge. In addition to the harmonicity condition for $\dot g$, we have the Lorentz gauge condition for $\phi$: $\delta_g \phi = 0$ which is (in local coordinates)
\beq
g^{\alpha\la}\p_\la \phi_\alpha = \hat \Gamma^\mu \phi_\mu. 
\eeq
The linearization gives 
\beq
\hat g^{\alpha\la} \p_\la\dot\phi_\alpha + P(\hat g, \hat \phi, \dot g, \dot \phi) = 0
\eeq
for a smooth function $P$. Similar to the discussion in the end of Section \ref{sec-muls}, let $\mcy_{y, \eta}$ be the symbol space of the source $(f, J)$ at $(y, \eta)\in N^*Y$. Then $\mcy_{y, \eta}$ is a $9$ dimensional vector space because of the conservation laws. Let $\mcx_{x, \xi}$ be the symbol space of $v = (\dot g, \dot \phi)$ at $(x, \xi)$ on the bicharacteristics from $(y, \eta)$. We see that $\mcx_{x, \xi}$ is also $9$ dimensional because of the microlocal linearized gauge conditions. Thus 
\beq
\bold{Q}: \mcy_{y, \eta} \rightarrow \mcx_{x, \xi} \text{ is a linear bijection. }
\eeq
Now we  can follow the same arguments in Theorem \ref{mainsca} to finish the proof. 
 \epf

 \section*{Acknowledgments} 
The authors sincerely thank Matti Lassas for several helpful conversations, which influenced our understanding of the problem. Gunther Uhlmann is partially supported by NSF and a Si-Yuan Professorship at HKUST.



\begin{thebibliography}{99}
\bibitem{Bar} C. B\"ar, N. Ginoux, F. Pf\"affle. {\em Wave equations on Lorentzian manifolds and quantization.} ESI Lectures in Mathematics and Physics. European Mathematical Society (EMS), Z\"urich, 2007. 

\bibitem{Bea} M. Beals. {\em Propagation and interaction of singularities in nonlinear hyperbolic problems.} Progress in Nonlinear Differential Equations and their Applications, 3. Birkh\"auser Boston, 1989.

\bibitem{BEE} J. Beem, P. Ehrlich, K. Easley. {\em Global Lorentzian geometry.} Second edition. Monographs and Textbooks in Pure and Applied Mathematics, 202. Marcel Dekker, Inc., New York, 1996.

\bibitem{BS0} A. Bernal, M. S\'anchez. {\em On smooth Cauchy hypersurfaces and Geroch's splitting theorem.} Communications in Mathematical Physics 243.3 (2003): 461-470.

\bibitem{BS} A. Bernal, M. S\'anchez. {\em Globally hyperbolic spacetimes can be defined as 'causal' instead of 'strongly causal'.} Classical and Quantum Gravity 24.3 (2007): 745.

\bibitem{Cb} Y. Choquet-Bruhat. {\em General relativity and the Einstein equations.} Oxford Mathematical Monographs. Oxford University Press, Oxford, 2009.

\bibitem{CbD} Y. Choquet-Bruhat,  S. Deser. {\em On the stability of flat space.} Annals of Physics 81.1 (1973): 165-178.

\bibitem{DUV} M. de Hoop, G.  Uhlmann, A.  Vasy. {\em Diffraction from conormal singularities.} Annales Scientifiques de l'\'Ecole Normale Sup\'erieure, 4e serie, t.48, (2015): 351-408.

\bibitem{DUW} M. de Hoop, G. Uhlmann, Y. Wang. {\em Nonlinear responses from the interaction of two progressing waves at an interface.} Annales de l'Institut Henri Poincar\'e, Analyse Non Lin\'eaire (to appear).

\bibitem{DUW1} M. de Hoop, G. Uhlmann, Y. Wang. {\em  Nonlinear interaction of waves in elastodynamics and an inverse problem.} arXiv:1805.03811, (2018).

\bibitem{Du} J. J. Duistermaat. {\em Fourier integral operators.} Progress in Mathematics, 130. Birkh\"auser Boston, 1996.

\bibitem{FM} A. E. Fischer, J. E. Marsden. {\em The Einstein evolution equations as a first-order quasi-linear symmetric hyperbolic system, I.} Communications in Mathematical Physics 28.1 (1972): 1-38.

\bibitem {FM1} A. E. Fischer, J. E. Marsden. {\em Linearization stability of the Einstein equations.} Bulletin of the American Mathematical Society 79.5 (1973): 997-1003.

\bibitem{Fr} F. G. Friedlander. {\em The wave equation on a curved space-time.} Cambridge Monographs on Mathematical Physics, No. 2. Cambridge University Press, Cambridge-New York-Melbourne, 1975.

\bibitem{GB} J. Girbau, L. Bruna. {\em Stability by linearization of Einstein's field equation.} Progress in Mathematical Physics, 58. Birkh\"auser Verlag, Basel, 2010.

\bibitem{HE}  S. W. Hawking, G. F. R. Ellis. {\em The large scale structure of space-time.} Cambridge Monographs on Mathematical Physics, No. 1. Cambridge University Press, London-New York, 1973.

\bibitem{Ho} L. H\"ormander. {\em Lectures on nonlinear hyperbolic differential equations.} Math\'ematiques \& Applications (Berlin), 26. Springer-Verlag, Berlin, 1997.

\bibitem{Ho3} L. H\"ormander. {\em The analysis of linear partial differential operators III: Pseudo-differential operators.}  Reprint of the 1994 edition. Classics in Mathematics. Springer, Berlin, 2007. 

\bibitem{Ho4} L. H\"ormander. {\em The analysis of linear partial differential operators IV: Fourier integral operators.} Reprint of the 1994 edition. Classics in Mathematics. Springer-Verlag, Berlin, 2009. 

\bibitem{HKM} T. Hughes, T. Kato, J. Marsden. {\em Well-posed quasi-linear second-order hyperbolic systems with applications to nonlinear elastodynamics and general relativity.} Archive for Rational Mechanics and Analysis 63.3 (1977): 273-294.

\bibitem{GrU90} A. Greenleaf, G. Uhlmann. {\em Estimates for singular Radon transforms and pseudodifferential operators with singular symbols.} Journal of Functional Analysis 89.1 (1990): 202-232.

\bibitem{GrU93} A. Greenleaf, G. Uhlmann. {\em Recovering singularities of a potential from singularities of scattering data.} Communications in Mathematical Physics 157.3 (1993): 549-572.

\bibitem{Gri} J. B. Griffiths. {\em Colliding plane waves in general relativity.} Oxford Mathematical Monographs. Oxford Science Publications. The Clarendon Press, Oxford University Press, New York, 1991.

\bibitem{GU} V. Guillemin, G. Uhlmann. {\em Oscillatory integrals with singular symbols.} Duke Math. J 48.1 (1981): 251-267.

\bibitem{KLU} Y.  Kurylev, M. Lassas, G.  Uhlmann.  {\em Inverse problems for Lorentzian manifolds and non-linear hyperbolic equations.}  Inventiones Mathematicae (to appear). 

\bibitem{KLU1} Y. Kurylev, M. Lassas, G. Uhlmann. {\em Inverse problems in spacetime I: Inverse problems for Einstein equations-Extended preprint version.} arXiv:1405.4503 (2014).

\bibitem{KLU2} Y. Kurylev, M. Lassas, G. Uhlmann. {\em  Inverse problems in spacetime II: Reconstruction of a Lorentzian manifold from light observation sets.} arXiv:1405.3386 (2014).

\bibitem{LUW} M. Lassas, G. Uhlmann, Y. Wang. {\em Inverse problems for semilinear wave equations on Lorentzian manifolds.} Communications in Mathematical Physics (2018).

\bibitem{LUW1} M. Lassas, G. Uhlmann, Y. Wang {\em Determination of vacuum space-times from the Einstein-Maxwell equations} arXiv:1703.10704 (2017). 

\bibitem{Lee} J. M. Lee. {\em Riemannian manifolds: an introduction to curvature.} Graduate Texts in Mathematics, 176. Springer-Verlag, New York, 1997.

\bibitem{Ab} LIGO Project. {\em Observation of gravitational waves from a binary black hole merger.} Physical Review Letters 116.6 (2016): 061102.

\bibitem{LR} H. Lindblad, I. Rodnianski. {\em The weak null condition for Einstein's equations.} Comptes Rendus Mathematique 336.11 (2003): 901-906.


\bibitem{LuR} J. Luk, I. Rodnianski. {\em Local propagation of impulsive gravitational waves.} Communications on Pure and Applied Mathematics 68.4 (2015): 511-624. 

\bibitem{LuR1} J. Luk, I. Rodnianski. {\em Nonlinear interaction of impulsive gravitational waves for the vacuum Einstein equations.} arXiv:1301.1072  (2013).

\bibitem{MU} R. Melrose, G. Uhlmann. {\em Lagrangian intersection and the Cauchy problem.} Communications on Pure and Applied Mathematics 32.4 (1979): 483-519.

\bibitem{MTW} C. W. Misner, K. S. Thorne, J. A. Wheeler. {\em Gravitation.} W. H. Freeman and Co., San Francisco, Calif., 1973.

\bibitem{Mo} V. Moncrief. {\em Space-time symmetries and linearization stability of the Einstein equations. II.} Journal of Mathematical Physics 17.10 (1976): 1893-1902.

\bibitem{NO} K. Nomizu,  H. Ozeki. {\em The existence of complete Riemannian metrics.} Proceedings of the American Mathematical Society 12.6 (1961): 889-891.

\bibitem{Tay} M. E. Taylor. {\em Partial differential equations III. Nonlinear equations. Second edition.} Applied Mathematical Sciences, 117. Springer, New York, 2011. 

\bibitem{WZ}  Y. Wang, T. Zhou. {\em Inverse problems for quadratic derivative nonlinear wave equations.}  Communications in Partial Differential Equations (to appear).
\end{thebibliography}
\end{document}